\documentclass{amsbook}
\usepackage[cp1251]{inputenc}
\usepackage[english,russian]{babel}
\pdfoutput=1
\usepackage{chapterbib}

\usepackage{bbm}
\usepackage[a5paper, total={113mm,178mm}, left=2cm, includehead]{geometry}

\usepackage{amsfonts,amssymb,amsmath,amstext,mathrsfs}
\usepackage{latexsym}

\usepackage{euscript}
\usepackage{eufrak}

\usepackage{graphicx}

\makeatletter

\renewcommand{\chapter}{\@startsection{chapter}{1}{0pt}{-3.25ex plus -1ex minus-.2ex}{1.5ex plus .2ex}{\bfseries}}

\renewcommand{\section}{\@startsection{section}{1}{0pt}{-3.25ex plus -1ex minus-.2ex}{1.5ex plus .2ex}{\bfseries}}

\renewcommand{\subsection}{\@startsection{subsection}{2}{0pt}{-3.25ex plus -1ex minus-.2ex}{-1em}{\bfseries}}

\renewcommand{\@evenhead}%
{\raisebox{0pt}[\headheight][0pt]{%
\vbox{\hbox to\textwidth{%
\thepage\hfil\strut\leftmark}\hrule}}%
}
\renewcommand{\@oddhead}%
{\raisebox{0pt}[\headheight][0pt]{%
\vbox{\hbox to\textwidth{%
\rightmark\hfil\strut\thepage}\hrule}}%
}
\renewcommand{\@evenfoot}{}
\renewcommand{\@oddfoot}{}

\makeatother

\def\R{\mathbb{R}}
\def\N{\mathbb{N}}
\def\Z{\mathbb{Z}}

\newcommand{\cA}{{\mathcal A}}
\newcommand{\cB}{{\mathcal B}}
\newcommand{\cC}{{\mathcal C}}
\newcommand{\cD}{{\mathcal D}}
\newcommand{\cE}{{\mathcal E}}

\newcommand{\cI}{{\mathcal I}}

\newcommand{\cL}{{\mathcal L}}
\newcommand{\cM}{{\mathcal M}}
\newcommand{\cN}{{\mathcal N}}
\newcommand{\cO}{{\mathcal 0}}

\newcommand{\cQ}{{\mathcal Q}}

\newcommand{\cS}{{\mathcal S}}

\newcommand{\cW}{{\mathcal W}}
\newcommand{\cX}{{\mathcal X}}

\def\0{\boldsymbol{0}}
\def\1{\boldsymbol{1}}

\newcommand {\Rmaxm}	{\R_{\max,\text{m}}}

\newcommand{\Rmaxmn}	{\Rmaxm^n}

\renewcommand{\thesection}{\arabic{section}}

\newtheorem{theorem}{Theorem}[section]
\newtheorem{lemma}[theorem]{Lemma}
\newtheorem{proposition}[theorem]{Proposition}
\newtheorem{corollary}[theorem]{Corollary}

\theoremstyle{definition}
\newtheorem{definition}[theorem]{Definition}

\theoremstyle{remark}
\newtheorem{remark}[theorem]{Remark}

\numberwithin{equation}{section}

\theoremstyle{theorem}
\newtheorem{teorema}          {Теорема}[section]

\newtheorem{predl}    {Предложение}[section]

\theoremstyle{definition}

\newtheorem{opred}	      {Определение}[section]

\def\Bbb{\mathbb}
\def\goth{\mathfrak}

\newcommand{\const}{\mathop{\rm const}\nolimits}
\newcommand{\gradd}{\mathop{\rm grad}\nolimits}

\newcommand{\arxiv}[1]{arXiv:#1}

\newcounter{assume}[section]
\def\theassume{(A\arabic{assume})}

\def\of#1,#2{(#1,#2)}

\def\tchk{\cdot}
\def\msn{\medskip\noindent}

\newcommand{\x}{\mathbf x}

\newcommand\closure{\operatorname{cl}}

\newcommand\hil[3]{\text{Hil}_{#3}(#1,#2)}
\newcommand\funk[3]{\text{Funk}_{#3}(#1,#2)}

\newcommand\rev[3]{\text{rev}_{#3}(#1,#2)}

\newcommand\dotprod[2]{\langle{#1}|{#2}\rangle}

\newcommand\poynta{x}
\newcommand\poyntb{y}

\newcommand\dist{d}
\newcommand\dihed{A_k}
\newcommand\prodd{\operatorname{prod}}

\newcommand\gens{S}
\newcommand\dualgens{\tilde S}

\newcommand\plusclass{+\hat\infty}
\newcommand\minusclass{-\hat\infty}

\newcommand{\union}{\cup}

\newcommand\indicator{I}
\newcommand\dualspace{V^*}
\newcommand\dualball{{B^\circ}}

\newcommand\hnorm{f}

\def\got{\mathfrak}

\DeclareMathOperator{\Br}{\mathrm  Br}

\newcommand{\wt}{\widetilde}

\newcommand{\bR}{{\bf R}}

\newcommand{\wh}{\widehat}

\newcommand{\pa}{\partial}

\newcommand{\la}{\lambda}

\newcommand{\grad}{\nabla}

\newcommand{\liiR} {I\!\!R}

\newcommand{\alphaoo}{{\alpha^{1,1}}}
\newcommand{\alphaot}{{\alpha^{1,2}}}

\newcommand{\alphatt}{{\alpha^{2,2}}}

\newcommand{\Aprime}{{A^\prime}}

\newcommand{\betasubto} {{\beta_{t_1}}}

\newcommand{\betasubtt} {{\beta_{t_2}}}

\newcommand{\betasubtoptt} {{\beta_{t_1+t_2}}}
\newcommand{\Bto} {{B_{t_1}}}
\newcommand{\Btt} {{B_{t_2}}}
\newcommand{\Btoptt} {{B_{t_1+t_2}}}
\newcommand{\capow}[1]{^{\circledast #1}}

\newcommand{\cIbarhatmk}{{\widehat{\overline \cI}^{\raise -0.3em \hbox{\scriptsize$m$}}_k}}
\newcommand{\cIbarhatmf}{{\widehat{\overline \cI}^{\raise -0.3em \hbox{\scriptsize$m$}}_f}}

\newcommand{\cIbartildemtk}{{\widetilde{\overline \cI}^{\raise -0.3em \hbox{\scriptsize$m,t$}}_k}}
\newcommand{\cIbartildemtf}{{\widetilde{\overline \cI}^{\raise -0.3em \hbox{\scriptsize$m,t$}}_f}}

\newcommand{\cWp}{\cW^+}

\newcommand{\demi} {{\textstyle\frac{1}{2}}}
\newcommand{\dualop}{{\cD_\psi}}
\newcommand{\dualopinv}{{\cD_\psi^{-1}}}

\newcommand{\etaoo}{{\eta^{1,1}}}
\newcommand{\etaot}{{\eta^{1,2}}}

\newcommand{\etatt}{{\eta^{2,2}}}

\newcommand{\gammaoo}{{\gamma^{1,1}}}
\newcommand{\gammaot}{{\gamma^{1,2}}}
\newcommand{\gammato}{{\gamma^{2,1}}}
\newcommand{\gammatt}{{\gamma^{2,2}}}
\newcommand{\gamovertwo} {{\frac{\gamma^2}{2}}}

\newcommand{\gammovertwo} \gamovertwo
\newcommand{\gammuovertwo} \gamovertwo
\newcommand{\gammuepstovertwo} \gamovertwo
\newcommand{\gammubarepstovertwo} \gamovertwo

\newcommand{\intmprn} {\int^\oplus_{\liiR^n}}

\newcommand{\lambdahat}{{\hat\lambda}}

\newcommand{\matdualop}{{D_\psi}}
\newcommand{\matdualopinv}{{D_\psi^{-1}}}

\newcommand{\parenprime}{{)^\prime}}

\newcommand{\Ptilde}{{\widetilde P}}
\newcommand{\psicdotz} {\psi(\cdot,z)}
\newcommand{\psixcdot} {\psi(x,\cdot)}
\newcommand{\qforall} {\qquad\forall\,}

\newcommand{\liiRn} {\liiR^n}

\newcommand{\solt} {{S_t}}

\newcommand{\soltau} {{S_\tau}}

\newcommand{\Tbar} {{\overline T}}

\newcommand{\That} {{\widehat T}}
\newcommand{\Ttilde} {{\widetilde T}}

\newcommand{\Uinv}{{U^{-1}}}

\newcommand{\Vttdel} {\widetilde{\widetilde{V}}^{\raise -0.3em \hbox{\scriptsize$\delta$}} }
\newcommand{\Vz}{V^z}

\newcommand{\Whatz} {\widehat W^z}
\newcommand{\xbar}{\overline x}

\newcommand{\xprime}{{x^\prime}}

\newcommand{\zbar}{\overline z}

\newcommand{\zprime}{{z^\prime}}

\newcommand{\midalign} {&&\hskip -0.6em}

\newcommand{\maxp} {{max-plus }}

\newcommand{\beasnum}{\begin{eqnarray}}
\newcommand{\eeasnum}{\end{eqnarray}}
\newcommand{\beas}{\begin{eqnarray*}}
\newcommand{\eeas}{\end{eqnarray*}}

\newcommand{\er}[1]{\hbox{(\ref{#1})}}

\newcommand{\lineika}[2]{\hbox to\textwidth{\normalsize {\em #1} \dotfill  #2}}

\newcommand{\statya}[3]{\noindent\parbox[l]{\textwidth} 
{\parbox[l]{0.9\textwidth}{\small\scshape #1}%
\par\vspace{2mm}%
\noindent\lineika{#2}{#3}}%
\vskip0.4cm}

\newcommand{\russtatya}[3]{\noindent\parbox[l]{\textwidth}
{\parbox[l]{0.9\textwidth}{{\small\scshape #1}
{\bfseries (in Russian)}}%
\par\vspace{2mm}%
\noindent\lineika{#2}{#3}}%
\vskip0.4cm}

\newcommand{\nestatya}[2]{\noindent\hbox to\textwidth{\normalsize #1\dotfill  #2}}

\newcommand{\avtor}[4]{\noindent%
\parbox[l]{\textwidth}{\small {\bfseries #1,}\\ {\upshape #2,}\ {\scshape #3.}\ {\ttfamily #4}
{\tolerance=200

}
}%
\vskip0.3cm}
\newcommand{\rusavtor}[5]{\noindent%
\parbox[l]{\textwidth}{\small {\bfseries #1}\ {\slshape (#2),}\\ {\upshape #3,}\ {\scshape #4.}\ {\ttfamily #5}}
{\tolerance=200

}
\vskip0.3cm}

\newcommand{\nachalo}[2]{\par%
\pagebreak[2]\vspace{1cm plus 3mm minus 0.5mm}%
\begin{center}
{\large\rmfamily\bfseries\upshape #2\par}%
\nopagebreak%
\vspace{2mm plus 1mm}
\nopagebreak%
{\large\itshape #1\par}%
\end{center}
\nopagebreak\vspace{2mm plus 1mm}
\nopagebreak}

\newcommand{\nachaloe}[3]{\par%
\pagebreak[2]\vspace{1cm plus 3mm minus 0.5mm}%
\begin{center}%
{\large\rmfamily\bfseries\upshape #2\footnote{#3}\par}%
\nopagebreak%
\vspace{2mm plus 1mm}
\nopagebreak%
{\large\itshape #1\par}%
\end{center}%
\nopagebreak%
\vspace{2mm plus 1mm}%
\nopagebreak}

\def\fizfak{M.V. Lomonosov Moscow State University, 
Fac\-ulty of Phy\-sics, Chair of Quan\-tum Sta\-tis\-tics
and Field The\-ory, 119992, Le\-nin\-skie Gory, \mbox{Moscow}}
\def\INRIA{INRIA, Dom\-aine de Volu\-ceau,
B.P.~105, 78153 Le Ches\-nay Cedex}
\def\nezavisimy{Ind\-ep\-en\-dent 
Uni\-ver\-sity of Mos\-cow, 119002, B. Vla\-sy\-ev\-sky per. 11, \mbox{Moscow}}
\def\novisadg{Uni\-ver\-sity of Novi Sad, Fac\-ulty of Tech\-ni\-cal 
Sci\-en\-ces,\\ Trg Do\-si\-teja Ob\-ra\-do\-vi\-ca 6, 21000 Novi Sad}
\def\novisadps{Uni\-ver\-sity of Novi Sad, Fac\-ulty of Sci\-en\-ces and Math\-ema\-tics,\\
Trg Do\-si\-teja Ob\-ra\-do\-vica 4, 21000 Novi Sad}
\def\mehmat{M.V. Lo\-mo\-no\-sov Mos\-cow State Uni\-ver\-sity, Dep\-art\-ment of Mech\-an\-ics and 
Math\-ema\-tics, Fac\-ulty of Higher Alg\-ebra, 119992 Le\-nin\-skie Go\-ry, \mbox{Moscow}}
\def\IRMA{IRMA, Uni\-ver\-site Louis Pas\-teur, 7~rue Rene Des\-car\-tes, 67084 Stras\-bourg
Ce\-dex}
\def\delfin{Del\-fin-In\-for\-ma\-ti\-ka Co., ul.~Tkac\-kaya 1, \mbox{Moscow}}
\def\CEMI{CEMI RAS, Na\-hi\-mov\-ski prosp.~47, 117418, \mbox{Moscow}}
\def\krasped{Kra\-sno\-yarsk State Ped\-ag\-ogi\-cal 
Uni\-ver\-sity, 660049, ul.~ A.~ Le\-be\-de\-voy~89,\\ \mbox{Krasnoyarsk}}
\def\krasfed{Siberian Fed\-eral Uni\-ver\-sity, Dep\-art\-ment of 
Math\-ema\-tics and Inf\-orma\-tics,
660041, pr.~Svo\-bod\-nyi~79, \mbox{Krasnoyarsk}}
\def\tambov{G.R. Der\-zha\-vin Tam\-bov State Uni\-ver\-sity, 392000,
ul.~In\-ter\-na\-cio\-nal\-naya 33, Tambov}
\def\romanian{Ins\-ti\-tute of Math\-ema\-tics of the Roma\-nian Aca\-demy,
P.O. ~Box 1-764,\\ 010702 Bucharest}
\def\sengalese{Dep\-art\-ment of Math\-ema\-tics \& Com\-pu\-ter Sci\-ence, UCAD, Dakar}
\begin{document}
\selectlanguage{english}
\thispagestyle{empty}
\vbox to \textheight{\centering
{\scshape Independent University of Moscow \\
French--Russian Laboratory ``J.-V.~Poncelet''}
\vskip6pt
\hrule

\vss 

{\large International Workshop}

\vss

{\Large
{\fontseries{b}\selectfont IDEMPOTENT AND TROPICAL MATHEMATICS AND~PROBLEMS
OF~MATHEMATICAL PHYSICS}}
\indent

\vss

{\large\slshape G.L.~Litvinov, V.P.~Maslov, S.N.~Sergeev (Eds.)}
%\begin{tabular}{rl}
%{\large Editors:} & {\large\slshape G.L.~Litvinov, V.P.~Maslov, S.N.~Sergeev}
%\end{tabular}

\vss

\begin{tabular}{rl}
{\large Organizing committee:} & {\large\slshape G.L.~Litvinov, V.P.~Maslov,}\\
& {\large\slshape S.N.~Sergeev, A.N.~Sobolevski\u{\i}}\\
&\\
{\large Web-site:} & {\verb"http://www.mccme.ru/tropical07"}\\
&\\
{\large E-mail:} & {\verb"tropical07@gmail.com"}
\end{tabular}

\vss

Moscow, August 25--30, 2007 

\vss

{\large\bfseries Volume II}

\vss

Moscow, 2007
}
\newpage
\setcounter{page}{2}
\thispagestyle{empty}
\begin{flushleft}
\parbox[l]{0.8\textwidth}{
\indent {\bfseries Litvinov G.L., Maslov V.P., Sergeev S.N. (Eds.)}\\
\indent Idempotent and tropical mathematics and problems of mathematical physics 
(Vol. II) -- M.: 2007 -- 116 pages
\vskip12pt
\indent {\slshape This volume contains the proceedings of an International Workshop on Idempotent and Tropical Mathematics
and Problems of Mathematical Physics, held at the Independent University of Moscow, Russia, on August 25-30, 
2007.}
}
\end{flushleft}
%Organizing committee: G.L.~Litvinov, V.P.~Maslov, S.N.~Sergeev, and A.N.~Sobolevski\u{\i}}
\vskip 8cm
\begin{flushright}
\copyright\ 2007 by the Independent University of Moscow. All rights reserved. 
\end{flushright}
\newpage
\section*{CONTENTS}
\vskip0.5cm

%\nestatya{{\scshape Preface} {\bfseries(in Russian)}}{\pageref{preface-rus}}
%\vskip0.3cm

\statya{Ultrasecond quantization of
a classical version of \mbox{superfluidity in nanotubes}}{Victor P.~Maslov}{\pageref{mas-abs}}

\statya{Policy iteration and max-plus finite element method}{David McCaffrey}{\pageref{mcc-abs}}

\statya{Using max-plus convolution to obtain fundamental solutions for differential equations with quadratic nonlinearities}
{William M.~McEneaney}{\pageref{mce-abs}}
\thispagestyle{empty}
\statya{Polynomial quantization on \mbox{para-hermitian}
\mbox{symmetric} \mbox{spaces} from the \mbox{viewpoint} of
\mbox{overgroups: an example}}{Vladimir F.~Molchanov}{\pageref{mol-abs}}

\statya{The structure of max-plus hyperplanes}{V.~Nitica and I.~Singer}%
{\pageref{sin-abs}}

\statya{Image processing based on a partial differential equation
satisfying the pseudo-linear superposition principle}{E.~Pap and M.~\v{S}trboja}%
{\pageref{str-abs}}

\statya{Tropical analysis on plurisubharmonic singularities}%
{Alexander Rashkovskii}{\pageref{ras-abs}}

\statya{Minimal elements and 
cellular closures over the max-plus semiring}{Serge{\u{\i}} Sergeev}%
{\pageref{ser-abs}}
\thispagestyle{empty}

\statya{Semiclassical quantization of field theories}{Oleg Yu.~Shvedov}%
{\pageref{shv-abs}}

\statya{Convex analysis, transportation and 
reconstruction of peculiar velocities of galaxies}{Andre\u{\i} Sobolevski\u{\i}}%
{\pageref{sob-abs}}

\statya{The Weyl algebra and quantization of fields}{Alexander V.~Stoyanovsky}%
{\pageref{sto-abs}}

\statya{Polynomial quantization on para-hermitian spaces with pseudo-orthogonal group of translations}%
{Svetlana V.~Tsykina}{\pageref{tsy-abs}}

\statya{The horofunction boundary}{Cormac Walsh}{\pageref{wal-abs}}

%\statya{Квантование как приближенное
%описание некоторого диффузионного процесса}{Е.М. Бениаминов}%
%{\pageref{ben-rus-abs}}

\russtatya{Quantization as approximate description of a diffusion process}%
{Evgeny M.~Beniaminov}{\pageref{ben-rus-abs}}
\thispagestyle{empty}

\russtatya{Idempotent systems 
of nonlinear equations 
and computational problems 
arising in electroenergetic networks}{A.M.~Gel'fand and B.Kh.~Kirshteyn}%
{\pageref{kir-rus-abs}}

\russtatya{Classical and nonarchimedean amoebas in extensions of fields}%
{Oksana V.~Znamenskaya}{\pageref{zna-rus-abs}}

%\statya{Идемпотентные системы нелинейных уравнений и задачи расчета 
%электроэнергетических сетей}{А.М.~Гельфанд и Б.Х.~Кирштейн}%
%{\pageref{kir-rus-abs}}

\russtatya{Generalization of ultra second quantization for fermions
at non-zero temperature}{G.V.~Koval and V.P.~Maslov}%
{\pageref{kov-rus-abs}}

\russtatya{Contact classification of Monge-Amp\`ere equations\\}{Alexey G.~Kushner}%
{\pageref{kus-rus-abs}}
%\statya{Обобщение ультравторичного квантования
%для фермионов при ненулевой температуре}{Г.В.~Коваль и В.П.~Маслов}%
%{\pageref{kov-rus-abs}}

%\russtatya{Ultrasecond quantization of a classical 
%version of superfluidity in nanotubes}{Victor P.~Maslov}{\pageref{mas-rus-abs}}

%\statya{Ультравторичное квантование и 
%классическая версия сверхтекучести в нанотрубках}{В.П.~Маслов}%
%{\pageref{mas-rus-abs}}

\russtatya{On amoeba of discriminant of algebraic equation}
{Evgeny N.~Mikhalkin}{\pageref{mik-rus-abs}}

%\statya{Об амебе дискриминанта 
%алгебраического уравнения}{Е.Н.~Михалкин}{\pageref{mik-rus-abs}}
\thispagestyle{empty}

\russtatya{Universal algorithms solving discrete 
Bellman systems of equations over semirings 
(a computer demonstration)\\}{A.V. ~Chourkin and S.N.~Sergeev}%
{\pageref{chu-rus-abs}}

%\statya{Программа для
%демонстрации универсальных алгоритмов решения уравнения
%Беллмана в различных полукольцах}%
%{С.Н.~Сергеев и А.В.~Чуркин}{\pageref{chu-rus-abs}}

\russtatya{The curves in $\mathbb{C}^2$ whose amoebas
determine the fundamental group of complement}{Roman Ulvert}
{\pageref{ulv-rus-abs}}

\nestatya{\scshape List of participants and authors}{\pageref{list-of-pars}}

%\statya{Кривые в $\mathbb{C}^2$, амебы которых
%определяют фундаментальную группу дополнения}{Роман Ульверт}%
%{\pageref{ulv-rus-abs}}
%\selectlanguage{russian}
%\include{preface-rus}
%\selectlanguage{english}
\newpage
\renewcommand{\theequation}{\arabic{equation}}
\setcounter{footnote}{0}
\setcounter{equation}{0}
%begin{center}
%\bf\large Ultrasecond quantization of}\\
%\bf\large a classical version of superfluidity in nanotubes}%
%footnote{Supported by the joint RFBR-CNRS grant 05-01-02807 and by the 
%FBR grant 05-01-00824}\\
%vskip0.5cm
%\it\large V.P.~Maslov}
%end{center}
%null
\nachaloe{Victor P.~Maslov}{Ultrasecond quantization of
a \mbox{classical} \mbox{version} of \mbox{superfluidity in nanotubes}}%
{The work has been supported 
by the joint RFBR/CNRS grant 05-01-02807 and by the 
RFBR grant 05-01-00824.}
\label{mas-abs}
\markboth{Victor P.~Maslov}{Ultrasecond quantization}

\textbf{1.}  In order to distinguish the classical theory in its modern understanding
from the quantum theory, it is necessary to modify (somewhat) the ideology habitual to physicists, for whom the classical theory is simply the whole body of physics as 
it existed in the 19th century 
before the appearance of quantum theory. Actually, the correct meaning is that the classical theory is  the limit of the quantum one as $h\to 0$.

Thus, Feynman correctly understood that spin is a notion of classical mechanics. Indeed,
it is obtained via a rigorous passage from quantum mechanics to classical mechanics
\cite{TeorVoz}. In a similar same way, the polarization of light does not disappear when the frequency is increased, and is therefore a property of geometric rather than wave optics, contrary to the generally accepted belief, which arose because the polarization of light was discovered as the result of the appearance of wave optics.

Consider a ``Lifshits hole'', i.e. a one-dimensional Schr\"odinger equation with potential symmetric with respect to the origin of coordinates with two troughs.
Its eigenfunctions are symmetric or antisymmetric with respect to the origin. As
$h\to 0$ this symmetry remains, and since the square of the modulus of the eigenfunction corresponds to the probability of the particle to remain in the troughs, it follows that in the limit as $h\to 0$, i.e., in the ``classical theory", for energies less than those required to pass over the barrier, the particle is simultaneously located in two troughs, although a classical particle cannot pass through the barrier. Nevertheless, this simple example shows how the ideology of the ``classical theory'' must be modified.

To understand this paradox, one must take into consideration the fact that the symmetry must be very precise, up to ``atomic precision'', and that stationary state means  a state that arises in the limit for ``infinitely long''  time.

When we deal with nanotubes whose width is characterized by ``atomic'' or ``quantum'' dimensions, then new unexpected effects occur in the ``classical'' theory. Thus, already in 1958  \cite{DAN_58}, I discovered a strange effect of the standing longitudinal wave type in a slightly bent infinite narrow tube, for the case in which its radius is the same everywhere with atomic precision. It was  was impossible at the time to implement this effect in practice, which would have allowed to obtain a unimode laser, despite   A.M.Prokhorov's great interest in the effect.

\textbf{2.} Now let us discuss the notion known as ``collective oscillations''   in classical physics and as ``quasiparticles''  in quantum physics. In classical physics, it is described by the Vlasov equation for selfcompatible (or mean) fields, in quantum physics, by the Hartrey (or the Hartrey-Fock) equation.

(1) Variational equations depend on where (i.e., near what solutions of the original equation) we consider the variations. For example, in \cite{QuasiPart5, QuasiPart6, QuasiPart7} we considered variations near a microcanonical distribution in an ergodic construction, while in \cite{QuasiPart1, QuasiPart2, QuasiPart3, QuasiPart4} this was done near a nanocanonical distribution concentrated on an invariant manifold of lesser dimension, i.e., not on a manifold of constant energy but, for example, on a Lagrangian manifold of dimension coinciding with that of the configuration space.

(2) Let us note the following crucial circumstance. The solution of the variational equation for the Vlasov equation   \textit{does not coincide} with the classical limit for variational equations for the mean field equations in  quantum theory.

Consider the mean field equation in the form

\begin{equation}
\label{Hart}
\begin{split}
ih\frac{\pa}{\pa t}\varphi^t (x)&= \bigg(-\frac{h^2}{2m}\Delta
+W_t(x)\bigg) \varphi^t(x), \\
W_t(x)&=U(x)+\int V(x,y)|\varphi^t(y)|^2dy,
\end{split}
\end{equation}
with the initial condition
 $\varphi|_{t=0}=\varphi_0$, where $\varphi_0$ belongs to 
$W_2^\infty (\bR^\nu)$ and satisfies
$\int dx|\varphi_0(x)|^2=1$.

In order to obtain asymptotics of the complex germ type \cite{VKB}
one must write out the system consisting of the Hartrey equation and its dual,
then consider the corresponding variational equation, and, finally, replace
the variations  $\delta\varphi$ and $\delta\varphi^*$ 
by the independent functions $F$ and $G$. 
For the functions $F$ and $G$, we obtain the following system of equations:

\begin{equation}
\label{Shv1}
\begin{split}
& i\frac{\pa F^t(x)}{\pa t} = \int dy
\bigg(\frac{\delta^2H}{\delta\varphi^*(x)\delta\varphi(y)}F^t(y)+
\frac{\delta^2H}{\delta\varphi^*(x)\delta^*\varphi(y)}G^t(y)\bigg);\\
& -i\frac{\pa G^t(x)}{\pa t} = \int dy
\bigg(\frac{\delta^2H}{\delta\varphi(x)\delta\varphi(y)}F^t(y)+
\frac{\delta^2H}{\delta\varphi(x)\delta^*\varphi(y)}G^t(y)\bigg).
\end{split}
\end{equation}

The classical equations are obtained from the quantum ones, roughly speaking,
by means of a substitution of the form $\varphi= \chi e^{\frac ih S}$
(the VKB method), $\varphi^*= \chi^* e^{\frac ih S^*}, where \ S=S^*, \
\chi=\chi(x,t)\in C^\infty, \ S=S(x,t)\in C^\infty$.

To obtain the variational equations, it is natural to take the variation not only of the
limit equation for  $\chi$ and  $\chi^*$, but also for the functions $S$ and
$S^*$. This yields a new important term of the equation for collective oscillations.

Let us describe this fact for the simplest example, which was studied in
N.N.Bogolyubov's famous paper concerning ``weakly ideal Bose gas" \cite{Bogol}.

Suppose $U=0$ in equation~(\ref{Hart}) in a three-dimensional cubical
box of edge $L$, the wave functions satisfying the periodicity condition (i.e., the
problem being defined on the 3-torus with generators of lengths  $L, L, L$).
Then the function

\begin{equation}\label{bog1}
  \varphi(x)=L^{-3/2} e^{i/h(px-\Omega t)},
\end{equation}
where $p=2\pi n/L$, $n$ is an integer,  satisfies the equation
(\ref{Hart}) for
\begin{equation}\label{bog2}
  \Omega=\frac{p^2}{2m}+L^{-3}\int dx V(x).
\end{equation}

For $\lambda=2\pi {n}/L$, $n$ a non\-zero integer, consider the functions $F^{(\lambda)}(x)$ and $G^{(\lambda)}(x)$
given by

\begin{equation}
\label{bog3}
\begin{split}
& F^{(\lambda)t}(x)= L^{-3/2}\rho_\lambda e^{\frac
ih|(p+\lambda)x +(\beta-\Omega)t|}, \\
& G^{(\lambda)t}(x)= L^{-3/2}\sigma_\lambda e^{\frac
ih|(-p+\lambda)x +(\beta+\Omega)t|};
\end{split}
\end{equation}
here
\begin{equation}
\label{bog4}
\begin{split}
&-\beta_\lambda\rho_\lambda=\bigg(\frac{(p+\lambda)^2}{2m}-\frac{p^2}{2m}
+\wt{V}_\lambda\bigg)\rho_\lambda+V_\lambda\sigma_\lambda, \\
&\beta_\lambda\rho_\lambda=\bigg(\frac{(p-\lambda)^2}{2m}-\frac{p^2}{2m}
+\wt{V}_\lambda\bigg)\sigma_\lambda+V_\lambda\rho_\lambda,\\
&|\sigma_\lambda|^2 -|\rho_\lambda|^2=1, \qquad
\wt{V}_\lambda=L^{-3}\int dx V(x)e^{\frac ih \lambda x}. 
\end{split}
\end{equation}

From the system (\ref{bog4}), we find
\begin{equation}\label{bog5}
\beta_\lambda = -p \lambda+
\sqrt{\bigg(\frac{\lambda^2}{2m}+\wt{V}_\lambda\bigg)^2-\wt{V}_\lambda^2}.
\end{equation}

In this example $u=e^{\frac ih s(x,t)}, \ u^*=e^{-\frac
{s(x,t)}{h}}$, where $s(x,t)= px+\beta t$, while the variation of the action for the vector
 $\big({\delta u}, {\delta u^*}\big)$ equals $\lambda x\pm\Omega t$.

Under a more accurate passage to the limit, we obtain $$\wt{V}_\lambda \to V_0
=L^{-3} \int dx V(x)$$.

Thus, in the classical limit, we have obtained the famous Bogolyubov relation
 (\ref{bog5}).  In the case under consideration
$u(x)=0$ and, as in the linear Schr\"odinger equation, the exact solution coincides with the quasiclassical one. In the paper  \cite{QuasiPart4}, the case
$u(x)\neq 0$ is investigated, and it turns out that the relation similar to
 (\ref{bog5}) is the classical limit as $h\to 0$ of the variational equation in this general case. 
The curve showing the dependence of $\beta_\lambda$ on $\lambda$
is known as the {\em Landau curve} and determines the superfluid state. The value
 $\lambda_{\text{cr}}$ for which superfluidity disappears is called the {\em Landau critical level}. Bogolyubov explains 
the superfluidity phenomenon in the following terms:
``the `degenerate condensate' can move 
without friction relatively to elementary 
perturbations with any sufficiently small velocity''  \cite[p. ~210]{QuasiPart4}.

However, there is no Bose-Einstein condensate 
whatever in these mathematical considerations, 
it is just that the spectrum defined for  $\lambda <\lambda_{\text{cr}}$  
is a positive spectrum of quasiparticles. 
This means it is  metastable (see \cite{Masl_Shved}). 
The Bose-Einstein condensate is not involved here, 
it is only needed only to show that it would be wrong to 
believe that this argument works for a classical liquid, 
as one might think from the considerations above.

Indeed, for example, the molecules of a 
classical nondischarged liquid are, as a rule, Bose particles. 
For such a liquid, one can write out the $N$-particle 
equation, having in mind that each particle (molecule) is 
neutral and consists of an even number $l$ of neutrons. 
Thus each $i$th particle is a point in $3(2k+l)$-dimensional space, 
where
$k$ is the number of electrons, 
$x_i \in R^{6k+3l}$, depends on the potential 
$u(x_i), \ x_i \in R^{6k+3l}$ and we can consider 
the  $N$-particle equation for $x_i, i=1, \dots, N$,
with pairwise interaction $V(x_i-x_j)$.

\textbf{3.} However, there is a purely mathematical explanation 
of this paradox. The thing is that Bogolyubov found only one series of points 
in the spectrum of the many particle problem. 
Landau wrote ``N.N.Bogolyubov recently succeeded, 
by means of a clever application of second quantization, 
in finding the general form of the energy spectrum of a Bose-Einstein gas 
with weak interaction between the particles''
(\cite[p.~43]{landau}). But this series is not unique, i.e., 
the entire energy spectrum was not obtained.

In 2001, the author proposed the method of ultra 
second quantization \cite{Book_Ultravt};
see also
\cite{FAN_2000}, \cite{Uspehi_2000}, \cite{RJ_2001}, \cite{RJ_2001_2}, \cite{Mas1},
%[16], [17], [18], [19], [20].
The ultra second quantization of the Schr\"odinger equation, as well as its ordinary second quantization, is a representation of the $N$-particle
Schr\"odinger equation, and this means that basically the ultra second quantization of the  equation is the same as the original $N$-particle equation: they coincide in $3N$-dimen\-sional space. However, the replacement of the creation and annihilation  operators by $c$-numbers, in contrast with the case of second quantization, does not yield the correct asymptotics, but it turns out that it coincides with the result of applying the Schroeder variational principle or the Bogolyubov variational method.
{\tolerance=1000

}
For the exotic Bardin potential, the correct asymptotic solution 
coincides with the one obtained by applying the ultra second quantization 
method described above.
In the case of general potentials, in particular for pairwise 
interaction potentials, the answer is not the same. Specifically, 
the ultra second quantization method gives other asymptotic 
series of eigenvalues corresponding to the $N$-particle Schr\"odinger equation, 
and these eigenvalues, unlike the Bogolyubov  ones (7), are not metastable.

%Но если тор не является однородным, т.е. его образующие не равны
%между собой (как у Боголюбова), а равны $L_1, L_2,L_2$, где $L_1
%\gg L_2$, то для скоростей меньших $\min(\lambda_{\text{кр}},
%\frac{h}{mL_2})$, где $m$ -- масса частицы (в частности,
%молекулы), сверхтекучесть будет иметь место.

It turns out that the main point is not related 
to the Bose-Einstein condensate,
but has to do with the width of the capillary (the nanotube) through 
which the liquid flows. If we consider a liquid in a capillary or a nanotube 
of sufficiently small radius the velocity corresponding to metastable states is not small. 
Hence at smaller velocities the flow will be without friction.

The condition that the liquid does not flow through the boundary of the nanotube 
is a Dirichlet condition. It yields a standing wave, which can be 
regarded as a pair particle--antiparticle: a particle 
with momentum $p$ orthogonal to the boundary of the tube, 
and an antiparticle with momentum $-p$.

We consider a short action pairwise 
potential $V(x_i -x_j)$. This means that as the number of 
particles tends to infinity, $N\to\infty$, interaction is 
possible for only a finite number of particles. Therefore, the potential 
depends on $N$ in the following way: $$V_N=V((x_i-x_j)N^{1/3}).$$ 
If $V(y)$ is finite with support
$\Omega_V$, then as  $N\to\infty$ the support 
engulfs a finite number of particles, and this number 
does not depend on $N$.

As the result, it turns out that for velocities less than
$\min(\lambda_{\text{cr}}, \frac{h}{2mR})$, 
where $\lambda_{\text{cr}}$
is the critical Landau velocity and  $R$ is the 
radius of the nanotube, superfluidity occurs.

Now let me present my own considerations, which are not related to the mathematical exposition. 
Viscosity is due to the collision of particles: 
the higher the temperature, the greater the number of collisions. 
In a nanotube, there are few collisions, 
and only with the walls, and those are taken into 
account by the author's series. 
It is precisely this circumstance, and not the Bose-Einstein condensate, 
which leads to the weakening of viscosity and so to superfluidity. 
What I am saying is that the main factor in the superfluidity 
phenomenon, even for liquid helium 4, is not the condensate, 
but the presence of an extremely thin capillary \cite{TMF_2005}, \cite{RJ_2005}.
%[21], [22].
It seems to me that a neutral gas like argon could be used 
for a crucial experiment.
\renewcommand{\theequation}{\thesection.\arabic{equation}}
\setcounter{section}{0}
\setcounter{footnote}{0}
%\begin{center}
%{\bf Policy iteration}\\
%{\bf and max-plus finite element method}\\
%\vskip0.5cm
%{\it\large D.~McCaffrey}
%\end{center}
%\null
\nachalo{David McCaffrey}{Policy iteration and \mbox{max-plus} 
\mbox{finite element method}}
\markboth{David McCaffrey}{Policy Iteration \& Max-Plus FEM}
\label{mcc-abs}
\section{Introduction}

We consider the finite horizon differential game 
\begin{equation}
v(x,T)=\inf_{a(.)}\sup_{b(.)}\int_{0}^{T}\left\{ \frac{1}{2}x(s)^{2}+\frac{1%
}{2}a(s)^{2}-\frac{\gamma ^{2}}{2}b(s)^{2}\right\} ds+\phi (x(T))
\label{diff_game}
\end{equation}
over trajectories $(x(.),a(.),b(.))$ satisfying $\dot{x}%
(s)=f(x(s))+g(x(s))a(s)+h(x(s))b(s)$, $x(0)=x,$ where $x(s)\in X\subseteq %
\mathbb{R}^{n},\;\;a(s)\in U\subseteq \mathbb{R}^{m},\;\;b(s)\in W\subseteq %
\mathbb{R}^{r}.$This problem arises, for example, as the differential game
formulation of a well-known class of non-linear affine $H_{\infty }$ control
problems - see \cite{sora,van1,van2,mccaf_hinf} for details. In particular
it is known that the value function $v(x,t)$ for the finite horizon problem
is a (possibly non-smooth) solution to the Hamilton-Jacobi-Isaacs equation 
{\sloppy

}
\begin{equation}
H(x,\partial v/\partial x)=\partial v/\partial t  \label{diff_Ham}
\end{equation}
with initial condition $v(x,0)=\phi (x)$ for $(x,t)\in X\times (0,T],$ where
the Hamiltonian is defined as 
\begin{equation*}
H(x,p)=\min_{a}\max_{b}\left\{ p\left( f(x)+g(x)a+h(x)b\right) +\frac{1}{2}%
x^{2}+\frac{1}{2}a^{2}-\frac{\gamma ^{2}}{2}b^{2}\right\} .
\end{equation*}
Note this Hamiltonian is non-convex in $p.$

Suppose we choose some feedback function $\hat{a}(x)$ and, on any solution
trajectory $(x(.),a(.),b(.)),$ define the control input $a(s)=\hat{a}(x(s))$
for all $s.$ We can then define $f_{\hat{a}}(x,b)=f(x)+g(x)\hat{a}(x)+h(x)b$
and $l_{\hat{a}}(x,b)=\frac{1}{2}x^{2}+\frac{1}{2}\hat{a}(x)^{2}-\frac{%
\gamma ^{2}}{2}b^{2},$ and consider the finite horizon optimal control
problem 
\begin{equation}
v_{\hat{a}}(x,T)=\sup_{b(.)}\int_{0}^{T}l_{\hat{a}}(x(s),b(s))ds+\phi (x(T))
\label{opt_contr}
\end{equation}
over trajectories $(x(.),b(.))$ satisfying $\dot{x}(s)=f_{\hat{a}%
}(x(s),b(s)),\;\;x(0)=x.$ In this case, the value function $v_{\hat{a}}(x,t)$
satisfies the Hamilton-Jacobi equation 
\begin{equation}
H_{\hat{a}}(x,\partial v_{\hat{a}}/\partial x)=\partial v_{\hat{a}}/\partial
t  \label{opt_Ham}
\end{equation}
with initial condition $v_{\hat{a}}(x,0)=\phi (x)$ for $(x,t)\in X\times
(0,T],$ where the Hamiltonian is defined as $H_{\hat{a}}(x,p)=\max_{b}\left%
\{ pf_{\hat{a}}(x,b)+l_{\hat{a}}(x,b)\right\} .$ Note that this Hamiltonian
is convex in $p$ for all $x.$

A max-plus analogue of the finite element method (FEM) is set out in \cite
{akian} for the numerical computation of the value function $v_{\hat{a}}$
solving this convex optimal control problem (\ref{opt_contr}). In this note,
we set out a policy iteration algorithm for the solution of the non-convex
differential game (\ref{diff_game}). This involves the use of the max-plus
FEM to solve (\ref{opt_Ham}) for a given fixed control feedback $a(x)$ in
the value determination step of the algorithm, and then a QP to improve the
control feedback in the policy improvement step.We show here that the
algorithm converges. It can also be shown that the approximation error on
the converged solution is of order $\sqrt{\Delta t}+\Delta x(\Delta t)^{-1},$
the same order as that obtained in \cite{akian} for the errors associated
with the max-plus FEM. We do not give details of this result here, due to
limited space.

\section{The Max-Plus Finite Element Method}

In the following, let $S^{t}$ denote the evolution semi-group of the PDE (%
\ref{diff_Ham}). This associates to any function $\phi ,$ the function $%
v^{t}=v(.,t)$ where $v$ is the value function of the differential game (\ref
{diff_game}). Similarly, let $S_{\hat{a}}^{t}$ denote the evolution
semi-group of the PDE (\ref{opt_Ham}) for some fixed feedback function $\hat{%
a}(.)$. This associates to any function $\phi ,$ the function $v_{\hat{a}%
}^{t}=v_{\hat{a}}(.,t)$ where $v_{\hat{a}}$ is the value function of the
otimal control problem (\ref{opt_contr}). Maslov \cite{maslov87} observed
that the semi-group $S_{\hat{a}}^{t}$ is max-plus linear. We now briefly
review the max-plus finite element method (FEM) set out in \cite{akian} for
the numerical computation of $v_{\hat{a}}$

Let $\mathbb{R}_{\max }$ denote the idempotent semi-ring obtained from $%
\mathbb{R},$ with its usual order $\leq ,$ by defining idempotent addition
as $a\oplus b:=\max (a,b)$ and multiplication as $ab:=a+b.$ Then let $\bar{%
\mathbb{R}}_{\max }:=\mathbb{R}_{\max }\cup \{+\infty \},$ with the
convention that $-\infty $ is absorbing for the mutiplication.

For $X$ a set, we consider the set $\bar{\mathbb{R}}_{\max }^{X}$ of $\bar{%
\mathbb{R}}_{\max }$ valued functions on $X.$ This is a semimodule over $%
\bar{\mathbb{R}}_{\max }$ with respect to componentwise addition $%
(u,v)\longmapsto u\oplus v,$ defined by $(u\oplus v)(x)=u(x)\oplus v(x),$
and componentwise scalar multiplication $(\lambda ,u)\longmapsto u\lambda $,
defined by $(u\lambda )(x)=u(x)\lambda ,$ where $u,v\in \bar{\mathbb{R}}%
_{\max }^{X}$, $\lambda \in \bar{\mathbb{R}}_{\max }$ and $x\in X.$ Note
that the natural order on $\bar{\mathbb{R}}_{\max }^{X}$ arising from the
idempotent addition, i.e. the order defined by $u\leq v\Longleftrightarrow
u\oplus v=v,$ corresponds to the componentwise partial order $u\leq
v\Longleftrightarrow u(x)\leq v(x)$ for all $x\in X.$

Now let $X$ and $Y$ be sets and consider an operator $A:\bar{\mathbb{R}}%
_{\max }^{Y}\rightarrow \bar{\mathbb{R}}_{\max }^{X}$ from $\bar{\mathbb{R}}%
_{\max }$ valued functions on $Y$ to $\bar{\mathbb{R}}_{\max }$ valued
functions on $X.$ Such an operator is called linear if , for all $u_{1},$ $%
u_{2}\in \bar{\mathbb{R}}_{\max }^{Y}$ and $\lambda _{1},$ $\lambda _{2}\in 
\bar{\mathbb{R}}_{\max },$ $A(u_{1}\lambda _{1}\oplus u_{2}\lambda
_{2})=A(u_{1})\lambda _{1}\oplus A(u_{2})\lambda _{2}.$ Given some $\bar{%
\mathbb{R}}_{\max }$ valued function $a\in \bar{\mathbb{R}}_{\max }^{X\times
Y}$ on $X\times Y,$ we are then interested in the linear operator $A:\bar{%
\mathbb{R}}_{\max }^{Y}\rightarrow \bar{\mathbb{R}}_{\max }^{X}$ with kernel 
$a$ which maps any function $u\in \bar{\mathbb{R}}_{\max }^{Y}$ to the
function $Au\in \bar{\mathbb{R}}_{\max }^{X}$defined, in terms of the normal
arithmetic operations on $\mathbb{R},$ by 
\begin{equation}
Au(x)=\sup_{y\in Y}\left\{ a(x,y)+u(y)\right\}  \label{ker_def}
\end{equation}
Then, as shown in the references cited in \cite{akian}, this kernel operator 
$A$ is residuated, i.e. for any $v\in \bar{\mathbb{R}}_{\max }^{X},$ the set 
$\{u\in \bar{\mathbb{R}}_{\max }^{Y}:Au\leq v\}$ has a maximal element. The
residual map $A^{\#}:\bar{\mathbb{R}}_{\max }^{X}\rightarrow \bar{\mathbb{R}}%
_{\max }^{Y}$ then takes any $v\in \bar{\mathbb{R}}_{\max }^{X}$ to this
maximal element in $\bar{\mathbb{R}}_{\max }^{Y}$ defined, again in terms of
the normal arithmetic operations on $\mathbb{R},$ as the function 
\begin{equation}
(A^{\#}v)(y)=\inf_{x\in X}\left\{ -a(x,y)+v(x)\right\}  \label{resid_form}
\end{equation}

The next notion to be introduced, for a kernel operator $B:\bar{\mathbb{R}}%
_{\max }^{Y}\rightarrow \bar{\mathbb{R}}_{\max }^{X},$ is that of projection
on the image im$B$ of $B.$ The projector is denoted $P_{\text{im}B}$ and is
a map $\bar{\mathbb{R}}_{\max }^{X}\rightarrow \bar{\mathbb{R}}_{\max }^{X}$
defined for all $v\in \bar{\mathbb{R}}_{\max }^{X}$ by $P_{\text{im}%
B}(v)=\max \{w\in $im$B:w\leq v\}.$ Again as shown in the references cited
in \cite{akian}, this projector on the subsemimodule im$B$ can be expressed
as a composition $P_{\text{im}B}=B\circ B^{\#}$ of $B$ and its residual $%
B^{\#}.$ If $b(x,y)$ denotes the kernel of $B,$ then this formula can be
expressed in the normal arithmetic of $\mathbb{R},$ as 
\begin{equation}
B\circ B^{\#}(v)(x)=\sup_{y\in Y}\left( b(x,y)+\inf_{\xi \in X}\left( -b(\xi
,y)+v(\xi )\right) \right)   \label{proj1}
\end{equation}

Given a kernel operator $C:\bar{\mathbb{R}}_{\max }^{X}\rightarrow \bar{%
\mathbb{R}}_{\max }^{Z}$ with kernel $c(z,x),$ we can consider the
transposed operator $C^{*}:\bar{\mathbb{R}}_{\max }^{Z}\rightarrow \bar{%
\mathbb{R}}_{\max }^{X}$ with kernel $c^{*}(x,z)=c(z,x).$ We can then define
a dual projector on the $\bar{\mathbb{R}}_{\min }$-subsemimodule $-$im$C^{*}$
in terms of $P^{-\text{im}C^{*}}(v)=\min \{w\in -$im$C^{*}:w\geq v\}$ for
all $v\in \bar{\mathbb{R}}_{\max }^{X}.$ Then, as above, this projector can
be expressed as a composition $P^{-\text{im}C^{*}}=C^{\#}\circ C$ which, in
the normal arithmetic of $\mathbb{R},$ has the form 
\begin{equation}
C^{\#}\circ C(v)(x)=\inf_{z\in Z}\left( -c(z,x)+\sup_{\xi \in X}\left(
c(z,\xi )+v(\xi )\right) \right)  \label{proj2}
\end{equation}

Now we can define the max-plus FEM for approximating the value function $v_{%
\hat{a}}^{t}=v_{\hat{a}}(.,t)$ for the optimal control problem (\ref
{opt_contr}). Let $Y=\{1,\ldots ,I\},$ $X=\mathbb{R}^{n}$ and $Z=\{1,\ldots
,J\}.$ Consider a family $\{w_{1},\ldots ,w_{I}\}$ of finite element
functions $w_{i}:X\rightarrow \bar{\mathbb{R}}_{\max },$ and a family $%
\{z_{1},\ldots ,z_{J}\}$ of test functions $z_{j}:X\rightarrow \bar{%
\mathbb{R}}_{\max }.$ The vectors $\lambda =(\lambda _{i})_{i=1,\ldots
,I}\in \bar{\mathbb{R}}_{\max }^{I}$ and $\mu =(\mu _{j})_{j=1,\ldots ,J}\in 
\bar{\mathbb{R}}_{\max }^{J}$can be considered as $\bar{\mathbb{R}}_{\max }$
valued functions on $Y$ and $Z$ respectively. So, as above in equation (\ref
{ker_def}), we can define max-plus kernel operators $W:\bar{\mathbb{R}}%
_{\max }^{Y}\rightarrow \bar{\mathbb{R}}_{\max }^{X}$ and $Z^{*}:\bar{%
\mathbb{R}}_{\max }^{Z}\rightarrow \bar{\mathbb{R}}_{\max }^{X}$ with
kernels $W=$col$(w_{i})_{1\leq i\leq I}$ and $Z^{*}=$col$(z_{j})_{1\leq
j\leq J}.$ The action of $W,$ which plays the role of operator $B$ above, is
as follows 
\begin{equation*}
W\lambda (x)=\sup_{i\in Y}\left\{ w_{i}(x)+\lambda _{i}\right\}
\end{equation*}
while $Z^{*}$ gives rise to the transposed operator $Z:\bar{\mathbb{R}}%
_{\max }^{X}\rightarrow \bar{\mathbb{R}}_{\max }^{Z}$ which plays the role
of operator $C$ above, and acts as follows 
\begin{equation*}
(Zv)_{j}=\sup_{x\in X}\left\{ z_{j}(x)+v(x)\right\} =\left\langle
z_{j}|v\right\rangle
\end{equation*}
where $\left\langle .|.\right\rangle $ denotes the max-plus scalar product.
Then from equations (\ref{proj1}) and (\ref{proj2}), we can give the
specific form of the corresponding two projectors 
\begin{eqnarray}
P_{\text{im}W}(v)(x) &=&\sup_{i\in Y}\left( w_{i}(x)+\inf_{\xi \in X}\left(
-w_{i}(\xi )+v(\xi )\right) \right)  \label{Proj_W} \\
P^{-\text{im}Z^{*}}(v)(x) &=&\inf_{j\in Z}\left( -z_{j}(x)+\sup_{\xi \in
X}\left( z_{j}(\xi )+v(\xi )\right) \right)  \label{Proj_Z}
\end{eqnarray}

To start the algorithm off, we approximate the initial data $v_{\hat{a}%
}^{0}=\phi $ with the maximal element $\leq v_{\hat{a}}^{0}$ in the space im$%
W$ spanned by the finite element functions. The approximation of $v_{\hat{a}%
}^{0}$ is denoted with a subscript $h$ and takes the form 
\begin{equation*}
v_{\hat{a}h}^{0}(x)=(W\lambda ^{0})(x)=\sup_{i\in Y}\left( w_{i}(x)+\lambda
_{i}^{0}\right)
\end{equation*}
where the coefficients $\lambda _{i}^{0}$ are determined from the
residuation of $W$ given in formula (\ref{resid_form}) as 
\begin{equation}
\lambda _{i}^{0}=\inf_{x\in X}\left( -w_{i}(x)+\phi (x)\right) .
\label{lambda_nought}
\end{equation}
As an induction assumption, suppose that at time step $q\Delta t$ we have a
vector of coefficients $\lambda _{i}^{q\Delta t}$ giving an approximation 
\begin{equation*}
v_{\hat{a}h}^{q\Delta t}(x)=\sup_{i\in Y}\left( w_{i}(x)+\lambda
_{i}^{q\Delta t}\right)
\end{equation*}
of $v_{\hat{a}}^{q\Delta t}$ by the maximal element $\leq $ $v_{\hat{a}%
}^{q\Delta t}$ in the space im$W$. Then the approximation $v_{\hat{a}%
h}^{(q+1)\Delta t}$ of $v_{\hat{a}}^{(q+1)\Delta t}$ at the next time step
can be calculated as 
\begin{equation*}
v_{\hat{a}h}^{(q+1)\Delta t}(.)=P_{\text{im}W}\circ P^{-\text{im}Z^{*}}\circ
S_{\hat{a}}^{\Delta t}\circ v_{\hat{a}h}^{q\Delta t}(.)
\end{equation*}
The coefficients of this approximation are given, from equations (\ref
{Proj_W}) and (\ref{Proj_Z}), by 
\begin{equation*}
\lambda _{i}^{(q+1)\Delta t}=\inf_{\xi \in X}\left( -w_{i}(\xi )+\inf_{j\in
Z}\left( -z_{j}(\xi )+\sup_{\eta \in X}\left( z_{j}(\eta )+S_{\hat{a}%
}^{\Delta t}\circ v_{\hat{a}h}^{q\Delta t}(\eta )\right) \right) \right)
\end{equation*}
It is shown in \cite{akian} that $v_{\hat{a}h}^{(q+1)\Delta t}$ is the
maximal element in the space im$W$ spanned by the finite element functions
which satisfies $$\left\langle z_{j}|v_{\hat{a}h}^{(q+1)\Delta
t}\right\rangle \leq \left\langle z_{j}|S_{\hat{a}}^{\Delta t}v_{\hat{a}%
h}^{q\Delta t}\right\rangle $$ for each test function $z_{j}.$ So $v_{\hat{a}%
h}^{(q+1)\Delta t}$ is the maximal solution to a max-plus variational
formulation of the semi-group equation.

If (see Section 3.3 of \cite{akian}) we further approximate the semi-group
action $S_{\hat{a}}^{\Delta t}v_{\hat{a}h}^{q\Delta t}$ by 
\begin{equation*}
\left( \tilde{S}_{\hat{a}}^{\Delta t}v_{\hat{a}h}^{q\Delta t}\right)
(x)=\sup_{i\in Y}\left( w_{i}(x)+\lambda _{i}^{q\Delta t}+\Delta tH_{\hat{a}%
}(x,\partial w_{i}/\partial x)\right)
\end{equation*}
then $\lambda _{i}^{(q+1)\Delta t}$ can be written explicitly as 
\begin{eqnarray}
\lambda _{i}^{(q+1)\Delta t} &=&\inf_{\xi \in X}\left( -w_{i}(\xi
)+\inf_{j\in Z}\left( -z_{j}(\xi )+\sup_{\eta \in X}\left( z_{j}(\eta
)\right. \right. \right.  \label{lamda_equ} \\
&&\left. \left. \left. +\sup_{k\in Y}\left( w_{k}(\eta )+\lambda
_{k}^{q\Delta t}+\Delta tH_{\hat{a}}\left( \eta ,\partial w_{k}/\partial
x|_{\eta }\right) \right) \right) \right) \right)  \notag
\end{eqnarray}

Finally, choose two sets $(\hat{x}_{i})_{i\in Y}$ and $(\hat{x}_{j})_{j\in
Z} $ of discretisation points, and take the finite element functions to be $%
w_{i}(x)=-\frac{c}{2}\left\| x-\hat{x}_{i}\right\| _{2}^{2},$ for some fixed
Hessian $c,$ and test functions to be $z_{j}(x)=-a\left\| x-\hat{x}%
_{j}\right\| _{1},$ for some fixed constant $a.$ Then it is shown in Theorem
22 of \cite{akian} that the error $\left\| v_{\hat{a}h}^{T}-v_{\hat{a}%
}^{T}\right\| _{\infty }=O(\Delta t+\Delta x(\Delta t)^{-1})$, where $\Delta
x$ is the maximal radius of the cells of the two Voronoi tessellations
centred on the points $(\hat{x}_{i})_{i\in Y}$ and $(\hat{x}_{j})_{j\in Z}$
respectively.

\section{Policy Iteration with Max-Plus FEM in the Value Determination Step}

Now let $p$ denote the cycle index within the policy iteration algorithm,
and let $q\in \left\{ 0,\ldots ,N-1\right\} $ denote the time step index, so
that the full time horizon $T$ is divided into $N$ equal steps of length $%
\Delta t,$ i.e. $T=N\Delta t,$ with the $q$th step running from $q\Delta t$
to $(q+1)\Delta t$. We restrict consideration of time-dependent feedback
control policies $a(x,t)$ to those in the form of sequences of $N$
constant-in-time policy components $\left( a^{0}(x),\ldots
,a^{N-1}(x)\right) ,$ and we then further restrict our choice of the
individual policy components to functions $a^{q}(.)$ chosen from the set $%
A=\{a(.):X\rightarrow U\}$ of functions which are locally constant with
respect to $x$ on cells of the Voronoi tessellation $V_{Y}$ centred on the
origins $(\hat{x}_{i})_{i\in Y}$ of the finite element functions $w_{i}.$

So suppose, as an induction hypothesis, that on iteration $p$, we have a set
of constants $\{a_{p}^{qi}\}$ for $q\in \left\{ 0,\ldots ,N-1\right\} $ and $%
i\in Y.$ These give rise to a fixed policy $a_{p}$ which, for a given $q,$
takes the form $a_{p}^{q}(x)=a_{p}^{q\mu (x)}$ $,$ where $\mu (x)\in Y$ is
the index of the cell of the Voronoi tessellation $V_{Y}$ containing $x.$
Note, the process can be initiated, for $p=0,$ by choosing some fixed value $%
a^{i}$ (say zero) such that $a_{0}^{q}(x)=a^{i}$ for all $q\in \left\{
0,\ldots ,N-1\right\} $ and for all $x$ $\in $ cell $i$ of $V_{Y},$ where
cell $i$ is the one centred on the origin $\hat{x}_{i}$ of finite element $%
w_{i}.$

\subsection{Value Determination Step}

The max-plus FEM outlined above can be applied to approximate the value
function $v_{a_{p}}^{t}$ solving the optimal control problem (\ref{opt_contr}%
) with fixed strategy $a_{p}.$ The coefficients of the expansion of this
approximation, with respect to the finite elements $w_{i},$ are obtained as
follows. For $q=0$ and $i\in Y,$ the coefficients $\lambda _{pi}^{0}=\lambda
_{i}^{0}$ defined in (\ref{lambda_nought}) above. Then, using (\ref
{lamda_equ}), for $q\in \left\{ 0,\ldots ,N-1\right\} $ and $i\in Y$ we get 
\begin{eqnarray*}
\lambda _{pi}^{(q+1)\Delta t} &=&\inf_{\xi \in X}\left( -w_{i}(\xi
)+\inf_{j\in Z}\left( -z_{j}(\xi )+\sup_{\eta \in X}\left( z_{j}(\eta
)\right. \right. \right. \\
&&\left. \left. \left. +\sup_{k\in Y}\left( w_{k}(\eta )+\lambda
_{pk}^{q\Delta t}+\Delta tH_{a_{p}^{q}}\left( \eta ,\partial w_{k}/\partial
x|_{\eta }\right) \right) \right) \right) \right)
\end{eqnarray*}
Note that in the Hamiltonian $H_{a_{q}^{p}}$ we apply the policy $%
a_{p}^{q}(\eta )=a_{p}^{q\mu (\eta )}$ where $\mu (\eta )\in Y$ is the index
of the cell of the Voronoi tessellation $V_{Y}$ containing $\eta .$ The
above can be re-arranged to give 
\begin{eqnarray*}
\lambda _{pi}^{(q+1)\Delta t} &=&\inf_{j\in Z}\left( -\left\langle
w_{i}|z_{j}\right\rangle +\sup_{k\in Y}\left( \lambda _{pk}^{q\Delta
t}+\sup_{\eta \in X}\left( z_{j}(\eta )\right. \right. \right. \\
&&\left. \left. \left. +w_{k}(\eta )+\Delta tH_{a_{p}^{q}}\left( \eta
,\partial w_{k}/\partial x|_{\eta }\right) \right) \right) \right)
\end{eqnarray*}
For a given policy $a,$ let 
\begin{equation*}
T_{jka}=\sup_{\eta \in X}\left( z_{j}(\eta )+w_{k}(\eta )+\Delta
tH_{a}\left( \eta ,\partial w_{k}/\partial x|_{\eta }\right) \right)
\end{equation*}
In the normal max-plus FEM, the $T_{jka}$ terms can be calculated offline.
This would be difficult in the application of max-plus FEM to policy
iteration, since we don't know the policies $a$ in advance. The relevant $a$
for each $p$ iteration is known at the start of that iteration and so, in
principle, the next set of $T_{jka}$ terms for a given $a$ could be
calculated at the start of that iteration. However, this would be slow. An
alternative is to approximate the $T_{jka}$ online by 
\begin{equation*}
\tilde{T}_{jka}=\left\langle z_{j}|w_{k}\right\rangle +\Delta tH_{a}\left(
\eta _{jk}^{opt},\partial w_{k}/\partial x|_{\eta _{jk}^{opt}}\right)
\end{equation*}
where $\eta _{jk}^{opt}=\arg \sup \left\langle z_{j}|w_{k}\right\rangle
=\arg \sup \left( z_{j}(\eta )+w_{k}(\eta )\right) .$ Note, this
approximation $\tilde{T}$ is presented in \cite{akian}, where it is shown in
Theorem 22 that the resulting error estimate on the max-plus FEM
deteriorates to $\left\| v_{ah}^{T}-v_{a}^{T}\right\| _{\infty }=O(\sqrt{%
\Delta t}+\Delta x(\Delta t)^{-1}).$ So, finally, the coefficients of the
expansion of the approximation to the value function $v_{a_{p}}^{t}$ for
fixed strategy $a_{p}$ are given by 
\begin{equation}
\lambda _{pi}^{(q+1)\Delta t}=\inf_{j\in Z}\left( -\left\langle
w_{i}|z_{j}\right\rangle +\sup_{k\in Y}\left( \lambda _{pk}^{q\Delta t}+%
\tilde{T}_{jka_{p}^{q}}\right) \right)  \label{lambda_value_det}
\end{equation}

\subsection{Policy Improvement Step}

For each $i$ and $q,$ there exists $\bar{j}(iq)\in Z$ which achieves the $%
\inf $ in (\ref{lambda_value_det}), so that 
\begin{equation}
\lambda _{pi}^{(q+1)\Delta t}=-\left\langle w_{i}|z_{\bar{j}}\right\rangle
+\sup_{k\in Y}\left( \lambda _{pk}^{q\Delta t}+\tilde{T}_{\bar{j}%
ka_{p}^{q}}\right)  \label{lambda_value_sup}
\end{equation}
For each $k,$ the Hamiltonian within $\tilde{T}_{\bar{j}ka_{q}^{p}}$ is
evaluated at $\eta _{\bar{j}k}^{opt},$ and so the strategy $a_{p}^{q}$
applied in the Hamiltonian term takes the value $a_{p}^{q\mu (\bar{j}k)},$
where $\mu (\bar{j}k)\in Y$ is the index of the cell of the Voronoi
tessellation $V_{Y}$ containing $\eta _{\bar{j}k}^{opt}.$

The policy improvement can be formulated for test functions given by $%
z_{j}(x)=-a\left\| x-\hat{x}_{j}\right\| _{1}$ for some constant $a.$ Here,
due to lack of space, we consider only the special case where the constant
term $a\rightarrow \infty $ in the test functions $z_{j}(x)$, so that they
are therefore defined as 
\begin{equation}
z_{j}=\left\{ 
\begin{array}{l}
0\text{ at }x=\hat{x}_{j} \\ 
-\infty \text{ otherwise}
\end{array}
\right.   \label{zj_simple}
\end{equation}

Then we have $\eta _{\bar{j}k}^{opt}=\hat{x}_{\bar{j}}$ for all $k\in Y$ and 
$\mu (\bar{j}k)=\mu (\bar{j})\in Y$ is the index of the cell of the Voronoi
tessellation $V_{Y}$ containing $\hat{x}_{\bar{j}}.$ It follows that $%
a_{p}^{q}(\hat{x}_{\bar{j}})=a_{p}^{q\mu (\bar{j})}$ is the policy value
applied in the Hamiltonian term in $\tilde{T}_{\bar{j}ka_{p}^{q}}$ for all $%
k.$ So every term $\tilde{T}_{\bar{j}ka_{p}^{q}}$ uses the same policy value 
$a_{p}^{q\mu (\bar{j})}$ for all $k\in Y$ within the $\sup_{k\in Y}$
operation in (\ref{lambda_value_sup}).

Now let $\bar{k}(iq)=\arg \sup_{k\in Y}$ in (\ref{lambda_value_sup}), so
that 
\begin{equation*}
\lambda _{pi}^{(q+1)\Delta t}=-\left\langle w_{i}|z_{\bar{j}}\right\rangle
+\lambda _{p\bar{k}}^{q\Delta t}+\tilde{T}_{\bar{j}\bar{k}a_{p}^{q}}
\end{equation*}
Then we can improve the policy $a_{p}^{q}$ in cell $\mu (\bar{j})$ of $V_{Y}$
by taking 
\begin{equation}
\min_{a\in U}\tilde{T}_{\bar{j}\bar{k}a}  \label{optim}
\end{equation}
subject to 
\begin{equation}
\lambda _{p\bar{k}}^{q\Delta t}+\tilde{T}_{\bar{j}\bar{k}a}\geq \lambda
_{pk}^{q\Delta t}+\tilde{T}_{\bar{j}ka}  \label{constraint}
\end{equation}
for all $k\in Y.$ This optimisation is feasible since the current policy
value $a_{p}^{q\mu (\bar{j})}$ satisfies 
\begin{equation*}
\lambda _{p\bar{k}}^{q\Delta t}+\tilde{T}_{\bar{j}\bar{k}a_{p}^{q\mu (\bar{j}%
)}}\geq \lambda _{pk}^{q\Delta t}+\tilde{T}_{\bar{j}ka_{p}^{q\mu (\bar{j})}}
\end{equation*}
for all $k\in Y.$

Let $\bar{a}=\arg \min_{a\in U}\tilde{T}_{\bar{j}\bar{k}a}$ subject to the
constraints (\ref{constraint}). In cell with index $\mu (\bar{j})$ of the
Voronoi tessellation $V_{Y},$ take new policy 
\begin{equation*}
a_{p+1}^{q}(x)=a_{p+1}^{q\mu (\bar{j})}:=\bar{a}
\end{equation*}
for all $x\in $ cell with index $\mu (\bar{j}).$ Note that for each $q,$
there may be some remaining cells of $V_{Y}$ whose indices $\neq \mu (\bar{j}%
(iq))$ for any $i\in Y.$ In these cells we leave the policy at time step $q$
unchanged, i.e. if $\mu ^{*}$ is the index of such a cell, then for all $%
x\in $ cell with index $\mu ^{*}$%
\begin{equation*}
a_{p+1}^{q}(x)=a_{p}^{q\mu ^{*}}
\end{equation*}
Then the resulting new policy $a_{p+1}=\{a_{p+1}^{qi}\}$ is an improvement
on the old one $a_{p}=\{a_{p}^{qi}\}$ in the sense that the corresponding $%
v_{a_{p}h}^{q\Delta t}$ and $v_{a_{(p+1)}h}^{q\Delta t},$ i.e. the
approximations to the value functions which solve the optimal control
problem (\ref{opt_contr}) with fixed policies $a_{p+1}$ and $a_{p}$
respectively, satisfy 
\begin{equation}
v_{a_{(p+1)}h}^{q\Delta t}\leq v_{a_{p}h}^{q\Delta t}  \label{imporve}
\end{equation}
for all $q\in \left\{ 0,\ldots ,N\right\} .$

To see this, note first that the policy improvement is unique. If, for a
given $q,$ there are two $i$ giving rise to the same $\bar{j}(iq),$ then
these both result in the same policy improvement $a_{p+1}^{q}(x)=\bar{a}$ in
cell $\mu (\bar{j})$ since the term $\tilde{T}_{\bar{j}ka_{p}^{q}}$ in (\ref
{lambda_value_sup}) does not depend on $i.$

Next, suppose with a view to induction on $q,$ that $\lambda
_{(p+1)i}^{q\Delta t}\leq \lambda _{pi}^{q\Delta t}$ for all $i.$ Then 
\begin{eqnarray*}
\lambda _{pi}^{(q+1)\Delta t} &=&-\left\langle w_{i}|z_{\bar{j}%
}\right\rangle +\lambda _{p\bar{k}}^{q\Delta t}+\tilde{T}_{\bar{j}\bar{k}%
a_{p}^{q}} \\
&=&-\left\langle w_{i}|z_{\bar{j}}\right\rangle +\lambda _{p\bar{k}%
}^{q\Delta t}+\tilde{T}_{\bar{j}\bar{k}a_{p}^{q\mu (\bar{j})}} \\
&\geq &-\left\langle w_{i}|z_{\bar{j}}\right\rangle +\lambda _{p\bar{k}%
}^{q\Delta t}+\tilde{T}_{\bar{j}\bar{k}\bar{a}} \\
&=&-\left\langle w_{i}|z_{\bar{j}}\right\rangle +\sup_{k\in Y}\left( \lambda
_{pk}^{q\Delta t}+\tilde{T}_{\bar{j}k\bar{a}}\right)  \\
&\geq &-\left\langle w_{i}|z_{\bar{j}}\right\rangle +\sup_{k\in Y}\left(
\lambda _{(p+1)k}^{q\Delta t}+\tilde{T}_{\bar{j}k\bar{a}}\right)  \\
&=&-\left\langle w_{i}|z_{\bar{j}}\right\rangle +\sup_{k\in Y}\left( \lambda
_{(p+1)k}^{q\Delta t}+\tilde{T}_{\bar{j}ka_{p+1}^{q}}\right) 
\end{eqnarray*}
since every term $\tilde{T}_{\bar{j}ka_{p+1}^{q}}$ uses the same policy
value in the same cell $\mu (\bar{j}).$ So 
\begin{equation*}
\lambda _{pi}^{(q+1)\Delta t}\geq \inf_{j}\left( -\left\langle
w_{i}|z_{j}\right\rangle +\sup_{k\in Y}\left( \lambda _{(p+1)k}^{q\Delta t}+%
\tilde{T}_{jka_{p+1}^{q}}\right) \right) =\lambda _{(p+1)i}^{(q+1)\Delta t}
\end{equation*}
Since this holds for all $i,$ then it follows that for any $x\in X$%
\begin{equation*}
\sup_{i}\left( \lambda _{(p+1)i}^{(q+1)\Delta t}+w_{i}(x)\right) \leq
\sup_{i}\left( \lambda _{pi}^{(q+1)\Delta t}+w_{i}(x)\right) 
\end{equation*}
i.e. $v_{a_{(p+1)}h}^{(q+1)\Delta t}(x)\leq v_{a_{p}h}^{(q+1)\Delta t}(x).$
So by induction, after noting that from (\ref{lambda_nought}), the
coefficients at the initial time step $q=0$ satisfy $\lambda
_{(p+1)i}^{0}=\lambda _{pi}^{0}=\lambda _{i}^{0},$ it follows that $\lambda
_{(p+1)i}^{q\Delta t}\leq \lambda _{pi}^{q\Delta t}$ and $%
v_{a_{(p+1)}h}^{q\Delta t}(x)\leq v_{a_{p}h}^{q\Delta t}(x)$ for all $q.$

Finally by taking $b=0$ in (\ref{opt_contr}), and by restricting our choice
of initial data to functions $\phi \geq 0,$ we can see that $%
v_{a_{p}}^{t}\geq 0$ for all $p.$ Since $\left\|
v_{a_{p}h}^{t}-v_{a_{p}}^{t}\right\| _{\infty }=O(\sqrt{\Delta t}+\Delta
x(\Delta t)^{-1}),$ it follows that 
\begin{equation}
v_{a_{p}h}^{t}\geq -K(\sqrt{\Delta t}+\Delta x(\Delta t)^{-1})
\label{lower_bdd}
\end{equation}
for some $K>0.$

Hence, the above policy iteration algorithm converges to a time-dis\-cre\-ti\-sed
finite element approximation $v_{h}^{t}(.)=\left\{ v_{h}^{\Delta
t}(.),\ldots ,v_{h}^{N\Delta t}(.)\right\} $ to the value function $v^{t}$
of the differential game (\ref{diff_game}).

\subsection{QP Optimisation}

The Hamiltonian appearing in (\ref{opt_Ham}) has the form 
\begin{eqnarray*}
H_{a}(x,p) &=&\max_{b}\left\{ pf_{a}(x,b)+l_{a}(x,b)\right\} \\
&=&p(f+ga)+\frac{1}{2}x^{2}+\frac{1}{2}a^{2}-\frac{1}{2\gamma ^{2}}phh^{T}p
\end{eqnarray*}
and for $z_{j}$ given by (\ref{zj_simple}), $\tilde{T}_{\bar{j}\bar{k}a}$
has the specific form 
\begin{equation*}
\tilde{T}_{\bar{j}\bar{k}a}=w_{\bar{k}}(\hat{x}_{\bar{j}})+\Delta
tH_{a}\left( \hat{x}_{\bar{j}},\partial w_{\bar{k}}/\partial x|_{\hat{x}_{%
\bar{j}}}\right)
\end{equation*}
The policy improvement optimisation set out in (\ref{optim}) and (\ref
{constraint}) can then be formulated as 
\begin{equation*}
\min_{a\in U}\Delta tH_{a}\left( \hat{x}_{\bar{j}},\partial w_{\bar{k}%
}/\partial x|_{\hat{x}_{\bar{j}}}\right)
\end{equation*}
subject to 
\begin{eqnarray*}
\Delta tH_{a}\left( \hat{x}_{\bar{j}},\partial w_{\bar{k}}/\partial x|_{\hat{%
x}_{\bar{j}}}\right) \geq \lambda _{pk}^{q\Delta t}-\lambda _{p\bar{k}%
}^{q\Delta t}+w_{k}(\hat{x}_{\bar{j}})-\\
-w_{\bar{k}}(\hat{x}_{\bar{j}})+\Delta
tH_{a}\left( \hat{x}_{\bar{j}},\partial w_{k}/\partial x|_{\hat{x}_{\bar{j}%
}}\right)
\end{eqnarray*}
for all $k\in Y.$ This can be simplified down to the following QP 
\begin{equation*}
\min_{a\in U}\left( \frac{\partial w_{\bar{k}}}{\partial x}ga+\frac{1}{2}%
a^{2}\right)
\end{equation*}
subject to 
\begin{eqnarray*}
\left( \frac{\partial w_{\bar{k}}}{\partial x}-\frac{\partial w_{k}}{%
\partial x}\right) ga &\geq &\frac{1}{\Delta t}\left( w_{k}-w_{\bar{k}%
}\right) +\frac{1}{\Delta t}\left( \lambda _{pk}^{q\Delta t}-\lambda _{p\bar{%
k}}^{q\Delta t}\right) \\
&&+\frac{1}{2\gamma ^{2}}\left( \frac{\partial w_{k}}{\partial x}-\frac{%
\partial w_{\bar{k}}}{\partial x}\right) ^{T}hh^{T}\left( \frac{\partial
w_{k}}{\partial x}-\frac{\partial w_{\bar{k}}}{\partial x}\right) \\
&&+\left( \frac{\partial w_{k}}{\partial x}-\frac{\partial w_{\bar{k}}}{%
\partial x}\right) f
\end{eqnarray*}
for all $k\in Y,$ and evaluated at $\hat{x}_{\bar{j}}.$

\setcounter{footnote}{0}
\setcounter{section}{0}
%\begin{center}
%{\bf\large Using Max-Plus Convolution to Obtain Fundamental Solutions  }\\
%{\bf\large for Differential Equations with Quadratic Nonlinearities}
%\footnote{Research partially supported by NSF grant DMS-0307229
%and AFOSR grant FA9550-06-1-0238.}\\
%\vskip0.5cm
%{\it\large William M. McEneaney}
%\end{center}
%\null
\nachaloe{William M.~McEneaney}{Using max-plus convolution to 
obtain fundamental solutions for differential equations
with quadratic nonlinearities}{Research partially supported by NSF grant DMS-0307229\\
and AFOSR grant FA9550-06-1-0238.}
\label{mce-abs}

\section{Introduction}
\label{sec:intro}

%The matrix differential Riccati equation (DRE) is ubiquitous  in control and systems theory.
We first consider time-invariant differential Riccati equations (DREs) of the form
\begin{equation}
\dot P_t=F(P_t)\doteq\Aprime P_t+P_tA+C+P_t\Sigma P_t
\label{eq:introric}\end{equation}
where $C$ is symmetric and $\Sigma=\sigma\sigma^\prime$ is symmetric, nonnegative definite
with at least one positive eigenvalue.
Throughout, we assume that all of the matrices are $n\times n$. We suppose one has initial condition, $P_0=p_0$
where $p_0$ is also symmetric.
The Daivson-Maki approach uses the Bernoulli substitution to create a linear system of two
matrices, each of the same size as $P_t$, thus leading to a fundamental solution.
%We denote these matrices as $V^1_t$ and $V^2_t$.
%The solution of the DRE is recovered as $P_t=V_t^2[V_t^1]^{-1}$.
%Thus, one may obtain a fundamental solution of the DRE.
We obtain a completely different form of fundamental solution, with 
a particularly clear control-theoretic motivation.
The new approach will be constructed through a finite-dimensional semigroup defined by
this fundamental solution. The forward propagation of the fundamental solution is naturally defined by
this operation through the semigroup property.

We will consider linear/quadratic control problems parameterized by $z\in\liiRn$,
and the value functions associated with these
control problems 
are propagated forward by
a \maxp linear semigroup, which we denote as $\soltau$.
The space of semiconvex functions
is a \maxp vector space (moduloid) \cite{Mcmpbook}, \cite{FMcmp}, \cite{BCOQ}, \cite{CGQ}, \cite{LMS}. 
Working in the semiconvex-dual space, $\soltau$ has a semiconvex-dual operator,
$\cB_\tau$ which takes the form of a \maxp integral operator with kernel, $B_\tau(x,y)$,
taking the form of a quadratic function. The matrix, $\beta_\tau$, defining this quadratic kernel function
will be the fundamental solution of the DRE. We will define a multiplication
operation ($\circledast$-multiplication)\markboth{William M.~McEneaney}{Max-plus convolution and fundamental solution}
with the semigroup property, specifically
$\beta_{t+\tau}=\beta_t\circledast\beta_\tau$,
where the $\circledast$ operation involves inverse, multiplication and addition 
$n\times n$-matrix operations (in the standard algebra).
We will also define an exponentiation operation ($\circledast$-exponentiation)
such that $\beta_t=\beta_1^{\circledast t}$.
The solution of \er{eq:introric} will be obtained by
$P_t=\matdualopinv\beta_t\matdualop p_0$
where the $\matdualop$ and $\matdualopinv$ operators are descended from
the semiconvex dual and its inverse.
It is important to note that the fundamental solution approach has the benefit that one only solves once for $\beta_t$, even if one wishes to solve the DRE for a variety of initial conditions.

This approach may be extended to a class of quasilinear, first-order PDEs, yielding a fundamental solution
for a class of such PDEs.
More specifically, we consider PDEs
\begin{equation}
0=-P_t+A^\prime P-BP_\lambda +\demi P\Sigma P
\end{equation}
on the domain
$[0,T]\times\cL$ where $\cL\doteq [0,L]$.
For simplicity, we consider the scalar case, and so $A,B,\Sigma \in\liiR$, $\Sigma>0$,
where we specifically require $B\not=0$ (otherwise this reduces to an ODE problem).
%The initial and boundary conditions are
%\beas
%P(0,\lambda)\midalign=p_0(\lambda)\qforall \lambda\in\cL\\
%P(t,\hat\lambda)\midalign=0\qforall (t,\hat\lambda)\in\cE^{B}_P
%\eeas
%where
%$$\cE^{B}_P=
%\begin{cases}
%(0,T]\times\{0\} & \mbox{ if } B>0\cr
%(0,T]\times\{L\} & \mbox{ if } B<0.\cr
%\end{cases}$$

We again create a linear/quadratic control problem, where in this case, the state takes values in $L_2(\cL)$.
The above PDE is essentially the ``Riccati'' equation for this virtual control problem. We again apply semiconvex duality, and \maxp vector space concepts. This leads to an extension of the $\circledast$ operator
to this infinite-dimensional context, and finally, a fundamental solution for this class of PDEs.

\section{The linear-quadratic control problem and semigroups}
\label{sec:lqprob}

\noindent The proofs of the results in the sections on the DRE may be found in
\cite{Mcricfund}.

As indicated above, the fundamental solution to the DRE will be obtained through an associated optimal control problem.
Recall that we are considering the DRE given by \er{eq:introric}.
Since we will be employing semiconvex duality (see below and \cite{FMcmp,Mcmpbook}),
we will choose some (duality-parametrizing) symmetric matrix, $Q$, such that
$F(Q)>0$,
where we note that,
for any square matrix $D$,
we will use the notation $D>0$ to indicate that matrix $D$ is positive definite throughout.
We will need to consider the specific solution of DRE \er{eq:introric} with initial condition
\begin{equation}
\Ptilde_0=Q.
\label{eq:fundPic}\end{equation}
We assume:

\noindent
\parbox{0.84\hsize}{\em
There exists a solution of DRE \er{eq:introric}, $\Ptilde_t$, with initial condition \er{eq:fundPic},
satisfying $\Ptilde_t>Q$ (i.e., $\Ptilde_t-Q$ positive-definite) for $t\in(0,\Tbar)$
with $\Tbar>0$, and we note specifically,
that we may have $\Tbar=+\infty$.
}
\parbox{0.11\hsize}{$\hskip 1.5em(A.e)$}

%\noindent This is the last assumption in the finite-dimensional (ordinary DRE) portion of the work.

\noindent We will be obtaining the fundamental solution $\beta_t$ for solutions with initial conditions, $P_0=p_0>Q$.
Note that we do not assume stability of the DRE, and finite-time blow-up is possible.
We will let
$\Ttilde=\Ttilde(p_0)=\sup\{t\ge 0\,\,\vert\,\,P_t\,\mbox{ exists, and }\, P_t>Q\}$,
and we let
$\That=\That(p_0)\doteq\Tbar\wedge\Ttilde$
where $\wedge$ indicates the minimum operation.

\begin{remark}
Note that with $\Sigma\ge 0$ and at least one positive eigenvalue,
we may take $Q=-kI$ for arbitrarily large $k$,
so that one can ensure $F(Q)>0$ (as well as for any $p_0>Q$).
\end{remark}

We will be using a control value function to motivate and develop the fundamental solution.
Consider the Hamilton-Jacobi-Bellman partial differential equation (HJB PDE) problems
on $[0,\Tbar)\times\liiRn$, indexed by $z\in\liiRn$, given by
\beasnum
&& \Vz_t=H(x,\grad\Vz)=(Ax)^\prime \grad\Vz+\demi\xprime C x+(\grad\Vz)^\prime\Sigma \grad\Vz\label{eq:hjbpde}\\
&& \Vz(0,x)=\psi(z,x)=\demi(x-z)^\prime Q(x-z).\label{eq:hjbic}
\eeasnum

\begin{theorem}
For any $z\in\liiRn$,
there exists a solution to \er{eq:hjbpde},\er{eq:hjbic} in $C^\infty([0,\Tbar)\times\liiRn)\cap
C([0,\Tbar)\times\liiRn)$,
and this is given by
\begin{equation}
\Vz=\demi(x-\Lambda_tz)^\prime \Ptilde_t(x-\Lambda_t z)+\zprime R_t z
\label{eq:vzquadform}\end{equation}
where $\Ptilde$ satisfies \er{eq:introric},\er{eq:fundPic},
and $\Lambda,r$ satisfy
$\Lambda_0=I$, $R_0=0$,
\begin{equation}
\dot\Lambda=\left[\Ptilde^{-1}C-A\right]\Lambda
\quad\mbox{ and }\quad
\dot R=\Lambda^\prime C\Lambda.
\label{eq:lamrodes}\end{equation}
\end{theorem}

For $\phi:\liiRn\rightarrow\liiR$ given by $\phi(x)=\demi(x-z)^\prime p_0(x-z)$
(and actually for a much larger set of functions), we define the \maxp linear semigroup, $\soltau$, by
\begin{equation}
\soltau[\phi](x)
=\Vz(\tau,x)=\demi(x-\Lambda_\tau z)^\prime \Ptilde_\tau(x-\Lambda_\tau z)+\zprime R_\tau z.
\label{eq:sgasquad}
\end{equation}

%\section{Solution via the semiconvex dual semigroup}
%\label{sec:scdualsol}

We let $\oplus,\otimes$ denote the \maxp addition and multiplication operations.
%We say that $\phi:\,\liiRn\rightarrow\rem\doteq\liiR\cup\{-\infty\}$ is semiconvex if
%given $R<\infty$, there exists finite, symmetric $C_R>0$ such that
%$\phi(x)+\demi\xprime C_R x$ is convex on $B_R(0)$.
We say that $\phi$ is uniformly semiconvex with (symmetric matrix) constant $K$ if
$\phi(x)+\demi\xprime K x$ is convex on $\liiRn$, and we denote this space as $\cS^K(\liiRn)$.
Recall that $\cS^K$ is a \maxp vector space.
% (c.f., \cite{Mcmpbook}).

%Semiconvex duality is parameterized by quadratic functions.
We will use the quadratic $\psi$ given in \er{eq:hjbic} to define our semiconvex duality.
The main duality result (c.f., \cite{Mcmpbook}, \cite{FMcmp}, where proofs may be found) is

\begin{theorem}\label{th:scduality}
Let $\phi\in\cS^K(\liiRn)$ where $-K>Q$.
Then, for all $x,z\in\liiRn$,
%\beasnum
\begin{equation}
\begin{split}
\phi(x)=\max_{z\in\liiRn}\left[\psi(x,z)+a(z)\right]
%\label{eq:semicdual1a}\\
%\nonumber\\
%\midalign
&\doteq\intmprn\psi(x,z)\otimes a(z)\,dz\\
%\label{eq:semicdual1b}\\
%\midalign
&\doteq\psi(x,\cdot)\odot a(\cdot)
\doteq\dualopinv[a]
\label{eq:semicdual1c}
\end{split}
\end{equation}
%\textline{where for all $z\in\liiRn$}
%%%\nonumber\\
\begin{equation}
\begin{split}
a(z)&=
%-\max_{x\in\liiRn}\left[\psi(x,z)-\phi(x)\right]
%\label{eq:semicdual2a}\\
%\nonumber\\
%\midalign= 
-\intmprn\psi(x,z)\otimes [-\phi(x)]\,dx\\
%\label{eq:semicdual2b}\\
%\midalign
&=-\left\{\psi(\cdot,z)\odot[-\phi(\cdot)]\right\}\doteq\dualop[\phi].
\label{eq:semicdual2c}
\end{split}
\end{equation}
%\\
%\textline{which using the notation of \cite{CGQ}}
%\midalign=\left\{\psi(\cdot,z)\odot[\phi^-(\cdot)]\right\}^-
%\doteq\dualop[\phi].
%\label{eq:semicdual2d}
%%%\eeasnum
\end{theorem}

%Recall that $\Ptilde_t>Q$ for all $t\in(0,\Tbar)$,
%and consequently, for any $t\in(0,\Tbar)$,
%there exists 
%$K_t$ such that $\Ptilde_t>-K_t>Q$
%(i.e., such that $\Ptilde_t+K_t>0$ and $-Q-K_t>0$),
%%$-K_t>Q$
%and such that (noting \er{eq:sgasquad})
%$$\solt[\psicdotz](\cdot)\in\cS^{K_t}\qforall z\in\liiRn.$$
%Therefore, by Theorem \ref{th:scduality},
Using  Theorem \ref{th:scduality} and some technical arguments,
for all $t\in(0,\Tbar)$ and all $x,z\in\liiRn$
\begin{equation}
\begin{split}
\solt[\psicdotz](x)&=
\intmprn\psi(x,y)\otimes B_t(y,z)\,dy\\
&=\psixcdot\odot B_t(\cdot,z),\label{eq:tildedualop1}
\end{split}
\end{equation}
where for all $y\in\liiRn$
\begin{equation}
\begin{split}
B_t(y,z)&=
-\intmprn\psi(x,y)\otimes
\bigl\{-\solt[\psicdotz](x)\bigr\}\,dx\\
&=\bigl\{\psi(\cdot,y)\odot
[\solt[\psicdotz](\cdot)]^-\bigr\}^-.
\label{eq:tildedualop2}
\end{split}
\end{equation}

%\begin{lemma}
%There exists symmetric $d_t<-Q$ such that
%$B_t(y,z)-\demi y^\prime d_t y$ is strictly concave for all $z\in\liiRn$.
%\end{lemma}

%In analogy to the spaces of uniformly semiconvex functions,
%we say that $\phi:\,\liiRn\rightarrow\liiR^+\doteq\liiR\cup\{+\infty\}$ is uniformly semiconcave
%with (symmetric matrix) constant $d$ if
%$\phi(x)-\demi\xprime d x$ is concave on $\liiRn$, and we denote this space as $\cS^d_-(\liiRn)$.
We define the time-indexed \maxp linear operators $\cB_t$ by
\begin{equation}
\cB_t[a](z)\doteq
B_t(\cdot,z)\odot a(\cdot)=\intmprn B_t(y,z)\otimes a(y)\,dy,
\label{eq:btopdef}\end{equation}
and one easily sees that these satisfy the semigroup property.
(We may use a space of uniformly semiconcave functions as the domain.)
We say that $B_t$ is the kernel of \maxp integral operator $\cB_t$.

\begin{theorem}
Let $\phi(x)\doteq\demi \xprime p_0 x$ and $a(z)=\dualop[\phi]$.
Then, for $t\in(0,\That)$, $x\in\liiRn$,
\begin{equation}
\solt[\phi](x)=\psi(x,\cdot)\odot\cB_t[a](\cdot)
=\dualopinv\cB_t[a](x)=\dualopinv\cB_t\dualop[\phi](x).
\label{eq:dualprop}\end{equation}
\end{theorem}

Now, note that by \er{eq:sgasquad} and \er{eq:tildedualop2},
\begin{equation}
B_t(x,y)= -\max_{x\in\liiRn}\biggl\{
\demi(x-y\parenprime Q(x-y)
%\nonumber\\ \midalign
-\left[\demi(x-\Lambda_tz\parenprime \Ptilde_t(x-\Lambda z)+\demi \zprime R_t z\right]
\biggr\}
\label{eq:basmax}\end{equation}
where $t<\Tbar$ guarantees strict concavity of the argument of the maximum.

\begin{lemma}\label{lem:maxofquads}
Let $\eta$ and $\alpha$ be $2n\times 2n$ matrices with block structure given by
\begin{equation}
\eta=
\left[\begin{matrix}
\etaoo&\etaot\cr
\etaot^\prime&\etatt\cr
\end{matrix}\right]
\quad\mbox{ and }\quad
\alpha=
\left[\begin{matrix}
\alphaoo&\alphaot\cr
\alphaot^\prime&\alphatt\cr
\end{matrix}\right],
\label{eq:etaalphanote}\end{equation}
Let
$$
F(x,z)\doteq
\max_{z\in\liiRn}
\left\{
\demi
\left(
\begin{matrix}x\cr
z\cr
\end{matrix}
\right)^\prime\eta
\left(
\begin{matrix}x\cr
z\cr
\end{matrix}
\right)
+\demi
\left(
\begin{matrix}z\cr
y\cr\end{matrix}
\right)^\prime
\alpha
\left(
\begin{matrix}z\cr
y\end{matrix}
\right)
\right\}.$$
Then,
$$F(x,y)=\demi
\left(
\begin{matrix}x\cr
y\end{matrix}
\right)^\prime
\gamma
\left(
\begin{matrix}x\cr
y\end{matrix}
\right)$$
where $\gamma$ has identical block structure to $\eta$ and $\alpha$, and is given by
$\gamma=\eta\circledast\alpha$ where the $\circledast$ operation is defined as
\beas
&&\gammaoo=\etaoo-\etaot S^{-1}\etaot^\prime,
%\\
%&&
\qquad
\gammaot=-\etaot S^{-1}\alphaot,\\
&&\gammato=\gammaot^\prime,
%\\
%&&
\qquad
\gammatt=\alphatt-\alphaot^\prime S^{-1}\alphaot,
\eeas
and $S\doteq\etatt+\alphaoo$.
\end{lemma}

Combining \er{eq:basmax} and Lemma \ref{lem:maxofquads}, one obtains the following.

\begin{theorem}\label{th:bandbeta}
\begin{equation}
B_t(x,y)=\demi
\left(
\begin{matrix}x\cr
y\end{matrix}
\right)^\prime
\beta_t
\left(
\begin{matrix}x\cr
y\end{matrix}
\right)
\label{eq:bbeta}\end{equation}
where $\beta_t$ has the same block structure as $\eta$ above.
% and in particular,
%
%\vskip -0.1in
%\parbox{0.84\hsize}
%{\beas
%&&\betatoo=Q\Deltatinv Q-Q=Q\Deltatinv\Ptilde_t,\\
%&&\betatot=-\Ptilde_t\Deltatinv Q\Lambda_t=-Q\Deltatinv\Ptilde_t\Lambda_t,\\
%&&\betatto=\betatot^\prime,\\
%&&\betattt=\Lambda_t^\prime\Ptilde_t\Lambda_t+R_t+
%\Lambda_t^\prime\Ptilde_t\Deltatinv\Ptilde_t\Lambda_t
%=R_t+\Lambda_t^\prime Q\Deltatinv\Ptilde_t\Lambda_t,
%\eeas}
%\parbox{0.12\hsize}{\beasnum\label{eq:betafromp}\eeasnum}
%and $\Delta_t\doteq Q-\Ptilde_t$.
\end{theorem}

\section{The DRE fundamental solution semigroup}
\label{sec:fundsolsg}

Now we will use the semigroup nature of the $\solt$ operators to obtain the semigroup nature of the
$\cB_t$ operators, and consequently the propagation of the $B_t$ and $\beta_t$.
The propagation of $\beta_t=(\beta_1)^{\circledast t}$ will be the dynamics of the
fundamental solution of the DRE.

\begin{lemma}\label{lem:dualforquads}
Let $a(z)=\demi(z-\zbar\parenprime q_a(z-\zbar)+r_a$ with $q_a<-Q$,
and $\phi=\dualopinv a$.
Then,
\begin{equation}
\phi(x)=\demi(x-\zbar\parenprime \left[Q \Uinv q_a\right](x-\zbar)+r_a
\label{eq;pfroma}\end{equation}
where $U=Q+q_a$.
Alternatively, let $\phi(x)=\demi(x-\xbar\parenprime q_p(x-\xbar)+r_p$
with $q_p>Q$,
and let $a=\dualop\phi$.
Then, with $\Delta\doteq Q-q_p$
\begin{equation}
a(z)=\demi(z-\xbar\parenprime \left[Q\Delta^{-1}q_p\right] (z-\xbar)+r_p.
\label{eq:afromphi}\end{equation}
%where $\Delta=Q-q_p$.
\end{lemma}

Based on this lemma, it is natural to make the following definitions,
which inherit notation from $\dualop$ and $\dualopinv$. 
For symmetric $q_p>Q$, define $D_\psi[q_p]\doteq Q(Q-q_p)^{-1}q_p$,
and for symmetric $q_a<-Q$, define $D_\psi^{-1}[q_a]= Q(Q+q_a)^{-1} q_a$.
One may show (see \cite{Mcricfund}):

\begin{theorem}\label{th:bsg}
For all $t_1,t_2\ge 0$ such that $t_1+t_2<\Tbar$,
$$\Btoptt(\zeta,x)=\intmprn\Bto(\zeta,z)\otimes\Btt(z,x)\,dz
\qforall x,\zeta\in\liiRn.$$
\end{theorem}

%Further, one obtains:

\begin{theorem}\label{th:betasgprop}
The forward propagation of semigroup $\beta_t$ is given by
\begin{equation}
\betasubtoptt=\betasubto\circledast\betasubtt
\label{eq:betasgpop}\end{equation}
where the $\circledast$ operation is given in Lemma  \ref{lem:maxofquads}.
%as above,
%\beas
%&&\left[\betasubto\circledast\betasubtt\right]^{1,1}=\betasubtooo-\betasubtoot U_{t_1,t_2}^{-1}\betasubtoot^\prime\\
%&&\left[\betasubto\circledast\betasubtt\right]^{1,2}=-\betasubtoot U_{t_1,t_2}^{-1}\betasubttot\\
%&&\left[\betasubto\circledast\betasubtt\right]^{2,1}=
%-\betasubttot^\prime U_{t_1,t_2}^{-1}\betasubtoot^\prime\\
%%{\left[\betasubto\circledast\betasubtt\right]^{1,2}}^\prime\\
%&&\left[\betasubto\circledast\betasubtt\right]^{2,2}=\betasubtttt-\betasubttot^\prime U_{t_1,t_2}^{-1}\betasubttot
%\eeas
%where $U_{t_1,t_2}\doteq\betasubtott+\betasubttoo$.
\end{theorem}

To summarize, suppose one wishes to obtain the solution of \er{eq:introric} at time $t$
with initial condition $P_0=p_0$.
Then, one performs the following steps:

\begin{itemize}

\item Obtain $q_0$ from $p_0$ via
$q_0=\matdualop p_0=Q(Q-p_0)^{-1} p_0$.

\item Obtain $q_t$ from $\beta_t$ and $q_0$ via
$
q_t
=
\beta_t^{1,1}-\beta_t^{1,2}\left(\beta_t^{2,2}+q_0\right)^{-1}{\beta_t^{1,2}}^\prime
%=\left\{
%\beta_t\circledast
%\left[
%\begin{matrix} q_0 & 0\cr
%0 & 0\end{matrix}
%\right]
%\right\}^{1,1}
\doteq
\beta_t\circledast^\prime q_0$.

\item Obtain $P_t$ from from $q_t$ via
$P_t=\matdualopinv q_t=Q(Q+q_t)^{-1}q_t$.

\end{itemize}

%\noindent This may be repeated for any number of appropriate initial conditions, $p_0$.
%The choice of symmetric $Q$ (which parametrizes semiconvex duality) is partially free, only needing to
%be sufficient to ensure existence of the semiconvex duals.

%\begin{remark}
%It is worth noting that if one worked in the dual space, one simply has
%$q_t=\beta_t\circledast^\prime q_0$ as the solution of the (dual of the) DRE.
%\end{remark}

\section{Exponentiation and a Semiring}\label{sec:fundsolprop}

%We have seen that this fundamental solution propagates according to the matrix operation
%$\betasubtoptt=\betasubto\circledast\betasubtt$.
Recall that for a standard-algebra linear system, one views the fundamental solution as
$e^{At}=(e^A)^t$.
%, and so $(e^A)^{t_1+t_2}=(e^A)^{t_1}\cdot(e^A)^{t_2}$.
%Thus, $\beta_t$ is analogous to $(e^A)^t$, and
We would like some similar exponential-type representation here.
Naturally, we define $\circledast$-exponentiation for positive integer powers through
$\beta\capow{2}=\beta\circledast\beta$,
$\beta\capow{3}=[\beta\capow{2}]\circledast\beta$,
et cetera.
Using Theorem \ref{th:betasgprop}, this immediately yields
$\beta_{nt}=\beta_t\capow{n}$.
However, this only works for integer powers. We will extend this to positive real powers so that we
may simply write
$\beta_t=(\beta_1)\capow{t}$ for any $t>0$.

Let $\cQ$ denote the set of rationals.
Given any $t\in(0,\infty)$, let $e_t\doteq\{s\in(0,\infty)\,|\,\exists p\in\cQ\,\mbox{ such that }\, s=pt\}$.
As is well-known, the collection of such $e_t$ forms an uncountable set of equivalence classes
covering $(0,\infty)$.
Suppose $s\in e_t$. Then, there exists $p=m/n$ with $m,n\in\cN$ such that $s=pt$.
Let $\tau=t/n$. Then, $t=n\tau$ and $s=m\tau$.
Consequently, by Theorem \ref{th:betasgprop},
$\beta_s=\beta_\tau\capow{m}$ and $\beta_t=\beta_\tau\capow{n}$.
With this in mind, we make the following extension of $\circledast$-exponentiation.
%to rationals,
%and the fact that this extension is well-defined will be proved immediately below.

\begin{definition}
Let $s=pt$ with $p=m/n$, $m,n\in\cN$. We define
$\beta_t\capow{p}\doteq\beta_\tau\capow{m}$ where $\tau=t/n$.
\end{definition}

We need to demonstrate that the definition is independent of the choice of $m,n\in\cN$.
That is, suppose $p=m_0/n_0=m_1/n_1$.
Let $\tau_0=t/n_0$ and $\tau_1=t/n_1$.
We must show $\beta_{\tau_0}\capow{m_0}=\beta_{\tau_1}\capow{m_1}$.
We will use the following, trivially-verified result.

\begin{lemma}\label{lem:astexpmult}
$\left[\beta_t\capow{n}\right]\capow{m}=\beta_t\capow{(nm)}$.
\end{lemma}

With this lemma, the above independence is easily proven.
Lastly, one extends the $\circledast$-exponentiation definition to exponents which may not be rational
using continuity.

%\section{$\langle\oplus,\circledast\rangle$-Semirings}\label{sec:semiring}

There are underlying semirings with the $\oplus,\circledast$ operations,
and this seems to be quite interesting.
These semirings are related to the convolution semiring of \cite{KM}.
We only touch on the matter here.
Let $a,b\in[0,+\infty)\cup\{+\infty\}\doteq\cWp$.
Then define
$a\circledast b\doteq ab/(a+b)$
which defines the $\circledast$ operation on $\cWp$.
Also, define $\oplus$ on $\cWp$ by $a\oplus b=\max\{a,b\}$.

\begin{theorem}\label{th:wpsemiring}
$\langle\cWp,\oplus,\circledast\rangle$ is a
commutative idempotent semiring.
\end{theorem}

\section{First-Order Quasilinear PDE}\label{sec:qlpde}

The same approach, which was used above in the case of the DRE, can be applied to
a first-order, quasilinear PDE with a quadratic nonlinearity. This PDE will take the form of
a Riccati equation, and we will refer to it as the fully-first-order Riccati PDE
(the FFOR PDE). The FFOR PDE will be
\begin{equation}
0=-P_t+A P-BP_\lambda +\demi P\Sigma P
\label{eq:fforpdedef}\end{equation}
where the domain will be
$[0,T]\times\cL$ where $\cL\doteq [0,L]$.
For simplicity, we consider the scalar case, and so $A,B,\Sigma \in\liiR$, $\Sigma>0$,
where we specifically require $B\not=0$ (otherwise this reduces to an ODE problem).
We let
$$\cE^{B}_P=
\begin{cases}
(0,T]\times\{0\} & \mbox{ if } B>0\cr
(0,T]\times\{L\} & \mbox{ if } B<0.\cr
\end{cases}$$
The initial and boundary conditions will be
\beasnum
P(0,\lambda)\midalign=p_0(\lambda)\qforall \lambda\in\cL
\label{eq:fforic}\\
P(t,\hat\lambda)\midalign=0\qforall (t,\hat\lambda)\in\cE^{B}_P.
\label{eq:fforbc}
\eeasnum

We will obtain a fundamental solution for \er{eq:fforpdedef}--\er{eq:fforbc}
using technology analogous to that used for the DRE.
In order to do so, we must devise a virtual control problem for which the time-reversed version 
of \er{eq:fforpdedef} is the associated Riccati equation.
%We will also need to generalize the
%semiconvex dual-space propagation to the infinite-dimensional domain case required here.

We begin by defining the dynamics of the virtual control problem.
The state will take values in $\cX=L_2(\cL;\liiR)$.
The control will take values in $\cW=L_2(\cL;\liiR)$.
In particular, we consider the control space
$\cW^s=L_2([-T,0],\cW)$.
The domain for the dynamics will be $[-T,0]\times\cL$,
and for $t\in[-T,0]$, we will have state, $\xi(t,\cdot)\in\cX$.
Let
$$\cE^{B}_X=
\begin{cases}
(-T,0]\times\{L\} & \mbox{ if } B>0\cr
(-T,0]\times\{0\} & \mbox{ if } B<0.\cr
\end{cases}$$
The virtual control problem dynamics is given by first-order PDE, initial condition and
boundary condition 
\beasnum
&&\xi_t(t,\lambda)=A\xi(t,\lambda)+B\xi_\lambda(t,\lambda)+\sigma w(t,\lambda),
\label{eq:vcpdyn}\\
&&\xi(-T,\cdot)=x_0(\cdot)\in\cX
\label{eq:vcpic}\\
&&\xi(t,\hat\lambda)=0\qforall (t,\hat\lambda)\in\cE^{B}_X.
\label{eq:vcpbc}
\eeasnum

Let the inner product and norm on $\cX,\cW$ be denoted by $\langle x,y\rangle$ and
$\|x\|$, respectively. Let $C,Q>0$.
The payoff and value are given by
\beasnum
&&J^z(-T,x,w)=\int_{-T}^0\demi \langle\xi(t,\cdot),C\xi(t,\cdot)\rangle
-\demi\|w(t,\cdot)\|^2\,dt+\psi(\xi(0,\cdot),z)
\label{eq:vcppayoff}\\
&&\Whatz(-T,x)=\sup_{w\in\cW^s}J^z(-T,x,w)
\label{eq:vcpvalue}
\eeasnum
where $\psi(x,z)\doteq\demi\langle\xi(0,\cdot)-z,Q(\xi(0,\cdot)-z)\rangle$.
%\textline{where}
%&&\psi(x,z)\doteq\demi\langle\xi(0,\cdot)-z,Q(\xi(0,\cdot)-z)\rangle
%\label{eq:idpsidef}
%\eeasnum
%and we note that $z\in\cX$.

The first step is to obtain the verification result.

\begin{theorem}\label{th:idver}
Suppose $V^z\in C([0,T]\times\cX)\cap C^1((0,T]\times\cX)$ satisfies
\begin{align}
&\begin{array}{l@{{}+{}}l}
0=-V^z_t&\langle(\grad_x V^z,Ax\rangle
-\langle(\grad_x v)_\lambda,Bx\rangle+\\[1.5ex]
&\demi\|\sigma^\prime\grad_x V^z\|^2
+\demi\langle x,Cx\rangle,
\end{array}\label{eq:iddpe}\\
&V^z(0,x)=\psi(x,z),\label{eq:iddpeic}\\
&\grad_x V^z(t,\lambdahat,x(\lambdahat))=0\qforall (t,\lambdahat)\in\cE^{B}_P,
\,x\in H^1(\cL).\label{eq:iddpebc}
\end{align}
Then, $V^z(T,x)\ge J^z(-T,x,w)$ for all $w\in\cW$, for all $x\in H^1(\cL)$.
Further, if there exists a solution, $\xi^*$ to \er{eq:vcpdyn}--\er{eq:vcpbc}
with $w^*(t,\lambda,\xi^*(t,\lambda))$, then letting
$\tilde w^*(t,\lambda)=w^*(t,\lambda,\xi^*(t,\lambda))$, one has
$V^z(T,x)= J^z(-T,x,\tilde w^*)$, and consequently $V^z(T,x)=\Whatz(-T,x)$.
\end{theorem}

%The proof is a straight-forward verification with some calculus of variations components.

The next step is to note that the solution for this problem has a simple form.

\begin{theorem}\label{id:form}
For any $z\in\cX$, there is a solution to \er{eq:iddpe}--\er{eq:iddpebc}
in $V^z\in C([0,T]\times\cX)\cap C^1((0,T]\times\cX)$ of the form
\begin{equation}
V^z(t,x)=\demi\langle (x-Z(t,\cdot)), \Ptilde(t,\cdot)(x-Z(t,\cdot))\rangle
+\demi\langle Z(t,\cdot),R(t)Z(t,\cdot)\rangle,\label{eq:vzform}
\end{equation}
where $\Ptilde$ satisfies \er{eq:fforpdedef}--\er{eq:fforbc} with
\begin{align}
&p_0(\lambda)=Q\qforall \lambda\in\cL,
\label{eq:pzQ}\\
\intertext{$Z\in C^1([0,T]\times\cL)$ satisfies}
&0=\Ptilde Z_t+[A\Ptilde+C]Z+B\Ptilde Z_\lambda,\label{eq:zpde}\\
&Z(0,\lambda)=z(\lambda)\qforall \lambda\in\cL,\label{eq:zic}\\
&Z(t,\hat\lambda)=0\qforall (t,\hat\lambda)\in\cE^{B}_X,\label{eq:zbc}\\
\intertext{and $R\in C^1([0,T])$ satisfies}
&R_t=C,
%\label{eq:rpde}\\
%&&
\qquad R(0)=0.\label{eq:ric}
%&&R(0,\lambda)=0\qforall\lambda\in\cL,\label{eq:ric}\\
%&&R(t,\hat\lambda)=0\qforall (t,\hat\lambda)\in\cE^{B}_P.\label{eq:rbc}
\end{align}
\end{theorem}

%The proof follows by simple substitution into \er{eq:iddpe}--\er{eq:iddpebc}.

Note that in the case where $C=0$, the $Z$ PDE takes the simpler form
\begin{equation}
0=Z_t+AZ+BZ_\lambda.
\label{eq:zpdealt}\end{equation}
Note also that $Z(t,\cdot)$ is given by a linear operator acting on $z$, denoted as
\begin{equation}
Z(t,\cdot)=\cM(t)[z](\cdot)\doteq\int_\cL M(t;\cdot,\eta)z(\eta)\,d\eta.
\label{eq:zlop}\end{equation}
Using this in \er{eq:vzform},
one has
\begin{equation}
V^z(t,x)=\demi\langle (x-\cM(t)z), \Ptilde(t,\cdot)(x-\cM(t)z)\rangle
+\demi\langle \cM(t)z,R(t)\cM(t)z\rangle.\label{eq:vzwithmop}\end{equation}

We may think of $V^z(t,\cdot)$ as given by the \maxp linear semigroup
$V^z(t,x)=\solt[\psi(\cdot,z)](x)$.
Introducing the semiconvex dual, one may propagate instead in the dual space.
The dual-space semigroup operator is naturally found in the form of a \maxp integral operator with some kernel,
which we denote by $B(t;x,z)$ for $t\in[0,T]$ and $x,z\in\cX$.
One obtains $B(t;x,z)$ from
\begin{multline*}
B(t;y,z)=
-\max_{x\in\cX}\biggl\{
\demi\langle (y-x),Q(y-x)\rangle\\
-\left[\demi\langle (x-\cM(t)[z],\Ptilde(x-\cM(t)[z]\rangle
+\demi\langle \cM(t)z,R(t)\cM(t)z\rangle
\right]\biggr\}.
\end{multline*}

\noindent Further details will appear in the full paper.

\setcounter{section}{0}
\setcounter{footnote}{0}
\nachaloe{Vladimir F.~Molchanov}{ Polynomial quantization on \mbox{para-hermitian}
\mbox{symmetric} \mbox{spaces} from the \mbox{viewpoint} of
\mbox{overgroups: an example}}{Supported by the Russian Foundation for Basic Research: grants 
No. 05-01-00074a, No. 05-01-00001a and 07-01-91209 YaF\_a, the Netherlands 
Organization for Scientific Research (NWO): grant 047-017-015, the Scientific 
Program "Devel. Sci. Potent. High. School":
project RNP.2.1.1.351 and Templan No. 1.2.02.}
\markboth{Vladimir F.~Molchanov}{Polynomial quantization on 
para-hermitian spaces}
\label{mol-abs}
Quantization in the spirit of Berezin on para-Hermitian symmetric spaces $G/H$ was 
constructed by the author in [2]. One of the variants of quantization is the so-called 
polynomial quantization (here for the initial algebra of operators, one has to 
take a representation of the universal enveloping algebra). A construction of 
polynomial quantization on para-Hermitian symmetric spaces $G/H$ was 
%old:
%offered
presented in [4]. 
For rank one, explicit formulas were given in [3]. In this paper we consider a new 
approach to the polynomial quantization using the notion of an "overgroup". This 
approach gives the Berezin covariant and contravariant symbols and the Berezin 
transform in a highly natural and transparent way. In the paper we restrict 
ourselves to a simple but crucial example: $G={\rm SL}(2,\Bbb R)$ with the diagonal 
subgroup $H$ and ${\widetilde G}=G\times G$.

%\null

\section{Groups, subgroups, a cone, sections}

The group $G={\rm SL}(2,\Bbb R)$ consists of real matrices
\begin{equation}
g=\left(
\begin{array}{cc}
\alpha&\beta \\
\gamma&\delta \\
\end{array}
\right), \ \ \ \alpha\delta-\beta\gamma=1. \nonumber
\end{equation} 
Its subgroups $H$, $Z$, $N$ of $G$ consist of matrices
\begin{equation}
h=\left(
\begin{array}{cc}
\alpha&0 \\
0&\alpha^{-1} \\
\end{array}
\right), \ \ 
z_\xi=\left(
\begin{array}{cc}
1&0 \\
\xi&1 \\
\end{array}
\right), \ \ 
n_\eta=\left(
\begin{array}{cc}
1&\eta \\
0&1 \\
\end{array}
\right),  \nonumber
\end{equation} 
respectively. 
%Old%
%There are the Gauss- and "anti-Gauss"-decompositions of $G$: 
%$G={\overline{NHZ}}$ and $G={\overline{ZHN}}$. 
%%%%%%%%%%%%%%%%%%
%New
The Gauss and "anti-Gauss" decompositions of $G$ are defined by
$G={\overline{NHZ}}$ and $G={\overline{ZHN}}$.
%%%%%%%%%%%%%%%%%
The group $G$ acts on $Z$ and $N$ by 
fractional linear transformations: 
\begin{equation}
\xi\mapsto \xi\cdot g= \frac{\alpha\xi+\gamma}{\beta\xi+\delta}, \ \ \ 
\eta\mapsto \eta\circ g=\frac{\delta\eta+\beta}{\gamma\eta+\alpha}.  \nonumber
\end{equation}
These actions are obtained when we decompose $z_\xi g$ "by Gauss" and $n_\eta g$ "by 
anti-Gauss". We can reduce the second action to the first one: 
$\eta\circ g=\eta\cdot\widehat g$ where
\begin{equation}
\widehat g=\left(
\begin{array}{cc}
\delta&\gamma \\
\beta&\alpha \\
\end{array}
\right).  \nonumber
\end{equation} 
We 
%Old
% consider
%New
assume
that the groups act from the right, in accordance with this we will write 
vectors in the row form.

%Old
%Let us take in the space ${\Bbb R}^4$ the following bilinear form:
%New
Let us take the following bilinear form in the space ${\Bbb R}^4$:
\begin{equation}
[x,y]=-x_1y_1-x_2y_2+x_3y_3+x_4y_4.  \nonumber
\end{equation}
Realize ${\Bbb R}^4$ as the space ${\rm Mat}(2,\Bbb R)$ of real $2\times 2$ 
matrices:
\begin{equation}
x=\frac{1}{2}\left(
\begin{array}{cc}
x_1-x_4&-x_2+x_3 \\
x_2+x_3&x_1+x_4 \\
\end{array}
\right).  \nonumber
\end{equation}
Denote the matrix corresponding to the vector $xJ$,  
$J={\rm diag}\{1,-1,-1,-1\}$, by $x^\natural$.  Then the form $[x,y]$ can be written in terms of 
matrices: $[x,y]= -2 \, {\rm tr}\left(x^\natural y\right)$.

As an overgroup for $G$, we take the direct product ${\widetilde G}=G\times G$.  
It acts on ${\rm Mat}(2,\Bbb R)$ as follows: to a pair 
$(g_1,g_2)\in {\widetilde G}$ we assign the transformation
\begin{equation}
x\mapsto g_1^{-1}xg_2.  
\end{equation}
This action preserves ${\rm det} \ x=-(1/4)[x,x]$. 
Therefore ${\widetilde G}$ covers the group ${\rm SO}_0(2,2)$ with multiplicity 2, 
%New
and
%%%%%%
the kernel of the homomorphism consists of two pairs: $(e,e)$ and $(-e,-e)$,   
$e$ being the unit matrix in $G$.

Let ${\mathcal C}$ be the cone in ${\Bbb R}^4$ defined by $[x,x]=0, x\ne 0$ (or  
${\rm det} \ x=0, \ x\ne 0$). Let us take the following two points in $\mathcal C$:
\begin{equation}
s^-=(1,0,0,-1)=
\left(
\begin{array}{cc}
1&0 \\
0&0 \\
\end{array}
\right), \ \ \ 
s^+=(1,0,0,1)=
\left(
\begin{array}{cc}
0&0 \\
0&1 \\
\end{array}
\right), \nonumber
\end{equation}
and two parabolic sections $\Gamma ^{-}=\{[x, s^{+}]=-2\}$ and 
$\Gamma ^{+}=\{[x, s^{-}]=-2\}$ containing $s^-$ and $s^+$ respectively.

Consider in ${\widetilde G}$ two unipotent subgroups $Q^-$ and $Q^+$ consisting of 
pairs $(z_\xi, n_\eta)$ and $(n_\eta, z_\xi)$ respectively. They act simply 
transitively on sections $\Gamma^{-}$ and $\Gamma^{+}$ respectively and transfer points 
$s^-$ and $s^+$ to the points 
\begin{eqnarray}
&u&=u(\xi, \eta)=(1-\xi \eta,\, -\xi-\eta,\, -\xi+\eta,\,-1-\xi \eta),  \nonumber\\
&v&=v(\xi, \eta)=(1-\xi \eta,\, \xi+\eta,\, \xi-\eta,\,1+\xi \eta),  \nonumber
\end{eqnarray}
respectively.  Let $u=u(\xi_1,\eta_1)$ and $v=v(\xi_2,\eta_2)$, then
\begin{equation}
[u,v]=-2N(\xi_1,\eta_2)N(\xi_2,\eta_1),
\end{equation}
where
\begin{equation}
N(\xi, \eta)=1-\xi \eta. \nonumber
\end{equation}
In terms of matrices, the vectors $u$ and $v$ are
%New
written as follows:
$$
u=\left(
\begin{array}{cc}
1&\eta \\
-\xi&-\xi\eta \\
\end{array}
\right)
=\left(
\begin{array}{c}
1 \\
-\xi \\
\end{array}
\right)(1 \ \  \ \eta),  \ \
v=
\left(
\begin{array}{cc}
-\xi\eta&-\eta \\
\xi&1 \\
\end{array}
\right)
=\left(
\begin{array}{c}
-\eta \\
1 \\
\end{array}
\right)(\xi \ \  \ 1),  \nonumber
$$
%Old
%They are connected by: $u=vJ$ or $u=v^{\natural}$. 
%%%%%%%%%
The relation between $u$ and $v$ is given by: $u=vJ$ or $u=v^{\natural}$. 

Let ${\mathcal X}$ be the section of the cone ${\mathcal C}$ by the plane $x_1=1$ (or 
 ${\rm tr} \ x=1)$. It is a hyperboloid of one sheet: $-x_2^2+x_3^2+x_4^2=1$, in 
${\Bbb R}^3$. 

Using maps of points along generating lines in ${\mathcal C}$, 
we obtain actions of 
$\widetilde G$ on sections.
%rasstavit' artikli?
Let $(g_1,g_2)\in\widetilde G$. For ${\mathcal X}$ we have:
\begin{equation}
x\mapsto {\widetilde x}=\frac{g_1^{-1}xg_2}
{{\rm tr}(g_1^{-1}xg_2)}.
\end{equation}
For $\Gamma^{-}$ and $\Gamma^{+}$, these actions are given by fractional linear 
transformations of $\xi$ and $\eta$:
\begin{equation}
u(\xi, \eta)\mapsto u(\xi\cdot g_1, \eta\circ g_2), \ \ \ 
v(\xi, \eta)\mapsto v(\xi\cdot g_2, \eta\circ g_1).  \nonumber
\end{equation}

Each of 
%New
the 
%%%
sections $\Gamma^{-}$ and $\Gamma^{+}$ is mapped on ${\mathcal X}$ along 
%New
the 
%%%
generating lines (almost everywhere):
\begin{equation}
u\mapsto x=\frac{u}{u_1}=\frac{u(\xi, \eta)}{N(\xi, \eta)}, \ \
v\mapsto y=\frac{v}{v_1}=\frac{v(\xi, \eta)}{N(\xi, \eta)}.  \nonumber
\end{equation}
These maps give 
%New
the following
%%%%
two systems of coordinates $\xi, \eta$ on ${\mathcal X}$:
\begin{eqnarray}
&&x=\Big(1, -\,\frac{\xi+\eta}{N(\xi, \eta)}, \ -\,\frac{\xi-\eta}{N(\xi, \eta)}, \
 -\,\frac{1+\xi \eta}{N(\xi, \eta)}\Big),  \nonumber\\
&&y=\Big(1, \frac{\xi+\eta}{N(\xi, \eta)}, \ \frac{\xi-\eta}{N(\xi, \eta)}, \
 \frac{1+\xi \eta}{N(\xi, \eta)}\Big).  \nonumber
\end{eqnarray}
Let us call these coordinates the {\it horospherical coordinates} corresponding to 
$\Gamma^{-}$ and $\Gamma^{+}$  respectively. 
%Old
%These systems are connected by: 
%%%%%%
%New
The relation between these two systems is given by:
$x=yJ$ or $x=y^\natural$.

Let us take the following measures on the sections  ${\mathcal X}$, $\Gamma^{-}$ and 
$\Gamma^{+}$:
\begin{equation}
dx=|x_4|^{-1}dx_2 \, dx_3, \ \ \ du=d\xi \,d\eta, \ \ \ dv=d\xi \,d\eta.  \nonumber
\end{equation}
Under the maps mentioned above, the measures are 
%Old
%linked
%%%%
related as follows: 
\begin{equation}
dx=dy=2N(\xi, \eta)^{-2}d\xi \,d\eta.  \nonumber
\end{equation}

%\null

\section{Representations of $G={\rm SL}(2,\Bbb R)$}

%\null 
%New
The 
representations $T_{\sigma,\varepsilon}$, $\sigma\in\Bbb C$, $\varepsilon=0,1$, 
of $G$ act on functions $\varphi (\xi)$ on ${\Bbb R}$ by: 
\begin{equation}
(T_{\sigma,\varepsilon}(g)\varphi )(\xi)=\varphi (\xi\cdot g)
(\beta\xi+\delta)^{2\sigma,\varepsilon},  \nonumber
\end{equation}
%New
where
%%%%
we use the notation:
\begin{equation}
t^{\lambda,\nu}=|t|^\lambda {\rm sgn}^\nu t, \ \ {\lambda\in\Bbb C}, \ \ \nu=0,1, \ \ 
t\in{\Bbb R}\setminus\{0\}.  \nonumber
\end{equation}
%Old
%Equally
%%%%
%New
Together
 with these representations, we consider the representations 
$${\widehat T}_{\sigma,\varepsilon}(g)=T_{\sigma,\varepsilon}(\widehat g),$$
so that
\begin{equation}
({\widehat T}_{\sigma,\varepsilon}(g)\psi )(\eta)=\psi (\eta\circ g)
(\gamma\eta+\alpha)^{2\sigma,\varepsilon}  \nonumber
\end{equation}
(notice that ${\widehat T}_{\sigma,\varepsilon}$ and $T_{\sigma,\varepsilon}$ are 
equivalent). The operator $A_{\sigma,\varepsilon}$ defined by 
\begin{equation}
(A_{\sigma,\varepsilon}\varphi)(\eta)=
\int^\infty_{-\infty}(1-\xi\eta)^{-2\sigma-2,\varepsilon}\varphi(\xi) \, d\xi,
 \nonumber
\end{equation}
intertwines $T_{\sigma,\varepsilon}$ 
with ${\widehat T}_{-\sigma-1,\varepsilon}$ 
and also ${\widehat T}_{\sigma,\varepsilon}$ 
with $T_{-\sigma-1,\varepsilon}$. The 
product $A_{-\sigma-1,\varepsilon}A_{\sigma,\varepsilon}$ 
is a scalar operator:
{\sloppy

}
\begin{equation}
A_{-\sigma-1,\varepsilon}A_{\sigma,\varepsilon}=
\omega(\sigma,\varepsilon)\cdot {\rm id},  \nonumber
\end{equation}
where 
\begin{equation}
\omega(\sigma,\varepsilon)=
\frac{2\pi}{2\sigma+1}{\rm tan}\frac{2\sigma-\varepsilon}{2}\pi.  \nonumber
\end{equation}
Notice that 
\begin{equation}
\omega(-\sigma-1,\varepsilon)=\omega(\sigma,\varepsilon).
\end{equation}

%\null
\section{Representations of ${\widetilde G}=G\times G$ associated with a cone}

For $\lambda\in {\Bbb C}$, $\nu=0,1$, let ${\mathcal D}_{\lambda, \nu}({\mathcal C})$ 
denote
the space of functions $f\in C^\infty(\mathcal C)$ satisfying the condition:
\begin{equation}
f(tx)=t^{\lambda,\nu} f(x), \ \ x\in{\mathcal C}, \ \ t\in{\Bbb R}\setminus\{0\}.  
\nonumber
\end{equation}
Let $R_{\lambda,\nu}$ be the representation of ${\widetilde G}$ on 
${\mathcal D}_{\lambda,\nu}({\mathcal C})$ by translations: 
\begin{equation}
(R_{\lambda,\nu}(g_1,g_2)f)(x)=f(g_1^{-1}xg_2).  \nonumber
\end{equation}
In fact, it is the representation of the group ${\rm SO}_0 (2,2)$ associated with 
a cone [1]. This representation can be realized on functions on sections 
of the cone $\mathcal C$, see \S \ 1. 
%Old
In the realization on ${\mathcal X}$,
%%%%
%Realized on ${\mathcal X}$,
%%%%%%
the representation 
$R_{\lambda,\nu}$ is given by (see (1.3)):
\begin{equation}
(R_{\lambda,\nu}(g_1,g_2)f)(x)=
f({\widetilde x})\left\{{\rm tr}(g_1^{-1}xg_2)\right\}^{\lambda,\nu},
 \ x\in {\mathcal X}.  \nonumber
\end{equation}
On $\Gamma^{-}$ and on $\Gamma^{+}$ we have respectively:
\begin{eqnarray}
&&(R_{\lambda,\nu}(g_1,g_2)f)(\xi,\eta)=
f(\xi\cdot g_1,\eta\circ g_2)
\Big\{(\beta_1\xi+\delta_1)(\gamma_2\eta+\alpha_2)\Big\}^{\lambda,\nu}, \\
&& (R_{\lambda,\nu}(g_1,g_2)f)(\xi,\eta)=
f(\xi\cdot g_2,\eta\circ g_1)
\Big\{(\beta_2\xi+\delta_2)(\gamma_1\eta+\alpha_1)\Big\}^{\lambda,\nu}.
\end{eqnarray}
Formulas (3.1) and (3.2) show that $R_{\lambda,\nu}(g_1,g_2)$ is the tensor product 
$T_{\sigma,\varepsilon}(g_1)\otimes {\widehat T}_{\sigma,\varepsilon}(g_2)$ and  
$T_{\sigma,\varepsilon}(g_2)\otimes {\widehat T}_{\sigma,\varepsilon}(g_1)$ 
respectively with $\sigma = \lambda/2$.

%Old
%Define an operator $B_{\lambda, \, \nu}$ in the ${\mathcal X}$-realization:
%%%%
%New
Define the operator $B_{\lambda, \, \nu}$ in the ${\mathcal X}$-realization by
\begin{equation}
(B_{\lambda, \, \nu}f)(x)=\int_{\mathcal X}[x,y]^{-\lambda-2,\nu}f(y) \, dy, 
\ x\in {\mathcal X}.
\end{equation}
It intertwines $R_{\lambda,\nu}$ with 
$R_{-\lambda-2,\nu}$. It acts from $\Gamma^{-}$ to $\Gamma^{+}$ by
\begin{equation}
(B_{\lambda ,\, \nu }f)(u)=2\int _{\Gamma^+}[u,v]^{-\lambda-2,\nu}f(v) \, dv, \ \ 
u\in\Gamma^{-},  \nonumber
\end{equation}
and similarly from $\Gamma^{-}$ to $\Gamma^{+}$. By (1.2) it can be written as
\begin{equation*}
\begin{split}
(B_{\lambda ,\, \nu }&f)(\xi_1,\eta_1)=\\
&(-1)^\nu 2^{-\lambda-1} \int _{\Gamma^+}\Big[N(\xi_1,\eta_2)N(\xi_2,\eta_1)\Big]^{-\lambda-2,\nu}
f(\xi_2,\eta_2) \, d\xi_2d\eta_2.  
\end{split}
\end{equation*}
It shows that 
\begin{equation}
B_{\lambda ,\, \nu}=(-1)^\nu 2^{-\lambda-1}A_{\sigma,\nu}\otimes A_{\sigma,\nu}, \ \ 
\sigma=\lambda/2.  \nonumber
\end{equation}
Therefore
\begin{equation}
B_{\lambda ,\, \nu}B_{-\lambda-2 ,\, \nu}=\big[\omega(\lambda/2,\nu)\big]^2
\cdot{\rm id}.
\end{equation}
Let us go back to the ${\mathcal X}$-realization and use 
%New
the
%%%%
{\it both} horospherical 
coordinate systems. Then 
\begin{equation}
(B_{\lambda ,\, \nu }f)(x)=(-1)^\nu 2^{-\lambda-2}
\int _{\mathcal X}\Big[\frac{N(\xi_1,\eta_2)N(\xi_2,\eta_1)}
{N(\xi_1,\eta_1)N(\xi_2,\eta_2)}\Big]^{-\lambda-2,\nu}
f(y) \, dy,
\end{equation}
where $x$ and $y$ have  coordinates $\xi_1,\eta_1$ and $\xi_2,\eta_2$ in 
%New
the
%%%%
horospherical 
coordinate systems corresponding to $\Gamma^{-}$ and $\Gamma^{+}$, respectively.

\section{The Berezin symbols and 
%New
the 
%%% you need "the" if it is the only one and well-known.
%%% if there are only 2 or 3 Berezin symbols, and also well-known,
%%% then you need "the" before them as well. 
Berezin transform}

The group $\widetilde G$ contains three subgroups isomorphic to $G$. The first one 
is the diagonal consisting of $(g,g), \ g\in G$. It preserves ${\mathcal X}$ under the action (1.1), hence ${\mathcal X}=G/H$. 
The measure $dx$ is invariant. The representation 
$R_{\lambda,\nu}$ is the representation by translations: 
\begin{equation}
R_{\lambda,\nu}(g,g)f(x)=f(g^{-1}xg). \nonumber
\end{equation}
Other two subgroups $G_1$ and $G_2$ consist of pairs $(g,e)$ and $(e,g)$, $g\in G$, 
respectively. By (3.2) we have on $\Gamma^{+}$:
\begin{equation}
(R_{\lambda,\nu}(e,g)f)(\xi,\eta)=
f(\xi\cdot g,\eta)
(\beta\xi+\delta)^{\lambda,\nu}.  \nonumber
\end{equation}
Therefore, in the horospherical coordinates on ${\mathcal X}$ corresponding to 
$\Gamma^{+}$, we have that
\begin{equation}
( R_{\lambda,\nu}(e,g)f)(\xi,\eta)
=\Big[ \frac{1}{N(\xi,\eta)}\Big]^{\lambda,\nu}
f(\xi\cdot g,\eta)N(\xi\cdot g,\eta)^{\lambda,\nu}
(\beta\xi+\delta)^{\lambda,\nu}. 
\end{equation}
This equation can be rewritten as follows. Denote
\begin{equation}
\Phi_{\lambda, \nu}(\xi,\eta)=N(\xi,\eta)^{\lambda,\nu}.  \nonumber
\end{equation}
It is the kernel of the intertwining operator for $G$ (see Sect. 1) and it is an 
analogue of the Berezin supercomplete system. Then (4.1) is
\begin{equation}
(R_{\lambda,\nu}(e,g)f)(\xi,\eta)=\frac{1}{\Phi_{\lambda,\nu} (\xi,\eta)}
(T_{\lambda/2,\nu}(g)\otimes 1)\big[f(\xi,\eta)\Phi_{\lambda,\nu} (\xi,\eta)\big]. 
 \nonumber
\end{equation}
Similarly, in the horospherical coordinates on ${\mathcal X}$ corresponding to 
$\Gamma^{-}$, we obtain that
\begin{equation}
(R_{\lambda,\nu}(e,g)f)(\xi,\eta)=\frac{1}{\Phi_{\lambda,\nu} (\xi,\eta)}
(1\otimes {\widehat T}_{\lambda/2,\nu}(g))
\big[f(\xi,\eta)\Phi_{\lambda,\nu}(\xi,\eta)\big]. 
 \nonumber
\end{equation}
(and similar formulas for $(g,e))$. Let us go from the group $G$ to its universal 
enveloping algebra ${\rm Env}(\goth g)$ and retain symbols for representations. 
Then 
%New
the 
%%%%
dependence of representations on $\nu$ disappears and we omit $\nu$ for them. 
Now take for $f$ the function $f_0$ equal to 1 identically, then for 
$X\in {\rm Env} (\goth g)$ we obtain that
\begin{align*}
(R_\lambda(0,X)f_0)(\xi,\eta)&=\frac{1}{\Phi_{\lambda,\nu} (\xi,\eta)}
(T_{\lambda/2}(X)\otimes 1)\Phi_{\lambda,\nu} (\xi,\eta),  \nonumber\\
(R_{-\lambda-2}(0,X)f_0)(\xi,\eta)&=\frac{1}{\Phi_{-\lambda-2,\nu} (\xi,\eta)}
(1\otimes {\widehat T}_{-\lambda/2-1}(X))\Phi_{-\lambda-2,\nu} (\xi,\eta).  \nonumber
\end{align*}
The right hand sides of these formulas are just 
%New
the 
%%%
{\it covariant} and {\it contravariant} 
 symbols of the operator $T_{\lambda/2}(X)$ in 
%New
the 
%%%%
polynomial quantization. 
\markboth{V.~Nitica, I.~Singer}{Polynomial quantization on para-hermitian spaces}

We can normalize the operator $B_{-\lambda-2,\nu}$ 
%Old
%such 
%New
so
that the normalized operator 
$Q_{\lambda,\nu}$ will satisfy the condition
\begin{equation}
Q_{\lambda,\nu}Q_{-\lambda-2,\nu}={\rm id}.  \nonumber
\end{equation}
Namely,
\begin{equation}
(Q_{\lambda ,\, \nu }f)(x)=c(\lambda,\nu)
\int _{\mathcal X}\Big[\frac{(1-\xi_1\eta_2)(1-\xi_2\eta_1)}
{(1-\xi_1\eta_1)(1-\xi_2\eta_2)}\Big]^{\lambda,\nu}
f(y) \, dy,
\end{equation}
where $x$ and $y$ have coordinates $\xi_1,\eta_1$ and $\xi_2,\eta_2$ as in (3.5) 
and  
\begin{equation}
c(\lambda,\nu)^{-1}=2\omega(\lambda/2,\nu),
\end{equation}
see (2.1), (3.4), (3.5). The kernel in (4.2) (with the factor $c(\lambda,\nu)$) 
is nothing but the Berezin kernel. Therefore, the operator $Q_{\lambda,\nu}$ is 
the Berezin transform. It transfers contravariant symbols to covariant ones. 
%New
Note that  if we 
want to write 
the Berezin transform using only {\it one} coordinate system, then 
we will have to change the operator (3.3), namely, we will have to write $[x,yJ]$ instead of $[x,y]$.
%Old
%If we would write the Berezin transform using only {\it single} coordinate system, then we
%should change the operator (3.3): write $[x,yJ]$ instead of $[x,y]$.

\setcounter{section}{0}
\setcounter{footnote}{0}
\nachaloe{V.~Nitica and I.~Singer}{The structure of max-plus hyperplanes}%
{Partially supported by NSF
grant DMS-0500832 and by grant nr. 2-CEx06-11-34/2006.}
\markboth{V.~Nitica and I.~Singer}{The structure of max-plus hyperplanes}
\label{sin-abs}
A max-plus hyperplane (briefly, a hyperplane) is the set of all points 
$x=(x_{1},...,x_{n})$ in $\R_{\max }^{n}$ 
satisfying an equation of the form 
$$
a_{1}x_{1}\oplus ...\oplus a_{n}x_{n}\oplus a_{n+1}=b_{1}x_{1}\oplus
...\oplus b_{n}x_{n}\oplus b_{n+1},
$$ 
that is, 
$$\max
(a_{1}+x_{1},...,a_{n}+x_{n},a_{n+1})=\max
(b_{1}+x_{1},...,b_{n}+x_{n},b_{n+1}),
$$ 
with $a_{i},b_{i}\in \R_{\max
}\;(i=1,...n+1),$ where each side contains 
at least one term, and where $a_{i}\neq b_{i}$ 
for at least one index $i.$ We show that the complements of
(max-plus) semispaces at finite points $z\in \R^{n}$ are 
``building blocks''
for the hyperplanes in $\R_{\max }^{n}$ (recall that a 
semispace at $z$ is a
maximal --with respect to inclusion-- max-plus convex subset of 
$\R_{\max
}^{n}\backslash \{z\}).$ Namely, observing that, up to 
a permutation of
indices, we may write the equation of any 
hyperplane $H$ in one of the
following two forms:
\begin{align*}
& a_{1}x_{1}\oplus ...\oplus a_{p}x_{p}\oplus a_{p+1}x_{p+1}\oplus ...\oplus
a_{q}x_{q} \\
& =a_{1}x_{1}\oplus ...\oplus a_{p}x_{p}\oplus a_{q+1}x_{q+1}\oplus
...\oplus a_{m}x_{m}\oplus a_{n+1},
\end{align*}
where $0\leq p\leq q\leq m\leq n$ and all $a_{i}\;(i=1,...,m,n+1)$ are
finite, or,
\begin{align*}
& a_{1}x_{1}\oplus ...\oplus a_{p}x_{p}\oplus a_{p+1}x_{p+1}\oplus ...\oplus
a_{q}x_{q}\oplus a_{n+1} \\
& =a_{1}x_{1}\oplus ...\oplus a_{p}x_{p}\oplus a_{q+1}x_{q+1}\oplus
...\oplus a_{m}x_{m}\oplus a_{n+1},
\end{align*}
where $0\leq p\leq q\leq m\leq n$, and all $a_{i}\;(i=1,...,m)$ are finite
(and $a_{n+1}$ is either finite or $-\infty ),$ we give a formula that
expresses a nondegenerate strictly affine hyperplane (i.e., with $m=n$ and $%
a_{n+1}>-\infty )$ as a union of complements of semispaces at a point $z\in
\R^{n},$ called the ``center'' of $H,$ with the boundary of a union of
complements of other semispaces at $z.$ Using this formula, we obtain
characterizations of nondegenerate strictly affine hyperplanes with empty
interior. We give a description of the boundary of a nondegenerate strictly
affine hyperplane with the aid of complements of semispaces at its center,
and we characterize the cases in which the boundary bd $H$ of a
nondegenerate strictly affine hyperplane $H$ is also a hyperplane. Next, we
give the relations between nondegenerate strictly affine hyperplanes $H$,
their centers $z$, and their coefficients $a_{i}.$ In the converse direction
we show that any union of complements of semispaces at a point $z\in \R^{n}$
with the boundary of any union of complements of some other semispaces at
that point $z$, is a nondegenerate strictly affine hyperplane. We obtain a
formula for the total number of strictly affine hyperplanes. We give
complete lists of all strictly affine hyperplanes for the cases $n=1$ and $%
n=2$. We show that each linear hyperplane $H$ in $\R_{\max }^{n}$ (i.
e., with $a_{n+1}=-\infty )$ can be decomposed as the union of four
parts, where each part is easy to describe in terms of complements
of semispaces, some of them in a lower dimensional space.

The paper in extenso will appear in Linear Algebra and its
Applications.
\markboth{E.~Pap, M.~\v{S}trboja}{The structure of max-plus hyperplanes}
\setcounter{footnote}{0}
\setcounter{section}{0}
\nachaloe{E.~Pap and M.~\v{S}trboja}{Image processing based on a partial differential equation
satisfying the pseudo-linear superposition principle}%
{Partially supported by the
Project MNZ\v ZSS $144012,$ grant of MTA HTMT, French-Serbian
project "Pavle Savi\'c", and by the project ''Mathematical Models
for Decision Making under Uncertain Conditions and Their
Applications'' of Academy of Sciences and Arts of Vojvodina
supported by Provincial Secretariat for Science and Technological
Development of Vojvodina.}
\markboth{E.~Pap, M.~\v{S}trboja}{Image processing based on PDE}
\label{str-abs}
We consider a general form of PDE-based methods for image restoration, and
give a short overview of the underlying models. In these models, the
original image is transformed through a process that can be represented by a
second-order partial differential equation. Typically, this role is played
by some nonlinear generalization of the heat equation, and it is possible to
analyse the solutions from the viewpoint of pseudo-linear (idempotent)
analysis. Our main result is that the generalization of the heat equation
proposed by Perona and Malik satisfies the pseudo-linear superposition
principle.

\section{Introduction}

The approach based on partial differential equations is well-known in image
processing (\cite{image,catte,peron}). In this approach, a restored image
can be seen as a version of the initial image at a special scale. Image $u$
is an instance of an evolution process, denoted by $u\left( t,\cdot \right) $%
. The original image is taken at time $t=0,$ $u\left( 0,\cdot \right)
=u_{0}\left( \cdot \right) .$ The original image is then transformed, and
this process can be described by the equation $\frac{ \partial u}{\partial t}%
\left( t,x\right) +F\left( x,u(t,x),\nabla u\left( t,x\right) ,\nabla
^{2}u\left( t,x\right) \right) =0$, where $x\in\Omega .$ Some possibilities
for $F$ to restore an image are considered in \cite{image}.

Pseudo-linear superposition principle means the following. Instead of the
field of real numbers, think of a semiring defined on a real interval $\left[
a,b\right] \subset $ $\left[ -\infty ,\infty \right]$. This is a structure
equipped with pseudo-addition $\oplus$, which is typically idempotent ($%
x\oplus x=x$), and with pseudo-multiplication $\odot$. The pseudo-linear
superposition principle says that some nonlinear equations (ODE, PDE,
difference equations, etc.) turn out to be linear over such structures,
meaning that if $u_1$ and $u_2$ are two solutions of the considered
equation, then $a_1\odot u_1\oplus a_2\odot u_2$ is also a solution for any
constants $a_1$ and $a_2$ from $\left[ a,b\right]$.

By pseudo-analysis we mean analysis over such semirings, in the framework of
\cite{90,g,hand,null,dec,gen,papral}. One of its key ideas is the
pseudo-linear superposition principle stated above. This (pseudo-) linear
intuition leads to the concepts of $\oplus $-measure, pseudo-integral,
pseudo-convolution, pseudo-Laplace transform, etc.

Similar ideas were developed independently by Maslov and his collaborators
in the framework of idempotent analysis and idempotent mathematics, with
some applications \cite{KoMa89,kolmasl,Lit05,maslSomb}. In particular,
idempotent analysis is fundamental for the theory of weak solutions to
Hamilton-Jacobi equations with non-smooth Hamiltonians, see \cite%
{KoMa89,kolmasl,maslSomb} and also \cite{gen,papral} (which use the language
of pseudo-analysis). In some cases, this theory enables one to obtain exact
solutions in the similar form as for the linear equations. Some further
developments relate more general pseudo-operations with applications to
nonlinear partial differential equations, see \cite{pv}. Recently, these
applications have become important in the field of image processing \cite%
{gen,papral}.

Our report is organized as follows. In Sect. \ref{s2} we consider a
general form of PDE for image restoration. The starting PDE in image
restoration is the heat equation. Because of its oversmoothing
property (edges get smeared), it is necessary to introduce some
nonlinearity. We consider then the following model
(\cite{image,peron})\begin{equation*} \frac{\partial u}{\partial
t}=\text{div}\left( c\left( \left\vert \nabla u\right\vert
^{2}\right) \nabla u\right) ,  \label{c}
\end{equation*}where we choose the function $c$ such that the equation remains to
be of the parabolic type. We take $c(s)\approx 1/\sqrt{s}$ as
$s\rightarrow \infty ,$ because we want to preserve the
discontinuities \cite{image}. Because of this behavior, it is not
possible to apply general results from parabolic equations theory.
An appropriate framework to study this equation is nonlinear semigroup theory (\cite%
{image,51,62}). In Sect. \ref{s3} we show that Perona and Malik
equation satisfies the pseudo-linear superposition principle.

\section{PDE-based method in image processing}

\label{s2}

PDE-methods for restoration can be written in the following general
form: \begin{equation*} \left\{
\begin{array}{l}
\frac{\partial u}{\partial t}\left( t,x\right) +F\left( x,u(t,x),\nabla
u\left( t,x\right) ,\nabla ^{2}u\left( t,x\right) \right) =0\text{ in }%
\left( 0,T\right) \times \Omega , \\
\frac{\partial u}{\partial N}\left( t,x\right) =0\text{ on }\left(
0,T\right) \times \partial \Omega \text{ (Neumann boundary condition),} \\
\text{ }u\left( 0,x\right) =u_{0}\left( x\right) \text{ (initial condition),}%
\end{array}%
\right\}  \label{opsti}
\end{equation*}
where $u\left( t,x\right) $ is the restored version of the initial degraded
image $u_{0}\left( x\right) $. The idea is to construct a family of
functions $\left\{ u\left( t,x\right) \right\} _{t>0}$ representing
successive versions of $u_{0}\left( x\right) $. As $t$ increases, the image $%
u\left( t,x\right) $ becomes more and more simplified. We would like to
attain two goals. The first is that $u\left( t,x\right) $ should represent a
smooth version of $u_{0}\left( x\right) $, where the noise has been removed.
The second, $u\left(t,x\right)$ should be able to preserve some features
such as edges, corners, which may be viewed as singularities.

The heat equation is the basic PDE for image restoration:
\begin{equation*} \left\{
\begin{array}{c}
\frac{\partial u}{\partial t}\left( t,x\right) -\Delta u\left( t,x\right) =0,%
\text{ }t\geq 0,\text{ }x\in \mathbb{R}^{2},\text{ \ } \\
u\left( 0,x\right) =u_{0}\left( x\right) .%
\end{array}%
\right.  \label{topl}
\end{equation*}
The heat equation has been successfully applied in image processing
but it has some drawback. It is too smoothing and because of that
edges can be lost or severely blurred. In \cite{image} authors
consider models that are generalizations of the heat
equation. Suppose that the domain image is a bounded open set $\Omega $ of $%
\mathbb{R}^{2}$. The following equation was initially proposed by
Perona and Malik \cite{peron}: \begin{equation} \left\{
\begin{array}{c}
\frac{\partial u}{\partial t}=\text{div}\left( c\left( \left\vert \nabla
u\right\vert ^{2}\right) \nabla u\right) \text{ \ \ in }\left( 0,T\right)
\times \Omega , \\
\frac{\partial u}{\partial N}=0\text{ on }\left( 0,T\right) \times \partial
\Omega \text{, \ \ \ \ \ \ \ \ \ \ \ \ \ \ \ \ \ \ \ \ \ \ } \\
u\left( 0,x\right) =u_{0}\left( x\right) \text{ in}\ \Omega \ \ \ \ \ \ \ \
\ \ \ \ \ \ \ \ \ \ \ \ \ \ \ \ \
\end{array}%
\right.  \label{dif}
\end{equation}%
where $c:\left[ 0,\infty \right) \rightarrow \left( 0,\infty \right) .$ If
we choose $c\equiv 1$, then it is reduced to the heat equation. If we assume
that $c\left( s\right) $ is a decreasing function satisfying $c\left(
0\right) =1$ and $\underset{s\rightarrow \infty }{\lim }c\left( s\right) =0$%
, then inside the regions where the magnitude of the gradient of $u$ is
weak, equation (\ref{dif}) acts like the heat equation and the edges are
preserved.

For each point $x$ where $\left\vert \nabla u\right\vert \neq 0$ we can
define the vectors $\mathbf{N}=\frac{\nabla u}{\left\vert \nabla u\right\vert }$ and $%
\mathbf{T}$ with $\mathbf{T}\cdot \mathbf{N}=0,$ $\left\vert
\mathbf{T}\right\vert =1.$ For the first and second
partial derivatives of $u$ we use the usual notation $u_{x_{1}},$ $%
u_{x_{2}}, $ $u_{x_{1}x_{1},\text{...}}$ We denote by
$u_{\mathbf{NN}}$ and $u_{\mathbf{TT}}$ the second derivatives of
$u$ in the $\mathbf{N}$-direction and $\mathbf{T}$-direction,
respectively:
\begin{eqnarray*}
u_{\mathbf{NN}} &=&\mathbf{N}^{t}\text{ }\nabla ^{2}u\text{
}\mathbf{N}=\frac{1}{\left\vert \nabla u\right\vert
^{2}}(u_{x}^{2}u_{xx}+u_{y}^{2}u_{yy}+2u_{x}u_{y}u_{xy}),\\
u_{\mathbf{TT}} &=&\mathbf{T}^{t}\text{ }\nabla ^{2}u\text{
}\mathbf{T}=\frac{1}{\left\vert \nabla u\right\vert ^{2}}\left(
u_{x}^{2}u_{yy}+u_{y}^{2}u_{xx}-2u_{x}u_{y}u_{xy}\right) .
\end{eqnarray*}

The first equation in (\ref{dif}) can be written as \begin{equation}
\frac{\partial u}{\partial t}\left( t,x\right) =c\left( \left\vert
\nabla u\left( t,x\right) \right\vert ^{2}\right)
u_{\mathbf{TT}}+b\left( \left\vert \nabla u\left( t,x\right)
\right\vert ^{2}\right) u_{\mathbf{NN}}, \label{nt}
\end{equation}%
where $b(s)=c(s)+2sc^{\prime }(s)$. Therefore, (\ref{nt}) is a sum
of a diffusion in the $\mathbf{T}$-direction and a diffusion in the
$\mathbf{N}$-direction. The functions $c$ and $b$ act as weighting
coefficients. Since $\mathbf{N}$ is normal to the edges, it would be
preferable to smooth more in the tangential direction $\mathbf{T}$
than in the normal direction. Because of that we impose
\begin{equation}
\underset{s\rightarrow \infty }{\lim }\frac{b(s)}{c(s)}=0\text{ \ or \ }%
\underset{s\rightarrow \infty }{\lim }\frac{sc^{\prime }(s)}{c(s)}=-\frac{1}{%
2}  \label{ponasanje}
\end{equation}%
If $c(s)>0$ with power growth, then (\ref{ponasanje}) implies that $%
c(s)\approx 1/\sqrt{s}$ as $s\rightarrow \infty $. The equation (\ref{dif})
is parabolic if $b(s)>0$.

The assumptions imposed on $c\left( s\right) $ are \begin{equation}
\left\{
\begin{array}{c}
c:\left[ 0,\infty \right) \rightarrow \left( 0,\infty \right) \text{
decreasing, \ } \\
c(0)=1,\text{ }c(s)\approx \frac{1}{\sqrt{s}}\text{ as }s\rightarrow \infty ,
\\
b(s)=c(s)+2sc^{\prime }(s)>0.\text{ \ \ \ \ \ \ \ \ \ \ }%
\end{array}%
\right.  \label{sum}
\end{equation}
Consider $c(s)=\frac{1}{\sqrt{1+s}}$, an often used function satisfying (\ref%
{sum}). Because of the behavior $c(s)\approx 1/\sqrt{s}$ as $s\rightarrow
\infty $, it is not possible to apply general results from parabolic
equations theory. An appropriate framework to study this equation is
nonlinear semigroup theory (see \cite{image,51,62}).

\section{Pseudo-linear superposition principle for Perona and Malik equation}

\label{s3} Let $\left[ a,b\right] $ be a closed (in some cases semiclosed)
subinterval of $\left[ -\infty ,+\infty \right] $. We consider here the
total order $\leq $ on $\left[ a,b\right] $. The operation $\oplus $
(pseudo-addition) is a commutative, non-decreasing, associative function $%
\oplus :\left[ a,b\right] \times \left[ a,b\right] \rightarrow \left[ a,b%
\right] $ with a zero (neutral) element denoted by $\mathbf{0}$. 
Denote $%
\left[ a,b\right] _{+}=\{x:x\in \left[ a,b\right] ,x\geq \mathbf{0}\}.$ The
operation $\odot $ (pseudo-multiplication) is a function $\odot :\left[ a,b%
\right] \times \left[ a,b\right] \rightarrow \left[ a,b\right] $ which is
commutative, positively non-decreasing, i.e., $x\leq y$ implies $x\odot
z\leq y\odot z,z\in \left[ a,b\right] _{+},$associative and for which there
exist a unit element $\mathbf{1}\in \left[ a,b\right] ,$ i.e., for each $%
x\in \left[ a,b\right] ,1\odot x=x.$ We assume $\mathbf{0}\odot
x=\mathbf{0}$ and that $\odot $ is distributive over $\oplus ,$
i.e.,$$ x\odot (y\oplus z)=(x\odot y)\oplus (x\odot z)$$ The
structure $\left( \left[ a,b\right] ,\oplus ,\odot \right) $ is
called a \textit{semiring} (see \cite{kuich,null}). In this paper we
shall
consider only the min-plus (or tropical) semiring. It is defined on the interval $%
\left( -\infty ,+\infty \right] $ and has the following continuous
operations: $x\oplus y=\min \left\{ x,y\right\} ,\text{ \ \ }x\odot y=x+y.$
Note that the pseudo-addition is idempotent, while the pseudo-multiplication
is not. We have $\mathbf{0}=-\infty $ and $\mathbf{1}=0$.

We show that the pseudo-linear superposition principle 
holds
for Perona and Malik equation.

\begin{theorem}
If $u_{1}=$ $u_{1}\left( t,x\right) $ and $u_{2}=$ $u_{2}\left(
t,x\right) $ are solutions of the equation \begin{equation}
\frac{\partial u}{\partial t}-\text{div}\left( c\left( \left\vert
\nabla u\right\vert ^{2}\right) \nabla u\right) =0, \label{strboja1}
\end{equation}
then $\left( \lambda _{1}\odot u_{1}\right) \oplus \left( \lambda _{2}\odot
u_{2}\right) $ is also a solution of (\ref{strboja1}) on the set
\begin{equation*}
D=\{\left( t,x\right) |t\in \left(0,T\right),x\in \mathbb{R}^{2},
u_{1}\left( t,x\right) \neq u_{2}\left( t,x\right) \},
\end{equation*}
with respect to the operations $\oplus =\min $ and $\odot =+$.
\end{theorem}

The obtained results will serve for further investigation of the
weak solutions of the equation (\ref{strboja1}) in the sense of Maslov \cite%
{KoMa89,kolmasl} and Gondran \cite{Gond,GondMin}, as well as some important
applications.

\setcounter{footnote}{0}
\setcounter{section}{0}
\nachalo{Alexander Rashkovskii}{Tropical analysis on 
plurisubharmonic singularities}
\label{ras-abs}

\section{Plurisubharmonic singularities} 

Recall that an upper
semicontinuous, real-valued function on an open set in $\Bbb C^n$ is
called {\it plurisubharmonic} (psh) if its restriction to every
complex line is a subharmonic function. A basic example is $\log|f|$
for an analytic function $f$. Moreover, by Bremermann's theorem
\cite{B}, every psh function $u$ can be written as
$$u(z)=\limsup_{y\to z}\limsup_{m\to\infty}\frac1m\log|f_m(y)|.$$

Let $\cO_0$ denote the ring of germs of analytic functions $f$
at $0\in{\Bbb C}^n$, and let $\mathfrak{m}_0=\{f\in \cO_0 :
f(0)=0\}$ be its maximal ideal. The log-transformation
$f\mapsto\log|f|$ maps $\cO_0$ into the collection of germs of
psh functions at $0$.

We will say that a psh germ $u$ has singularity at $0$ if
$u(0)=-\infty$. For functions $u=\log|f|$ this means
$f\in\mathfrak{m}_0$; asymptotic behaviour of arbitrary psh
functions can be much more complicated. By $PSHG_0$ we denote the
collections of all psh germs singular at $0$.
\markboth{Alexander Rashkovskii}{Tropical analysis on 
plurisubharmonic singularities}

The operations on $\cO_0$ induce a natural tropical structure
on $PSHG_0$ with the addition $u\oplus v:=\max\{u,v\}$ (which is
based on Maslov's dequantization: $f+g\mapsto
\frac1N\log|f^N+g^N|\to \log|f|\oplus\log|g|$ as $N\to\infty$) and
multiplication $u\otimes v:=u+v$ (simply by $fg\mapsto \log|f
g|=\log|f|\otimes\log|g|$). Thus $PSHG_0$ becomes a tropical
semiring, closed under (usual) multiplication by positive constants.

A partial order on $PSHG_0$ is given as follows: $u\preceq v$ if
$u(z)\le v(z)+O(1)$ as $z\to 0$, which leads to the equivalence
relation $u\sim v$ if $u(z)=v(z)+O(1)$. The equivalence class ${\rm
cl}(u)$ of $u$ is called the {\it plurisubharmonic singularity} of
the germ $u$. The collection of psh singularities
$PSHS_0=PSHG_0/\sim$ has the same tropical structure
$\{\oplus,\otimes\}$ and the partial order: ${\rm cl}(u)\le {\rm
cl}(v)$ if $u\preceq v$. It is endowed with the following topology:
${\rm cl}(u_j)\to {\rm cl}(u)$ if there exists a neighbourhood
$\omega$ of $0$ and psh functions $v_j\in {\rm cl}(u_j)$, $v\in {\rm
cl}(u)$ in $\omega$ such that $v_j\to v$ in $L^1(\omega)$.

By abusing the notation, we will right occasionally $u$ for ${\rm
cl}(u)$.

\section{Characteristics of singularities} 

The main characteristic
of an analytic germ $f\in\mathfrak{m}_0$ is its multiplicity
(vanishing order) $m_f$: if $f=\sum P_j$ is the Taylor expansion of
$f$ in homogeneous polynomials, $P_j(tz)=t^jP(z)$, then
$m_f=\min\{j: P_j\not\equiv 0\}$.

The basic characteristic of singularity of $u\in PSHG_0$ is its {\it
Lelong number}
$$ \nu(u)=\lim_{t\to -\infty}\frac1t M(u,t)=
\liminf_{z\to 0}\frac{u(z)}{\log|z|}=dd^cu\wedge
(dd^c\log|z|)^{n-1}(0);$$ here $M(u,t)$ is the mean value of $u$
over the sphere $\{|z|=e^t\}$, $d=\partial + \bar\partial$, $d^c= (
\partial -\bar\partial)/2\pi i$. If $f\in\mathfrak{m}_0$, then
$\nu(\log|f |)=m_f$. This characteristic of singularity gives an
important information on the asymptotic behaviour of $u$ at $0$:
$u(z)\le\nu(u)\log|z|+O(1)$.

Since $\nu(v)=\nu(u)$ for all $v\in {\rm cl}(u)$, Lelong number can
be considered as a functional on $PSHS_0$ with values in the
tropical semiring $\Bbb R_{+}(\min,+)$ of non-negative real numbers
with the operations $x\bar\oplus y=\min\{x,y\}$ and $x\otimes
y=x+y$. As such, it is
\par\noindent (i) positive homogeneous: $\nu(cu)=c\,\nu(u)$ for all
$c>0$, \par\noindent (ii) additive: $\nu(u\oplus v)=
\nu(u)\bar\oplus \nu(v)$, \par\noindent (iii) multiplicative:
$\nu(u\otimes v)= \nu(u)\otimes \nu(v)$, and \par\noindent (iv)
upper semicontinuous: $\nu(u)\ge\limsup \nu(u_j)$ if $u_j\to u$.

Lelong numbers are independent of the choice of coordinates. Let us
now fix a coordinate system (centered at $0$). The {\it directional
Lelong number} of $u$ in the direction $a\in\Bbb R_+^n$ (introduced
by C.~Kiselman \cite{Kis2}) is \begin{equation}\label{eq:dir} \nu(u,
a)=\lim_{t\to -\infty}\frac1t M(u,ta)= \liminf_{z\to
0}\frac{u(z)}{\phi_a(z)},\end{equation} where $M(u,ta)$ is the mean
value of $u$ over the distinguished boundary of the polydisk
$\{|z_k|<\exp(ta_k)\}$ and $\phi_a(z)=\oplus_k \,a_k^{-1}\log|z_k|$.
It has the same properties (i)--(iv), and the collection
$\{\nu(u,a)\}_{a}$  gives a refined information on the singularity
$u$. In particular, $\nu(u)=\nu(u,(1,\ldots,1))$.

A general notion of Lelong number with respect to a plurisubharmonic
weight was introduced by J.-P.~Demailly \cite{D1}. Let $\varphi\in
PSHG_0$ be continuous and locally bounded outside $0$. Then the
mixed Monge--Amp\`ere current $dd^c u\wedge (dd^c\varphi)^{n-1}$ is
well defined for any psh function $u$ and is equivalent to a
positive Borel measure. Its mass at $0$, $\nu(u,\varphi)= dd^c
u\wedge (dd^c\varphi)^{n-1}(\{0\})$, is called {\it the generalized
Lelong number}, or {\it the Lelong--Demailly number}, of $u$ with
respect to the weight $\varphi$. Since it is constant on ${\rm
cl}(u)$, we have a different kind of functional on $PSHS_0$. It
still has the above properties (i), (iii), and (iv), however in
general is only subadditive: $\nu(u\oplus v, \varphi)\le \nu(u,
\varphi)\bar\oplus \nu(v, \varphi)$.

Note that $\nu(u,a)=a_1\ldots a_n\,\nu(u,\phi_a)$.

\section{Additive functionals} 

Another generalization of the notion
of Lelong number was introduced in \cite{R7}. Let $\varphi\in
PSHG_0$ be locally bounded and {\it maximal} outside $0$ (that is,
satisfies $(dd^c\varphi)^n=0$ on a punctured neighbourhood of $0$);
the collection of all such germs (weights) will be denoted by
$MW_0$. The {\it type of $u\in PSHS_0$ relative to} $\varphi\in
MW_0$,
$$\sigma(u,\varphi)=\liminf_{z\to 0}\frac{u(z)}{\varphi(z)},$$
gives the bound $u\le\sigma(u,\varphi)\varphi$.

This functional is positive homogeneous, additive,
supermultiplicative, and upper semicontinuous. Actually, relative
types give a general form for all "reasonable" additive functionals
on $PSHS_0$:

\begin{theorem} 
{\rm (\cite{R7})}  
Let a functional $\sigma:\:
PSHS_0\to [0,\infty]$ be such that
\begin{enumerate} \item[1)] $\sigma(cu)=c\,\sigma(u)$ for all $c>0$;
\item[2)] $\sigma(\oplus u_k)=\bar\oplus \sigma(u_k)$, $k=1,2$;
\item[3)] if $u_j\to u$, then $\limsup\,
\sigma(u_j)\le\sigma(u)$;
\item[4)] $\sigma(\log|z|)>0$;
\item[5)] $\sigma(u)< \infty$ if
$u\not\equiv-\infty$.
\end{enumerate}
Then there exists a weight $\varphi\in MW_0$ such that
$\sigma(u)=\sigma(u,\varphi)$ for every $u\in PSHS_0$. The
representation is essentially unique: if two maximal weights
$\varphi$ and $\psi$ represent $\sigma$, then ${\rm
cl}(\varphi)={\rm cl}(\psi)$.
\end{theorem}

\section{Relative types and valuations} 

Recall that a {\it
valuation} on the analytic ring $\cO_0$  is a nonconstant
function $\mu:\: \cO_0\to [0,+\infty]$ such that
$$\mu(f_1f_2)=\mu(f_1)+\mu(f_2),\quad
\mu(f_1+f_2)\ge\min\,\{\mu(f_1),\mu(f_2)\}, \quad\mu(1)=0;$$ a
valuation $\mu$ is {\it centered} if $\mu(f)>0$ for every
$f\in\mathfrak{m}_0$, and {\it normalized} if $\min\,
\{\mu(f):f\in\mathfrak{m}_0\}=1$. Every weight $\varphi\in MW_0$
generates a functional $\sigma_\varphi$ on $\cO_0$,
$\sigma_\varphi(f)=\sigma(\log|f|,\varphi)$, with the properties
\begin{align*}
&\sigma_\varphi(f_1f_2)\ge \sigma_\varphi(f_1)+\sigma_\varphi(f_2), \\
&\sigma_\varphi(f_1+f_2)\ge\min\,\{\sigma_\varphi(f_1),\sigma_\varphi(f_2)\},\quad\sigma_\varphi(1)=0.
\end{align*}

It is a valuation, provided
$\sigma(u,\varphi)$ is tropically multiplicative; $\sigma_\varphi$
is centered iff $\sigma(\log|z|,\varphi)>0$, and normalized iff
$\sigma(\log|z|,\varphi)=1$.

One can thus consider linear (both additive and multiplicative)
functionals on $PSHS_0$ as tropicalizations of certain valuations on
$\cO_0$.

For example, the (usual) Lelong number is the tropicalization of the
multiplicity valuation $m_f$. The types relative to the directional
weights $\phi_{a}$ are multiplicative functionals on $PSHS_0$, and
$\sigma_{\phi_{a}}$ are monomial valuations on $\cO_0$; they
are normalized if $\min_ka_k=1$. It was shown in \cite{FaJ1} that an
important class of valuations (quasi-monomial valuations, or
Abhyankar valuations of rank $1$) can be realized as
$\sigma_\varphi$ with certain weights $\varphi\in MW_0$; when $n=2$,
all other centered valuations are limits of increasing sequences of
the quasi-monomial ones \cite{FaJ}.

\section{Local indicators as Maslov's dequantizations} 

Consideration
of psh germs is the first step of Maslov's dequantization of
analytic functions $\cO_0\ni f\mapsto\log|f|\in PSHG_0$. One
can perform the next step -- namely, passage to the logarithmic
scale in the arguments $z$.

For a fixed coordinate system at $0$, let $\nu(u,a)$ be the
directional Lelong numbers of $u\in PSHS_0$ in the directions $a\in
\Bbb R_+^n$ (\ref{eq:dir}). Then the function
$\psi_u(t)=-\nu(u,-t)$, $t\in\Bbb R_-^n$, is convex and increasing
in each $t_k$, so $\psi_u(\log|z_1|,\ldots,|z_n|)$ can be extended
(in a unique way) to a function $\Psi_u(z)$ plurisubharmonic in the
unit polydisk $\Bbb D^n=\{z\in\Bbb C^n:\:|z_k|<1,\:1\le k\le n\}$.
This function is called the {\it local indicator} of $u$ at $0$
\cite{LeR}. It is easy to see that it has the homogeneity property
\begin{equation}\label{eq:indic} \Psi_u(z_1,\ldots,z_n)=\Psi_u(|z_1|,
\ldots,|z_n|)= c^{-1}\Psi_u(|z_1|^c,\ldots,|z_n|^c) \quad \forall
c>0.\end{equation}

It was shown in \cite{R} that $\Psi_u(z)$ can be represented as the
(unique) weak limit of the functions $m^{-1}u(z_1^m,\ldots, z_n^m)$
as $m\to\infty$, so the indicator can be viewed as the tangent (in
the logarithmic coordinates) for the function $u$ at $0$. This means
that for $u=\log|f|$, $f\in \mathfrak{m}_0$, the sublinear function
$\psi_u(t)$ on $\Bbb R_-^n$ is just a Maslov's dequantization of
$f$.

The indicator is a psh characteristic of asymptotic behaviour near
$0$. Namely, if $u$ is psh in the unit polydisk $\Bbb D^n$, then
$u(z)\le \Psi_u(z)+\sup_{\Bbb D}u$. When $u$ has isolated
singularity at $0$, this implies the following relation between the
residual Monge-Amp\`ere masses: $(dd^cu)^n(0)\ge (dd^c\Psi_u)^n(0)$.

Since $\Psi_u$ is much simpler than the original function $u$, one
can compute explicitly the value of its residual mass. The first
equation in (\ref{eq:indic}) suggests us to pass from
plurisubharmonic functions to convex ones and from the complex
Monge-Amp\`ere operator to the real one, while the second equation
allows us to calculate the real Monge-Amp\`ere measure in terms of
volumes of gradient images. Denote
$$
\Theta_{u,x}=\{b\in\Bbb R_+^n:\: \sup_{\sum a_k=1}[\nu(u,a)-\langle
b,a \rangle]\ge 0\}.
$$
The convex image $\psi_u(t)$, $t\in\Bbb R_-^n$, of the indicator is
just the support function to the convex set $\Gamma_u=\Bbb
R_+^n\setminus\Theta_{u,x}$: $\psi_u(t)=\sup\,\{\langle t, a\rangle:
\: a\in \Gamma_u\}$. This gives

\begin{theorem}
{\rm (\cite{R})}
The residual Monge-Amp\`ere mass of
$u\in PSHG_0$ with isolated singularity at $0$ has the lower bound
$$(dd^cu)^n(0)\ge (dd^c\Psi_u)^n(0))=
n!\,{\rm Vol}(\Theta_{u,x}).$$
\end{theorem}

If $F=(f_1,\ldots, f_n)$ is a holomorphic mapping with isolated zero
at $0$, then its multiplicity at $0$ equals $(dd^c\log|F|)^n(0)$ and
the set $\Gamma_{\log|F|}$ is the convex hull of the union of the
Newton polyhedra ${\rm conv} \{\alpha+\Bbb R_+^n: \:
D^{(\alpha)}f_j(0)\neq 0\}$ of $f_j$ at $0$, $1\le j\le n$. In this
case, Theorem~2 gives us Kushnirenko's theorem on multiplicity of
holomorphic mappings \cite{Ku1}.

%\bigskip
The results on local indicators have global counterparts concerning
psh functions of logarithmic growth in $\Bbb C^n$ (i.e., $u(z)\le
A\log(1+|z|)+B$ everywhere in $\Bbb C^n$, a basic example being
logarithm of modulus of a polynomial), see \cite{R1} and \cite{R2}.
Similar notions concerning Maslov's dequantization in $\Bbb C^n$ and
generalized Newton polytops were also introduced and studied in
\cite{LS}.

\markboth{Serge\u{\i} Sergeev}{Tropical analysis on plurisubharmonic singularities}

\newpage
\setcounter{footnote}{0}
\setcounter{section}{0}
\nachaloe{Serge{\u{\i}} Sergeev}{Minimal elements and 
cellular closures over the max-plus semiring}%
{Supported by
the RFBR grant 05-01-00824 and the joint
RFBR/CNRS grant 05-01-02807.}
%\markboth{Serge{\u{\i}} Sergeev}{Tropical analysis on plurisubharmonic singularities}
\label{ser-abs} 
\markboth{Serge{\u{\i}} Sergeev}{Minimal elements and cellular closures}
This report is based on the publications
\cite{BSS-07} (part 1) and \cite{Ser-07<DC>} (part 2). In part 1,
I outline some simple consequences of the observation that extremals
are minimal elements with respect to the certain preorder relation.
Part 2 is occupied with some extensions of algebraic closure operation,
which arise from the cellular decomposition considered in \cite{DS-04}. 
The common feature of these works is a bit of interplay between 
max-algebra \cite{But-03} and tropical convexity \cite{DS-04},
\cite{Jos-05}. All results are obtained in the setting of $\Rmaxmn$,
the $n$-dimensional free semimodule over
the semiring $\Rmaxm=(\R_+,\oplus=\max,\odot=*)$.

\section{Extremals as minima} 
An element $u$ of a (sub)semimodule $K\subseteq\Rmaxmn$ is an {\em extremal},
if $u=x\oplus y$, $x,y\in K$ implies that $u=x$ or $u=y$.
The preorder relation $\leq_j$ is defined by
\begin{equation}
\label{leqj}
u\leq_j v\Leftrightarrow u_j\ne 0,\,v_j\ne 0,\,u/u_j\leq v/v_j.
\end{equation}
The role of $\leq_j$ is explained in the following.
\begin{proposition}
\label{spanit}
{\rm \cite{Jos-05,BSS-07}}
The following are equivalent:
\begin{itemize}
\item[(1)] $y$ is a (max-)linear combination of $x^1,\ldots,x^m\in\Rmaxmn$;
\item[(2)] for any $j\in\text{supp}(y)$, there exists some $x^l$ from
$x^1,\ldots,x^m$ such that $x^l\leq_j y$.
\end{itemize}
\end{proposition}

\begin{proposition} 
\label{min-extr}
{\rm \cite{BSS-07}}
Let a semimodule $K$ be generated by a subset
$S$ of $\Rmaxmn$. The following are equivalent.
\begin{itemize}
\item[(1)] $y$ is an extremal of $K$;
\item[(2)] for some $j$, this $y$ is a minimal element of $S$ 
(and, equivalently, of $K$) with respect to $\leq_j$.
\end{itemize}
\end{proposition} 

Proposition \ref{min-extr} enables to treat idempotent extremals as
minima. The problem of finding partial maxima (and minima) in 
$n$-dimensional real space was
investigated by F. Preparata et al. \cite{PS:85}. The following
estimate is derived from their results.

\begin{theorem}
\label{compcomp}
{\rm \cite{Jos-05}(for $n=3$),\cite{BSS-07}}
Let $K$ be a semimodule in $\Rmaxmn$ generated by $k$ elements.
The problem of finding all extremals of $K$ requires 
not more than $O(k\log_2 k)$ operations, if $n=3$, and not more than
$O(k(\log_2 k)^{n-3}$ operations, if $n>3$ (with $n$ fixed).
\end{theorem}

Propositions \ref{spanit} and \ref{min-extr} imply a number of statements
for generators of idempotent 
semimodules in $\Rmaxmn$, see \cite{BSS-07} for details. Here I mention
two of them. 

\begin{theorem}
\label{mainres}
{\rm \cite{Wag-88,BSS-07}}
Let $K$ be a semimodule in $\Rmaxmn$ generated by $S$, and let
$E$ be the set of extremals of $K$ such that
$||u||=1$\footnote{the choice of norm does not matter} 
for all $u\in E$. Then $S=E\cup F$, where $F$ is redundant in the
sense that $S-\{u\}$ generates $K$ for any $u\in F$.
\end{theorem}

As a corollary, the weak basis of a semimodule 
is essentially unique whenever it exists.
The following is a tropical version of Minkowski's theorem.

\begin{theorem}
\label{mink}
{\rm \cite{GK-07<LAA>,BSS-07}}
A closed%\footnote{in Euclidean topology} 
semimodule in $\Rmaxmn$ is generated by its extremals.
\end{theorem}

\section{Cellular closures}

Algebraic closure of a square matrix
$A$ is the series $I\oplus A\oplus A^2\oplus\ldots$, where $I$ is the
identity matrix. This series converges iff $\lambda(A)\leq 1$,
where $\lambda(A)$ is the maximal cycle mean of $A$. This $\lambda(A)$
is also the maximal eigenvalue of the problem $Ax=\lambda x$. The
corresponding eigenspace will be denoted by $\text{eig}(A)$.

The following theorem and its corollary
are the ``ground stone'' of this section.
\begin{theorem}
\label{*=*1}
{\rm \cite{Ser-07<DC>}}
Let $A$ and $B$ be two square matrices such that $\lambda(A)\leq 1$
and $\lambda(B)\leq 1$. Then $A^*=B^*$ iff the spaces generated by 
columns of $A^*$ and $B^*$ coincide.
\end{theorem}
A square matrix $A$ is {\em definite}, if $\lambda(A)=1$ and all the diagonal
entries equal $1$.
\begin{corollary}
\label{*=*2}
Let $A$ and $B$ be two definite matrices. Then $A^*=B^*$ iff 
$\text{eig}(A)=\text{eig}(B)$.
\end{corollary}
I consider now the concepts of \cite{DS-04}.
Let $A$ be an $n\times m$ matrix over $\Rmaxm$ and $y$ an $n$-component
vector. Denote the collection $S=\{S_j\colon j\in\text{supp}(y)\}$,
where $S_j=\{i\colon y\geq_j A_{\tchk i}\}$, by
$\text{type}(y\mid A)$ and call it the {\em combinatorial
type} of $y$ with respect to $A$. Combinatorial types can be formally
defined as arbitrary collections of not more than $n$ possibly empty
subsets of $\{1,\ldots,m\}$. Denote the set of indices $i$, whose $S_i$ are
present in the type, by $\text{supp}(S)$. If $S=\text{type}(y\mid A)$ for some $y$, then
$\text{supp}(S)=\text{supp}(y)$. The types are partially ordered by the 
rule $S\subseteq S'$ if $\text{supp}(S')\subseteq\text{supp}(S)$ and
$S_i\subseteq S'_i$ for all $i\in\text{supp}(S)$. The set
$$
X^S=\{z\colon S\subseteq\text{type}(z\mid A)\}
$$
is the {\em region} of $S$. If $A_{ik}\ne 0$ for all $i\in S_k$, then
$S$ is {\em compatible} and we introduce the matrix $A^S$ by
$$
A_{\tchk i}^S=
\begin{cases}
\bigoplus_{k\in S_i} A_{\tchk k}/A_{ik}, & \text{if $i\in$supp($S)$ and 
$S_i\ne\emptyset$;}\\
e_i, & \text{if $i\in$supp($S)$ and 
$S_i=\emptyset$;}\\ 
\textbf{0}, & \text{if $i\notin$supp($S)$}.
\end{cases}
$$
If the region $X^S$
is not empty, then $X^S=\text{eig}(A^S)$.
Hence any region is (essentially) the eigenspace of a definite matrix, and
we use Corollary \ref{*=*2}.\markboth{Oleg Yu.~Shvedov}
{Minimal elements and cellular closures}
\begin{theorem}
\label{cellclos}
If $S$ and $T$ are (compatible) types such that $X^S$ and $X^T$ are not
empty and $X^S=X^T$, then $(A^S)^*=(A^T)^*$.
\end{theorem}
Theorem \ref{cellclos} enables to define various {\em cellular closures} of
$A$ to be $(A^S)^*$. This operation is correctly defined for every region,
being independent of the type. 

\msn Consider now the case when $A$ is a square $n\times n$ matrix with a
permutation $\sigma$ whose weight $\odot_{i=1}^n A_{i\sigma(i)}$ is 
nonzero. A permutation with maximal weight is called {\em maximal}.
We define $D^{\sigma}$ to be the matrix such that $D_{ij}^{\sigma}=A_{ij}$
if $j=\sigma(i)$ and $D_{ij}=0$ otherwise. If $\sigma$ is maximal,
then $(D^{\sigma})^{-1}A$ is definite and is called the 
{\em definite form} of $A$ \cite{Ser-07<DC>}. Different maximal permutations
lead to different definite forms. But we have that eigenspaces of all
definite forms coincide
(see \cite{Ser-07<DC>}), and by Corollary \ref{*=*2} closures
of all definite forms are equal. 

\msn Thus, for any square matrix $A$ with
nonzero permutations, we can define its {\em definite closure} to be
$(D^{\sigma})^{-1}A)^*$, where $\sigma$ is a maximal permutation.
Definite closure is a cellular closure, since 
$\text{eig}((D^{\sigma})^{-1}A)$ is the same as $X^S$ with 
$S=(\{\sigma(1)\},\ldots,\{\sigma(n)\})$, where $\sigma$ is any 
maximal permutation.

\renewcommand{\theequation}{\arabic{equation}}
\setcounter{section}{0}
\setcounter{footnote}{0}
\nachaloe{Oleg Yu.~Shvedov}{Semiclassical quantization of field theories}%
{Supported by the RFBR grant  05-01-00824 
and the joint RFBR/CNRS  grant 05-01-02807.}
\markboth{Oleg Yu.~Shvedov}{Semiclassical quantization of field theories}
\label{shv-abs}

{\bf \S1.}
It is well-known that equations of  quantum  field  theory  (QFT)  are
ill-defined \cite{5}. One  usually investigates the perturbative QFT instead of
"exact" QFT:  all quantities are presented as formal series in a small
perturbation parameter;  the  QFT  divergences are eleminated within a
perturbation framework only.

Semiclassical approximation \cite{13}
may be also viewed as an  expansion  in  a
small parameter.  Since  the  well-defined results are obtained within
the perturbation theory,  it seems to be more reasonable to talk about
semiclassical quantization rather than semiclassical approximation.

Consider the field theory model with the Lagrangian
$\mathcal L$ depending on the small parameter
$h$ ("Planck constant") as follows \cite{10}
(the scalar case is considered  for
the simplictiy):
\begin{equation}
{\mathcal L} = \frac{1}{2} \partial_{\mu} \varphi \partial^{\mu} \varphi -
\frac{1}{h} V(\sqrt{h} \varphi).
\label{2.1}
\end{equation}
with $V(\Phi)$ being a scalar potential.

In the formal quantum theory, field
$\hat{\varphi}({\bf  x})$ and momentum
$\hat{\pi}({\bf x})$
are viewed  as  operators   satisfying   the   canonical   commutation
relations. Semiclassical states depend on the small parameter $h$ as:
\begin{equation}
\Psi(t) \simeq e^{\frac{i}{h}S(t)}
e^{\frac{i}{\sqrt{h}} \int d{\bf x}
[\Pi({\bf x},t)  \hat{\varphi}({\bf  x})  - \Phi({\bf x},t) \hat{\pi}({\bf
x})]} f(t).
\label{2.2}
\end{equation}
Here $S(t)$ is a real c-number finction of $t$,
$\Phi({\bf
x},t)$ and  $\Pi({\bf  x},t)$  are  classical  fields  and  canonucally
conjugated momenta,
$\hat{\varphi}({\bf x})$  and  $\hat{\pi}({\bf  x})$ are quantum field
and momentum operators, $f(t)$ is a regular as $h\to 0$ state vector.

Superpositions of states \eqref{2.2}  are  also  viewed  as  semiclassical
states.

Presentation of  semiclassical  form  in  the  form  \eqref{2.2}  is  not
manifestly covariant. There are space and time coordinates. It happens
that the manifestly covariant  form  of  the  state  \eqref{2.2}  is  the
following:
\begin{equation}
\Psi \simeq
e^{\frac{i}{h}\overline{S}}
Texp\{\frac{i}{\sqrt{h}} \int dx J(x) \hat{\varphi}_h(x)\} \overline{f}
\equiv
e^{\frac{i}{h}\overline{S}}
T_J^h \overline{f}.
\label{2.3}
\end{equation}
Here
$\overline{S}$ is a real number,  $J(x)$ is a real function (classical
Schwinger   source),   $\hat{\varphi}_h(x)$   is  a  Heisenberg  field
operator,  $\overline{f}$ is a state vector being regular as $h\to 0$.
The  Schwinger  source  $J(x)$  should be rapidly damping at space and
time infinity \cite{9}.

{\bf \S2.}  Investigate  properties  of the semiclassical state \eqref{2.3}.
First of all,  note that the state $T^h_{J+\sqrt{h}\delta  J}  f$  can  be
expressed via  the  operator $T^h_J$.  To do this,  it is necessary to
investigate the operator
\begin{equation}
\underline{\Phi}_R(x|J)
\equiv - ih (T^h_J)^+ \frac{\delta T^h_J}{\delta J(x)}.
\label{3.1}
\end{equation}
It happens  to  coincide  with  the  well-known  LSZ  R-function
\cite{8}.

Notice that $\underline{\Phi}_R(x|J)$ is expanded in
$\sqrt{h}$; one writes
\begin{equation}
\underline{\Phi}_R(x|J) =
\Phi_R(x|J) + \sqrt{h} \Phi_R^{(1)}(x|J) + ...
\label{3.1aa}
\end{equation}
The c-number function
$\Phi_R(x|J)$ is called as a retarded classical field generated by the
Schwinger source $J$.
It is shown in
\cite{20} that for the model
\eqref{2.1}  $\Phi_R(x|J)$ is a solution of the equation
\begin{equation}
\partial_{\mu} \partial^{\mu} \Phi_R(x|J)
+ V'(\Phi_R(x|J)) = J(x),
\quad
\Phi_R\vert_{x{<\atop\sim} supp J} = 0.
\label{3.2a}
\end{equation}
which vanishes as $x^0 \to -\infty$.

The
following properties are corollaries of \eqref{3.1}.

1. The Hermitian property
\begin{equation}
\underline{\Phi}_R^+(x|J) = \underline{\Phi}_R(x|J).
\label{3.3o}
\end{equation}

2. The Poincare invariance property
\begin{equation}
\underline{U}_{g^{-1}}
\underline{\Phi}_R(x|u_gJ) \underline{U}_g =
\underline{\Phi}_R(w_gx|J).
\label{3.3a}
\end{equation}

3. The Bogoliubov causality property \cite{5}: R-function
$\underline{\Phi}_R(x|J)$
depends only on the source $J$ at the preceeding time moments.
Making use of the standard notations
$x>y$ iff $x^0-y^0 \ge |{\bf x} - {\bf y}|$,
$x<y$ iff $x^0-y^0 \le |{\bf x} - {\bf y}|$,
$x\sim y$ iff $|x^0-y^0| < |{\bf x} - {\bf y}|$, one rewrites
the Bogoliubov condition as
\begin{equation}
\frac{\delta \underline{\Phi}_R (x|J)}{\delta J(y)} = 0,
\quad y {>\atop\sim} x.
\label{3.3b}
\end{equation}

4. Commutation relation
\begin{equation}
[\underline{\Phi}_R(x|J);\underline{\Phi}_R(y|J)] = -ih \left(
\frac{\delta \underline{\Phi}_R(x|J)}{\delta J(y)} -
\frac{\delta \underline{\Phi}_R(y|J)}{\delta J(x)} \right).
\label{3.3c}
\end{equation}

5. Boundary condition at
$-\infty$.  If
$x{<\atop\sim} y$ for all $y\in supp J$,  the LSZ R-functiion does not depend
on the source:
\begin{equation}
\underline{\Phi}_R(x|J) = \hat{\varphi}_h(x) \sqrt{h},
\quad
x {<\atop\sim} supp J.
\label{3.3d}
\end{equation}
In particular, the classical retarded field vanishes as
$x^0 \to -\infty$.

Making use  of  the  operator   \eqref{3.1},   one   can   construst   the
semiclassical field:
\begin{displaymath}
\underline{\Phi}(x|J) = (T^h_J)^+ \hat{\varphi}_h(x) T^h_J
\end{displaymath}
coincides with
$\underline{\Phi}_R(x|J)$ at $x^0 \to +\infty$:
\begin{equation}
\underline{\Phi}(x|J) = \underline{\Phi}_R(x|J), \quad
x {>\atop\sim} supp J.
\label{3.4}
\end{equation}

{\bf \S 3.  } Another interesting feature of semiclassical states is that
some of them are approximately equal each other.
We say that
$J\sim 0$ iff
\begin{equation}
T^h_J \overline{f}      \simeq      e^{\frac{i}{h}     \overline{I}_J}
\underline{W}_J \overline{f}
\label{3.5}
\end{equation}
for some  number  $\overline{I}_J$  and   operator   $\underline{W}_J$
presented as a formal asymptotic series.

It is shown in \cite{20} for the model \eqref{2.1} that  the  source  $J$  is
equivalent to zero iff the retarded field generated by $J$ vanishes at
$+\infty$.

Analogously to   \cite{20},   one   derives the  following
properties.

1. Poincare invariance.
\begin{equation}
\underline{U}_g \underline{W}_J        \underline{U}_{g^{-1}}        =
\underline{W}_{u_gJ}, \qquad \overline{I}_{u_gJ} = \overline{I}_J;
\label{3.6}
\end{equation}

2. Unitarity
\begin{equation}
\underline{W}_J^+ =
\underline{W}_J^{-1};
\label{3.7}
\end{equation}

3. Bogoliubov causality:
as $J+\Delta J_2 \sim 0$,  $J+\Delta
J_1 + \Delta J_2 \sim 0$ and $supp \Delta J_2  {>\atop\sim}  supp  \Delta  J_1$,
the operator
{\sloppy

}
$$(\underline{W}_{J+\Delta J_2})^+ 
\underline{W}_{J+\Delta J_1
+ \Delta   J_2}$$   
and number $$-   \overline{I}_{J+\Delta   J_2}    +
\overline{I}_{J+\Delta J_1 + \Delta J_2}$$ 
do not depend on
$\Delta J_2$.

4. Variational property:
\begin{equation}
\delta \overline{I}_J - ih \underline{W}_J^+ \delta \underline{W}_J =
\int dx \underline{\Phi}_R(x|J) \delta J(x),
\label{3.8}
\end{equation}
which is valid as
$J\sim 0$ and $J+\delta J \sim 0$.

5. Boundary condition at $+\infty$:
\begin{equation}
\underline{\Phi}_R(x|J) =       \underline{W}_J^+       \hat{\varphi}(x)       \sqrt{h}
\underline{W}_{J}, \quad x {>\atop\sim} supp J.
\label{3.9}
\end{equation}

It follows from
\eqref{3.9} that the retarded classical field generated by the source
$J\sim  0$ will vanish at
$+\infty$.
For the model
\eqref{2.1},  an inverse statement is also valid:
{\it   for any field configuration $\Phi_c(x)$ with the compact support
one can uniquely choose a source
$J\sim  0$
(denoted as $J = J_{\Phi_c} = J(x|\Phi_c)$;
for example \eqref{2.1}, it is found from the relation
\eqref{3.2a})
generating
$\Phi_c(x)$ as a retarded classical field:
$\Phi_c(x) = \Phi_R(x|J)$;
it satisfies the locality condition
$\frac{\delta J(x|\Phi_c)}{\delta \Phi_c(y)} = 0$ as $x\ne y$.
}

It is possible to  treat  this  statement  as  a  basic  postulate  of
semiclassical field  theory.  Then the theory may be developed without
additional postulating  classical  stationary  action  principle   and
canonical commutation relation.

Namely, it  follows  from  eq.\eqref{3.8} in the leading order in $h$ that
the functional
\begin{equation}
I[\Phi_c] = \overline{I}_{J_{\Phi_c}} -
\int dx J_{\Phi_c}(x) \Phi_c(x)
\label{3.10}
\end{equation}
satisfies the "classical equation of motion"
\begin{equation}
J_{\Phi_c}(x) = - \frac{\delta I[\Phi_c]}{\delta \Phi_c(x)}.
\label{3.11}
\end{equation}
The functional $I[\Phi]$ should satisfy the locality condition
\begin{equation}
\frac{\delta^2I}{\delta \Phi_c(x) \delta \Phi_c(y)} = 0.
\quad x\ne y
\label{3.11a}
\end{equation}
This means that it is presented as an integral of a local Lagrangian.

Relation \eqref{3.11}  allows  us  to  reconstruct  the classical retarded
field,  making use of known  source  $J\sim  0$,  since  the  boundary
condition at  $-\infty$  is  known.  It  follows  from  
the Bogoliubov
causality condition that the retarded field
depends only on
$J$ at  the  preceeding  time  moments.  If  the  sourse  $J(x)$  is  not
equivalent to zero, it can be modified at $+\infty$ and transformed to
the sourse equivalent to zero.
Therefore, the relation
\begin{equation}
\frac{\delta I}{\delta \Phi_c} [\Phi_R(\cdot|J)] = - J(x),
\quad
\Phi_R\vert_{x<supp J} = 0
\label{3.11b}
\end{equation}
is valid for all sourses
$J$.
For the case $x{>\atop\sim} supp J$,
the property
\eqref{3.11b}
is taken to the classical field equation
\begin{equation}
\frac{\delta I}{\delta \Phi_c(x)} [\Phi(\cdot|J)] = 0.
\label{3.11c}
\end{equation}
This is  a  classical  stationary  action  principle.  It  is viewed a
coroollory of other general principles of semiclassical field theory.

Thus, we see that classical action  $I[\Phi_c]$  in  field  theory  is
related  with the phase of the state $T_J^h \overline{f}$ as $J\sim 0$
according to eq.\eqref{3.10}.

Let us rewrite the properties of the  operator  $\underline{W}_J$  via
the    field    $\Phi_c$.    Denote    $\underline{W}[\Phi_c]   \equiv
\underline{W}_{J_{\Phi_c}}$.

1. Poincare invariance.
\begin{equation}
\underline{U}_g \underline{W}[{\Phi}_c]  \underline{U}_{g^{-1}}
= \underline{W}[u_g{\Phi}_c].
\label{3.12a}
\end{equation}

2. Unitarity.
\begin{equation}
\underline{W}^+[{\Phi}_c] = (\underline{W}[{\Phi}_c] )^{-1}.
\label{3.12b}
\end{equation}

3. Bogoliubov causality.
\begin{equation}
\frac{\delta}{\delta {\Phi}_c(y)}
\left(
\underline{W}^+[{\Phi}_c]
\frac{\delta \underline{W}[{\Phi}_c]}
{\delta {\Phi}_c(x)}
\right) = 0, \quad y {>\atop\sim} x;
\label{3.12c}
\end{equation}

4. Yang-Feldman relation \cite{7}.
\begin{equation}
\int dy
\frac{\delta^2 I}{\delta {\Phi}_c(x) \delta {\Phi}_c(y)}
[\underline{\Phi}_R(y|J) - {\Phi}_c(y|J)] =
ih \underline{W}^+[{\Phi}_c]
\frac{\delta
\underline{W}[{\Phi}_c]}{\delta {\Phi}_c(x)}.
\label{3.12d}
\end{equation}

5. Boundary condition.
\begin{equation}
\underline{W}^+[{\Phi}_c] \hat{\varphi}_h(x) \sqrt{h}
\underline{W}[{\Phi}_c] =
\underline{\Phi}_R(x|J_{\Phi_c}), \quad
x {>\atop\sim} supp {\Phi}_c,
\label{3.12e}
\end{equation}
Here $\hat{\varphi}_h(x)  =  \underline{\Phi}_R(x|0)$
is the field operator without source.

{\bf \S 4.}
The covariant axioms of semiclassical field theory are as follows.

{\bf C1.} {\it A Hilbert state space $\mathcal F$ is given.}

{\bf C2.} {\it
An unitary represatation of the Poincare group is given. The operators
of the representation
$\underline{U}_g: {\mathcal  F}  \to  {\mathcal  F}$ are asymptoitc
series in $\sqrt{h}$.
}

{\bf C3.} {\it
To each classical source
$J(x)$   with compact support one assignes a retarded field
(LSZ R-function). It is an operator-valued distribution
$\underline{\Phi}_R(x|J)$ expanded in
$\sqrt{h}$ according to
\eqref{3.1aa}.  It satisfies the properties
\eqref{3.3o}, \eqref{3.3a}, \eqref{3.3b}, \eqref{3.3c}.
}

{\bf C4.} {\it
To each classical field configuration
$\Phi_c(x)$ with  compact  support  one  assigns  a c-number.  It is a
classical action
$I[\Phi_c]$ satisfying the locality condition
\eqref{3.11a}.  The property
$\Phi_c(x) = \Phi_R(x|J)$ is valid iff
\begin{equation}
J(x) = - \frac{\delta I[\Phi_c]}{\delta \Phi_c(x)}.
\label{3.11f}
\end{equation}
}
\markboth{Andre\u{\i} Sobolevski\u{\i}}{Semiclassical quantization of field theories}

{\bf C5.} {\it
To each classical field configuration $\Phi_c(x)$ with compact support
one  assigns  the   operator   $\underline{W}[\Phi_c]$   expanded   in
$\sqrt{h}$.   It   satisfies   the   relations  \eqref{3.12a},  \eqref{3.12b},
\eqref{3.12c}, \eqref{3.12d}, \eqref{3.12e}.
}

It is possible to develop a semiclassical perturbation theory,  making
use of these properties.

\renewcommand{\theequation}{\thesection.\arabic{equation}}
\setcounter{footnote}{0}
\nachaloe{Andre\u{\i} Sobolevski\u{\i}}{Convex analysis, transportation and 
reconstruction of peculiar velocities of galaxies}%
{Supported by the joint RFBR/CNRS grant 05-01-02807.}
\label{sob-abs}
\markboth{Alexander V. Stoyanovsky}{Convex analysis, transportation and galaxies}
We show how the problem of reconstruction of peculiar velocities of
galaxies starting from redshift-space catalogues can be rendered as a
convex quadratic optimization problem, invoking optimal transport
techniques for efficient large-scale astrophysical data processing. 
Connection with tropical algebra is briefly discussed.
\setcounter{section}{0}
\setcounter{footnote}{0}
\nachalo{Alexander V.~Stoyanovsky}{The Weyl algebra and quantization of fields}
\label{sto-abs}

In this talk we present a logically self-consistent procedure of quantization of fields.
In more detail our approach is exposed in the papers [1,2] and in the book [3].
As a basic example we use the $\varphi^4$ model in $4$-dimensional space-time.

\section{Difficulties of traditional approaches to quantum field theory}

Let us briefly discuss the logical contradictions in the known procedure of quantization
of fields. Usually one starts with the action functional
\begin{equation}
J=\int F(x^0,\ldots,x^n,u^1,\ldots,u^m,u^1_{x^0},\ldots,u^m_{x^n})\,
dx^0\ldots dx^n,
\end{equation}
where $x^0=t,x^1,\ldots,x^n$ are the independent variables, $u^1,\ldots,u^m$ are the dependent
variables, and $u^i_{x^j}=\frac{\partial u^i}{\partial x^j}$. For the $\varphi^4$ model
the action has the form
\begin{equation}
\label{phi4-J}
J=\int\left(\frac12\left(u_t^2-\sum_{j=1}^3 u_{x^j}^2-m^2u^2\right)
-\frac1{4!}gu^4\right)\,dtdx^1dx^2dx^3.
\end{equation}
One writes down the quantum field theory Schroedinger equation
\begin{equation}
\label{QFT-Sch}
ih\frac{\partial\Psi}{\partial t}=\int\widehat H\left(t,\x,u^i(\x),
\frac{\partial u^i}{\partial\x},-ih\frac{\delta}{\delta u^i(\x)}\right)
\Psi\,d\x,
\end{equation}
where $\x=(x_1,\ldots,x_n)$; the density of the 
Hamiltonian $H$ is the Legendre transform
of the Lagrangian $F$ with respect to 
the variables $u^i_t$; $\frac{\delta}{\delta u^i(\x)}$
is the variational derivative operator.
\markboth{Alexander V.~Stoyanovsky}{The Weyl algebra and quantization of fields}

Note that the Schroedinger equation
is not well defined in quantum field theory even
for the free scalar field ($g=0$ in (2)). For
example, if we consider mathematical 
equation \eqref{QFT-Sch} literally, then it is not difficult to
check that this equation has no nonzero 
four times differentiable solutions (the expression
for the derivative $\frac{\partial^2\Psi}{\partial t^2}$ 
has no sense). The traditional
approach is to ``subtract infinity'' from the RHS
of equation \eqref{QFT-Sch} and to solve it in the
Fock space of functionals. However, this approach contradicts
physical as well as mathematical
considerations. Physically, if states were functionals and energy were finite, then,
in principle, we could measure some quantities related with these functionals (such as energy).
However, it is known that quantum mechanical quantities like energy and momentum are
theoretically non-measurable in relativistic quantum dynamics, and the only measurable
quantities are the scattering sections. Mathematically, equation \eqref{QFT-Sch}
in the Fock space
does not admit a relativistically invariant generalization (usually called the
Tomonaga--Schwinger equation [4]), as shown in the important paper [5]. In this paper
it is shown that the evolution operators of the Klein--Gordon equation from one space-like
surface to another, which are symplectic transformations of the phase space
of the field, do not belong to the version of the infinite dimensional symplectic group
which acts on the Fock space.

So quantization of free fields, for example, following the
lines of the book [4], meets difficulties of the logical kind and is therefore not
completely satisfactory. Due to this fact, the renormalization procedure
for interacting fields, defined using the Bogolyubov--Parasyuk theorem,
gives us a model in which it is difficult to say how quantum field theory turns into
the classical one as $h\to0$.

\section{The infinite dimensional Weyl algebra}

The proposed way to overcome these difficulties is to replace the algebra of variational
differential operators by the infinite dimensional Weyl algebra defined below, which
admits an explicit action of the infinite dimensional group of continuous symplectic
transformations of the phase space of a field. This phase space is the Schwartz space
of functions $(u^i(s)$,$p^i(s))$, where $p^i(s)$ are the variables conjugate to $u^i(s)$,
and $s=(s_1,\ldots,s_n)$ are parameters on a spacelike surface,
with the Poisson bracket
\begin{equation}
\{\Phi_1,\Phi_2\}=\sum_i\int\left(\frac{\delta \Phi_1}{\delta u^i(s)}
\frac{\delta \Phi_2}{\delta p^i(s)}-\frac{\delta \Phi_1}{\delta p^i(s)}
\frac{\delta \Phi_2}{\delta u^i(s)}\right)ds
\end{equation}
of two functionals $\Phi_l(u^i(\cdot),p^i(\cdot))$, $l=1,2$.
Let us write this bracket in the form
\begin{equation}
\{\Phi_1,\Phi_2\}=\int\sum_{i,j}\omega^{ij}\frac{\delta \Phi_1}{\delta y^i(s)}
\frac{\delta \Phi_2}{\delta y^j(s)}\,ds,
\end{equation}
where $y^i=u^i$ for $1\le i\le m$ and $y^i=p^{i-m}$ for $m+1\le i\le 2m$, and
$\omega^{ij}=\delta_{i,j-m}-\delta_{i-m,j}$.
The Weyl algebra is defined as the algebra of weakly infinite differentiable functionals
$\Phi(u^i(\cdot),p^i(\cdot))$
with respect to the Moyal $*$-product
\begin{equation}
\begin{aligned}{}
&(\Phi_1*\Phi_2)(y^i(\cdot))\\
&=\left.\exp\left(-\frac{ih}2\int\sum_{i,j}\omega^{ij}
\frac{\delta}{\delta y^i(s)}\frac{\delta}{\delta z^j(s)}\,ds\right)
\Phi_1(y^i(\cdot))\Phi_2(z^i(\cdot))\right|_{z^i(\cdot)=y^i(\cdot)}.
\end{aligned}
\end{equation}
This product is not everywhere defined: for example, $u^i(s)*p^i(s)$ is undefined.
Note only that if all necessary series and integrals are absolutely convergent, then
the $*$-product is associative. Below we will be interested only in some concrete
computations in the Weyl algebra.

Let us replace the Schrodinger equation \eqref{QFT-Sch} 
by the Heisenberg equation in the Weyl algebra
\begin{equation}
\label{Hei-Weyl}
ih\frac{\partial\Phi}{\partial t}=\left[\int\!H(t,\x,u^i(\x),
\frac{\partial u^i}{\partial\x},p^i(\x))d\x,\Phi\right]
\end{equation}
and by its relativistically invariant generalization, where
\begin{equation}
[\Phi_1,\Phi_2]=\Phi_1*\Phi_2-\Phi_2*\Phi_1
\end{equation}
is the commutator in the Weyl algebra. The classical limits of equation \eqref{Hei-Weyl}
are the Hamilton equations
\begin{equation}
\frac{\partial\Phi}{\partial t}=\{\Phi,\int\!H\,d\x\}
\end{equation}
equivalent to the Euler--Lagrange equations.

\section{Quantization of free scalar field}

Put $g=0$ in \eqref{phi4-J}. Since the obtained Hamiltonian
\begin{equation}
H_0=\int\frac12(p(\x)^2+(\gradd u(\x))^2+m^2u(\x)^2)d\x
\end{equation}
is quadratic, we have
\begin{equation}
\frac1{ih}[H_0,\Phi]=\{\Phi,H_0\},
\end{equation}
hence,
$\Phi(t_1; u(\cdot),p(\cdot))$, subject to the Heisenberg equation, is obtained from
$\Phi(t_0;u(\cdot),p(\cdot))$ by the linear symplectic change of variables
\begin{equation}
(u(t_0,\x),p(t_0,\x)=u_t(t_0,\x))\to(u(t_1,\x),p(t_1,\x)=u_t(t_1,\x)),
\end{equation}
given by the evolution operator of the Hamilton equations, i.~e.,
of the Klein--Gordon equation, from the Cauchy surface $t=t_0$ to the Cauchy surface $t=t_1$.
Hence we can identify the Weyl algebras of various Cauchy surfaces by means of evolution
operators of the Klein--Gordon equation. In other words, we can consider the Weyl algebra
$W_0$ of the symplectic vector space of solutions $u(t,\x)$ of the Klein--Gordon equation
on the whole space-time. The symplectic form on this vector space is obtained by taking
Cauchy data on any spacelike surface (for example, on the surface $t=\const$).
Below we will fix this identification of the Weyl algebras of various spacelike surfaces.

Define the {\it vacuum average} linear functional
\begin{equation}
\Phi\to\langle0|\Phi|0\rangle
\end{equation}
on the Weyl algebra $W_0$ as the unique (not everywhere defined) functional with the
properties
\begin{equation}
\langle0|\Phi*u_-(t,\x)|0\rangle=\langle0|u_+(t,\x)*\Phi|0\rangle=0,\ \ \langle0|1|0\rangle=1.
\end{equation}
Here $u=u_++u_-$ is the decomposition of a solution $u(t,\x)$ of the Klein--Gordon
equation into the positive and negative frequency parts (we assume $m>0$ so that this
decomposition is unique). For $\Phi\in W_0$, define an operator in the standard Fock space
with the matrix elements
\begin{equation}
\langle0|\tilde u_-(-p'_{(1)})\ldots\tilde u_-(-p'_{(N')})*\Phi*\tilde u_+(p_{(1)})\ldots
\tilde u_+(p_{(N)})|0\rangle.
\end{equation}
Here $\tilde u_\pm(p)$ is the Fourier transform (the momentum representation) of $u_\pm$,
$p=(p_0,\ldots,p_n)$.

One can check the following two properties of this correspondence:

i) $*$-product of functionals $\Phi$ goes to composition of operators in the Fock space,
so that this correspondence is a (not everywhere defined) homomorphism from the algebra
$W_0$ to the algebra of operators in the Fock space;

ii) complex conjugation of functionals $\Phi$ goes to Hermitian conjugation of operators in
Hilbert space.

\section{Quantization of interacting fields}

{\bf Statement.} {\it There exists a map from the set of smooth functions $g=g(t,\x)$
with compact support to the set of functionals $P(g)\in W_0$ with the following properties.

1) $P(g)$ is a formal series in $g$ with the first two terms
\begin{equation}
P(g)=1+\frac1{ih}\int g(t,\x)u(t,\x)^4/4!\,dtd\x+\ldots.
\end{equation}

2) Classical limit:
$P(g)=a(g,h)\exp(iR(g)/h)$ where $a(g,h)$ is a formal series in $h$, and conjugation
by $\exp(iR(g)/h)$ in the Weyl algebra $W_0$ up to $O(h)$ yields the
perturbation series for the evolution operator of the nonlinear classical field equation
\begin{equation}
\Box u(x)-m^2u(x)=g(x)u^3(x)/3!
\end{equation}
from $t=-\infty$ to $t=\infty$.

3) The Lorentz invariance condition:
\begin{equation}
LP(L^{-1}g)=P(g)
\end{equation}
for a Lorentz transformation $L$.

4) The unitarity condition:
\begin{equation}
P(g)*\overline{P(g)}=1,
\end{equation}
where $\overline{P(g)}$ is complex conjugate to $P(g)$.

5) The causality condition: for two functions $g_1,g_2$ equal for $t\le t_0$, the product
$P(g_1)*P(g_2)^{-1}$ does not depend on the behavior of the functions $g_1,g_2$ for $t<t_0$.

6) The quasiclassical dynamical evolution (cf. with the Maslov--Shvedov quantum field
theory complex germ [6]): for any spacelike surfaces $\cC_1$, $\cC_2$
there exists a limit $P_{\cC_1,\cC_2}$ of
$P(g)$ modulo $o(h)$ as the function $g(x)$ tends to $1$ if $x$ belongs to the strip between
the spacelike surfaces and to $0$ otherwise. This limit possesses the property
\begin{equation}
P_{\cC_1,\cC_3}=P_{\cC_2,\cC_3}*P_{\cC_1,\cC_2}+o(h).
\end{equation}

7) The $S$-matrix: there exists a limit $P$ of $P(g)$ as $g(x)$ tends to the function
$g=\const$. This $P$ is a formal power series in $g$.

Any other choice of $P(g)$ with the properties 1--7 above is equivalent to some change of
parameters $m,g(x)$.
}

\markboth{Svetlana V.~Tsykina}{The Weyl algebra and quantization of fields}
This statement is completely similar to the Bogolyubov--Parasyuk theorem. Moreover,
if we denote by $S(g)$ the operator in the Fock space corresponding to $P(g)$ and by $S$
the operator corresponding to $P$, then $S(g)$ is exactly the Bogolyubov $S$-matrix
and $S$ is the physical $S$-matrix. The elements $P(g)$ are constructed in the same way
as $S(g)$ in [4], using the renormalization procedure, the main difference
being that composition of operators is replaced
by $*$-product of functionals, and the normally ordered product of operators is replaced
by the usual (commutative) product of functionals.

Note that the conditions on $P(g)$ (in particular, the causality condition)
are the natural analogs of conditions on dependence of
the evolution operator of a partial differential equation on the coefficient functions of this
equation. Therefore the above apparatus is similar to the scattering theory in the theory
of partial differential equations.

Note also that the presence of interaction cutoff function $g(x)$ is necessary
from the physical point of view, since the scattering particles are considered as
non-interacting at infinity (which means that $g=0$ at infinity).

\setcounter{section}{0}
\setcounter{footnote}{0}
\nachaloe{Svetlana V.~Tsykina}{\mbox{Polynomial quantization on 
para-hermitian spaces} \mbox{with pseudo-orthogonal group of translations}}%
{Supported by the Russian Foundation for Basic
Research (grant No. 05-01-00074a), the Netherlands Organization for
Scientific Research (NWO) (grant No. 047-017-015), the Scientific
Programs "Devel. Sci. Potent. High. School" (project RNP 2.1.1.351
and Templan, No. 1.2.02).}
\markboth{Svetlana V.~Tsykina}{Polynomial quantization with 
pseudo-orthogonal group}
\label{tsy-abs}
We construct polynomial quantization, which is a variant of quantization in 
spirit of Berezin, on para-Hermitian symmetric spaces $G/H$ with the
pseudo-orthogonal group $G={\rm SO}_0(p,q)$. For all these spaces,
the connected component $H_e$ of the subgroup $H$ containing the
identity of $G$ is the direct product ${\rm SO}_0(p-1,q-1)\times
{\rm SO}_0(1,1)$, so that $G/H$ is covered by $G/H_e$ (with
multiplicity 1, 2 or 4). The 
dimension of $G/H$ is equal to ${2n-4}$,
where ${n=p+q}$. We restrict ourselves to the spaces $G/H$ that are
$G$-orbits in the adjoint representation of $G$.

%Old
%A construction of quantization on arbitrary para-Hermitian symmetric
%spaces was offered in [2]. In the case of polynomial quantization
%covariant and contravariant symbols are polynomials on $G/H$. We
%introduce a multiplication of covariant symbols, establish the
%correspondence principle, study the Berezin transform. Polynomial
%quantization on rank one para-Hermitian symmetric spaces was
%constructed in [3]. Our spaces $G/H$ with $G={\rm SO}_0(p,q)$ have
%rank 2.

%New
A construction of quantization on arbitrary para-Hermitian symmetric
spaces was given in [2]. 
The term "polynomial quantization" means in particular that both
covariant and contravariant symbols are polynomials on $G/H$. 
Following the general
scheme of [2], we introduce multiplication of 
covariant symbols, establish the
correspondence principle, and study the Berezin transform. 

The polynomial
quantization on rank one para-Hermitian symmetric spaces has been
constructed in [3]. In this paper, we consider the spaces $G/H$ with $G={\rm SO}_0(p,q)$. Note that these spaces 
have rank 2.

\section{The pseudo-orthogonal group and its Lie algebra}

%Old
%Let us equip the space ${\Bbb R}^n$ with the following bilinear
%form:
%New
Consider the space ${\Bbb R}^n$ equipped with the following bilinear
form:
\begin{equation}
[x,y]=\sum\limits_{i=1}^n{\lambda}_i x_i y_i,\nonumber
\end{equation}
where ${\lambda_1=\ldots =\lambda_p=-1, \ \lambda_{p+1}=\ldots
=\lambda_n=1}$, and $x=(x_1,\ldots ,x_n)$, $y=(y_1,\ldots ,y_n)$ are
vectors in ${\Bbb R}^n$.
{\sloppy

}
%Old
%Let $G$ denote the group ${\rm SO}_0\,(p,q)$, the connected
%component of the identity in the group of linear transformations
%with determinant 1 of the space ${\Bbb R}^n$, preserving the
%bilinear form $[x,y]$. We consider that $G$ acts linearly on ${\Bbb
%R}^n$ from the right: $x\mapsto xg$. In accordance with that, we
%write vectors in the row form. We also assume that ${p>1,q>1}$, the
%general case.

%New
Let $G$ denote the group ${\rm SO}_0\,(p,q)$. This group is the connected
component of the identity, in the group of linear transformations
of ${\Bbb R}^n$ that preserve $[x,y]$ and have determinant equal to $1$. 
We assume that $G$ acts linearly on ${\Bbb
R}^n$ from the right: $x\mapsto xg$. In accordance with that, we
write vectors in the row form. We also assume that ${p>1,q>1}$.

Let us write matrices $g\in G$ in the block form corresponding to
the partition $n=1+(n-2)+1$. Denote by $H$ the subgroup of $G$
consisting of matrices
\begin{equation}
\label{h}
h=\left(
\begin{array}{ccc}
\alpha &0&\beta \\
0&v&0 \\
\beta&0&\alpha
\end{array}
\right),
\end{equation}
where ${\alpha^2-\beta^2=1}$, ${v\in{\rm SO}(p-1,q-1)}$. The
subgroup $H$ consists of two connected components. The connected
component $H_e$, containing the unit matrix $E$ of $G$ consists of
matrices~(\ref{h}), where $\alpha={\rm ch} t$, $\beta={\rm sh} t$.
Thus, it is ${\rm SO}_0(p-1,q-1)\times{\rm SO}_0(1,1)$. The second
connected component of $H$ (which does not contain $E$) 
%Old
%has as its
%representative the matrix ${\rm diag}\,\{-1,-1,1,\ldots,1,-1,-1\}$.
%New
contains the matrix ${\rm diag}\,\{-1,-1,1,\ldots,1,-1,-1\}$, as a representative.

The Lie algebra $\got g$ of $G$ consists of real matrices $X$ of
order $n$ satisfying the condition ${X'I+IX=0}$, where ${I={\rm
diag}\,\{\lambda_1,\ldots,\lambda_n\}}$, the prime denotes matrix
transposition.

Let
\begin{equation}
\label{z0} Z_0=\left(
\begin{array}{ccc}
0&0&1 \\
0&0&0 \\
1&0&0
\end{array}\right).
\end{equation}
The stabilizer of $Z_0$ in the adjoint representation is exactly the
group $H$, therefore, the manifold $G/H$ is just the $G$-orbit of
the matrix $Z_0$ in ${\got g}$.

The operator ${\rm ad}\,Z_0$ has three eigenvalues: $-1,0,+1$. Respectively,
the Lie algebra $\got g$ is decomposed into the direct sum of
eigenspaces
$$
\got g={\got q}^-+\got h+{\got q}^+,
$$
where $\got h$ is the Lie algebra of $H$. The subspaces ${\got
q^-,\got q^+}$ consist of matrices
$$
X_{\xi}:\ \left(
\begin{array}{ccc}
0&\xi&0 \\
{\xi}^*&0&{\xi}^* \\
0&-\xi&0
\end{array}\right),\qquad
Y_{\eta}:\ \left(
\begin{array}{ccc}
0&\eta&0 \\
{\eta}^*&0&-{\eta}^* \\
0&\eta&0
\end{array}\right)
$$
respectively, where ${\xi,\eta}$ are rows in ${\Bbb R}^{n-2}$. Both
spaces ${\got q}^{\pm}$ are Abelian subalgebras of ${\got g}$, they
have dimension $n-2$. The subgroup $H$ preserves both subspaces
$\got q^-$ and $\got q^+$ in the adjoint action:
\begin{equation}
\label{actz}
Z\mapsto h^{-1}Z h, \ h\in H.
\end{equation}
Let  $h\in H$ have the form (\ref{h}). For simplicity, we identify
matrices $X_{\xi}$ and $Y_{\eta}$ with vectors $\xi$ and $\eta$,
respectively.  Under the action (\ref{actz}) vectors $\xi\in\got
q^-$ and $\eta\in\got q^+$ are transformed as follows:
 \begin{equation}
\xi\mapsto{\widetilde\xi}= (\alpha+\beta)\xi v, \ \ \eta\mapsto
{\widehat\eta}=(\alpha-\beta)\eta v.
\end{equation}

%Old
%Equip the space ${\Bbb R}^{n-2}$ with a bilinear form 
%by means of
%the matrix 
%New
Consider the space ${\Bbb R}^{n-2}$ with bilinear form 
defined by the matrix 
$I_1={\rm diag}\,\{\lambda_2,\ldots,\lambda_{n-1}\}$:
{\sloppy

}
\begin{equation}
\langle\xi,\eta\rangle=\sum_{i=2}^{n-1}\lambda_i\xi_i\eta_i.\nonumber
\end{equation}

\section{Representations of $G$ associated with a cone}

The group ${G={\rm SO}_0(p,q)}$ preserves manifolds ${[x,x]=c}$,
${c\in\Bbb R}$, in ${\Bbb R}^n$. Let $\cC$ be the cone
${[x,x]=0}$, ${x\ne 0}$, in ${\Bbb R}^n$. Let us fix two points in
the cone: ${s^+=(1,0,\ldots,0,1)}$, ${s^-=(1,0,\ldots,0,-1)}$.
Consider the following two sections of the cone:
{\sloppy

}
\begin{eqnarray}
\Gamma^+&=&\{x_1+x_n=2\}=\{[x,s^-]=-2\},\nonumber \\
\Gamma^-&=&\{x_1-x_n=2\}=\{[x,s^+]=-2\}.\nonumber
\end{eqnarray}

The points $s^+,\ s^-$ belong to $\Gamma^+,\ \Gamma^-$ respectively.
They are eigenvectors of the maximal parabolic subgroups
${P^+=Q^+H}$ and ${P^-=Q^-H}$ respectively, with eigenvalues
${\alpha-\beta}$ and ${\alpha+\beta}$, where $\alpha,\beta$ are
parameters of  $h\in H$, see (\ref{h}). Here $Q^-=\exp{\got q}^-$,
$Q^+=\exp{\got q}^+$.

The section $\Gamma^\pm$ meets almost all
%??????
%rulings 
%??????
generatrices
of the cone $\cC$. The linear action of $G$ on the cone 
%Old
%gives rise to
%New
induces
the following
actions of $G$ on $\Gamma^-$ and $\Gamma^+$ respectively:
\begin{eqnarray}
x\longmapsto{\widetilde x}=-\frac2{[xg,s^+]}\cdot
xg,\ x\in \Gamma^-,\label{gam-}  \\
x\longmapsto{\widehat x}=-\frac2{[xg,s^-]}\cdot xg,\ x\in
\Gamma^+,\label{gam+}
\end{eqnarray}
defined almost everywhere on $\Gamma^\pm$. For the subgroups $Q^-$
and $Q^+$ respectively, these actions turn out to be linear:
$x\mapsto xg$. Moreover, the subgroups $Q^\pm$ act on $\Gamma^\pm$
%Old
%simply-transitive.
%New
simply transitively.
This allows to define the coordinates
${\xi=(\xi_2,\ldots,\xi_{n-1})}$ on $\Gamma^-$ and
${\eta=(\eta_2,\ldots,\eta_{n-1})}$ on $\Gamma^+$ transferring them
from  ${\got q}^-$ on ${\got q}^+$ respectively, namely, for
$u\in\Gamma^-$ and $v\in\Gamma^+$ we set:
\begin{eqnarray}
u&=&u(\xi)=s^-e^{X_\xi}=(1+\langle\xi,\xi\rangle,2\xi,-1+\langle\xi,\xi\rangle),\label{u} \\
\label{v}
v&=&v(\eta)=s^+e^{Y_\eta}=(1+\langle\eta,\eta\rangle,2\eta,1-\langle\eta,\eta\rangle).
\end{eqnarray}
The stabilisers in $G$ of the points ${s^-\in\Gamma^-}$ and
${s^+\in\Gamma^+}$ under the actions (\ref{gam-}) and (\ref{gam+})
are the subgroups ${P^+=Q^+H}$ and ${P^-=Q^-H}$ respectively.

Let ${\sigma\in{\Bbb C},\ \varepsilon=0,1}$. Let ${\cD}_{\sigma,\,\varepsilon}({\cC})$ be the space of $C^\infty$
functions $f$  on the cone ${\cC}$ with homogeneity ${\sigma}$
and parity ${\varepsilon}$, i.e.
$$
f(tx)=t^{\sigma,\,\varepsilon}f(x),\ x\in{\cC},\ t\in{\Bbb
R}^*={\Bbb R}\setminus\{0\},
$$
where we denote ${t^{\sigma,\varepsilon}=|t|^{\sigma}{\rm sgn}
^{\varepsilon}t}$. 
%Old
%A representation $T_{\sigma,\varepsilon}$ of the
%group $G$ acts in ${\cal D}_{\sigma,\,\varepsilon}({\cC})$ by
%translations:
%${\left(T_{\sigma,\,\varepsilon}(g)f\right)(x)=f(xg)}$.
%New
Denote by $T_{\sigma,\varepsilon}$ the representation of $G$ which acts on
${\cD}_{\sigma,\,\varepsilon}({\cC})$ by
translations:
${\left(T_{\sigma,\,\varepsilon}(g)f\right)(x)=f(xg)}$.

%Old
%Realize now the representation
%$T_{\sigma,\,\varepsilon}$ on the
%functions on the sections $\Gamma^\pm$. 
%New
Consider now the restrictions of functions from
${\cD}_{\sigma,\,\varepsilon}({\cC})$ to  the sections $\Gamma^\pm$.
Such restrictions
form a space ${\cD}_{\sigma,\,\varepsilon}(\Gamma^\pm)$ of
functions $f$ on $\Gamma^\pm$. This space is contained in
$C^\infty(\Gamma^\pm)$ and contains ${\cD}(\Gamma^\pm)$. In the
coordinates ${\xi,\eta}$, the representation
$T_{\sigma,\,\varepsilon}$ of the group $G$ acts on the space of restrictions 
${\cD}_{\sigma,\,\varepsilon}({\cC})$ by
\begin{eqnarray}
\left(T_{\sigma,\,\varepsilon}(g)f\right)(\xi)&=&f(\widetilde{\xi})\left\{-\frac12[ug,s^+]\right\}^{\sigma,\,\varepsilon},\\
\left(T_{\sigma,\,\varepsilon}(g)f\right)(\eta)&=&f(\widehat{\eta})\left\{-\frac12[vg,s^-]\right\}^{\sigma,\,\varepsilon},
\end{eqnarray}
where ${u=u(\xi)},\ {v=v(\eta)}$ are defined by~(\ref{u}),~(\ref{v}), actions ${\xi\mapsto\widetilde{\xi}}$
and ${\eta\mapsto\widehat{\eta}}$ are defined by~(\ref{gam-}),~(\ref{gam+}).

Define the operator $A_{\sigma,\,\varepsilon}$ on ${\cD}_{\sigma,\,\varepsilon}(\Gamma^\pm)$  by:
\begin{equation}
(A_{\sigma,\,\varepsilon}f)(\xi)=\int\limits_{{\Bbb
R}^{n-2}}N(\xi,\eta)^{2-n-\sigma,\,\varepsilon}f(\eta)d\eta,
\label{spletop}
\end{equation}
where
$$
N(\xi,\eta)=-\frac12[u(\xi),v(\eta)]=1-2\langle\xi,\eta\rangle+\langle\xi,\xi\rangle\langle\eta,\eta\rangle.
$$
The function ${N(\xi,\eta)}$ is a polynomial in $\xi,\eta$. The
operator $A_{\sigma,\,\varepsilon}$ intertwines the representations
$T_{\sigma,\,\varepsilon}$ and $T_{2-n-\sigma,\,\varepsilon}$. These
representations act on functions on {\it different} sections. We can
change the position of $\xi$ and $\eta$ in~(\ref{spletop}). The
product $A_{2-n-\sigma,\,\varepsilon}A_{\sigma,\,\varepsilon}$ is a
scalar operator:
$$
A_{2-n-\sigma,\,\varepsilon}A_{\sigma,\,\varepsilon}=\omega_0(\sigma,\varepsilon)E,
$$
where
\begin{eqnarray}
\omega_0(\sigma,\varepsilon)&=&2^3\pi^{n-3}\frac{\Gamma(\sigma\!+\!1)\Gamma(3\!-\!n\!-\!\sigma)}{(2\sigma\!+\!n\!-\!2)\sin\left(\sigma\!+\!\frac{n}2\right)\pi}\times\nonumber\\
&\times&\sin\!{\frac{\sigma\!-\!\varepsilon}2\pi} \cdot
\sin\!{\frac{\sigma\!-\!\varepsilon\!+\!p}2\pi} \cdot
\sin\!{\frac{\sigma\!+\!\varepsilon\!+\!q}2\pi} \cdot
\sin\!{\frac{\sigma\!+\!\varepsilon\!+\!n}2\pi}.\nonumber
\end{eqnarray}

\section{The space $G/H$}

Consider the following realization of the space $G/H$. Let ${\Omega}$ be the set
of matrices:
\begin{equation}
\label{matr-z}
z=\frac{y^*x}{[x,y]},
\end{equation}
where ${x,y\in{\cC}}$, ${y^*=Iy'}$. For these matrices, rank and
trace are equal to $1$. The adjoint action $z\mapsto g^{-1}zg$
preserves ${\Omega}$. The stabilizer of the matrix $z^0$,
corresponding to the pair ${x=s^-,y=s^+}$, is the subgroup $H$, so
that $\Omega$ is just $G/H$.

%Old
%As vectors $x,y\in{\cC}$ in (\ref{matr-z}), let us take vectors
%$u=u(\xi)$ and $v=v(\eta)$ in the sections $\Gamma^-$ and $\Gamma^+$
%of the cone $\cC$ respectively. 
%New
Take vectors $u=u(\xi)$ and $v=v(\eta)$ in the sections $\Gamma^-$ and $\Gamma^+$
of the cone $\cC$, respectively, for $x$ and $y$ in (\ref{matr-z})
We obtain an embedding
$\Gamma^-\times\Gamma^+\to\Omega$ given by
{\sloppy

}
\begin{equation}
\label{vlogen}
z=z(\xi,\eta)=\frac{v^* u}{[u,v]},\quad u=u(\xi), \
v=v(\eta),
\end{equation}
The map ${(u,v)\longmapsto z}$ given by formula (\ref{vlogen}) is
defined for ${\xi,\eta\in{\Bbb R}^{n-2}}$ such that
${N(\xi,\eta)\ne0}$, since $[u,v]=-2N(\xi,\eta)$. Therefore,
vectors $\xi,\eta\in{\Bbb R}^{n-2}$ with the condition
$N(\xi,\eta)\ne0$ are local coordinates on ${\Omega}$. The adjoint
action of the group $G$ on ${\Omega}$ is generated by its actions on
$\xi$ and $\eta$. For each ${g\in G}$, this action is defined on a
dense set of $\Omega$.

We can identify the tangent space of $G/H$ at the initial point
$z^0$ with the space ${{\goth q}={\goth q}^-+{\goth q}^+}$ in the Lie
algebra ${\goth g}$. Let $S(\got q)$ denote the algebra of
polynomials on $\got q$. The action (\ref{actz}) of the group $H$ on
$\got q$ 
%Old
%gives rise to
%New
induces
an action of $H$ on $S(\got q)$. Let $S(\got
q)^H$ denote the algebra of polynomials invariant with respect to
$H$. This algebra is generated by two polynomials
${\langle\xi,\eta\rangle}$ and
${\langle\xi,\xi\rangle\langle\eta,\eta\rangle}$.

Let ${\Bbb D}(G/H)$ denote the algebra of differential operators on
$G/H$ invariant with respect to $G$. This algebra is in the
one-to-one correspondence with the algebra $S(\got q)^H$. Let
$\Delta_2$ and $\Delta_4$ denote operators in ${\Bbb D}(G/H)$
corresponding to generators $\langle\xi,\eta\rangle$ and
$\langle\xi,\xi\rangle\langle\eta,\eta\rangle$ of $S(\got q)^H$
respectively. Let us call these operators $\Delta_2$ and $\Delta_4$
the Laplace operators on  $G/H$. The operator $\Delta_2$ is the
Laplace-Beltrami operator. These operators are differential
operators of the second and the fourth order respectively, they are
generators in ${\Bbb D}(G/H)$. Explicit expressions of them are very
cumbersome. We write explicit expressions for their radial parts
$\buildrel 0\over{\Delta}_2$ and $\buildrel 0\over{\Delta}_4$ in
horospherical coordinates.

These coordinates are defined as follows. Let us take in ${\got
q=\got q^++\got q^-}$ the Cartan subspace ${\got a}$, consisting of
matrices
$$
A_t=\left(
\begin{array}{ccccc}
0&0&0&t_1&0\\
0&0&0&0&t_2\\
0&0&0&0&0\\
t_1&0&0&0&0\\
0&t_2&0&0&0
\end{array}
\right),
$$
where $t=(t_1,t_2)\in{\Bbb R}^2$. 
Introduce in ${\got a^*}$ the
%Old
%alphabethical
%New
lexicographical
order in coordinates.
Let $\got n$ denote the
subalgebra of $\got g$ formed by the corresponding positive root spaces.
Let ${A=\exp {\got a},\ N=\exp{\got n}}$. Consider the set of points
$z$ in $\Omega$ obtained from $z^0$ via the translation by
$a=a(t_1,t_2)\in A$ and then by $n\in N$, i.e.
${z=n^{-1}a^{-1}z^0an}$. It is a neighbourhood $U$ of the point
$z^0$. 
%??????????????
%For coordinates in this neighbourhood (horospherical
%coordinates) we obtain $t_1,t_2$ and also parameters of the subgroup
%$N$.
%??????????????
Parameters $t_1,t_2$ of the subgroup $A$ and also parameters of the
subgroup $N$ are coordinates in this neighbourhood (horospherical
coordinates).

Let $f$ be a function defined on $U$ that does not depend
 on ${n\in N}$. Then it is a function of ${t=(t_1,t_2)}$: $f(z)=F(t)$. Let $D$ be
a differential operator in ${\Bbb D}(G/H)$. Then $Df$ also does not
depend on ${n\in N}$:
$$
 Df =\buildrel 0\over{D}F,
$$
where $\buildrel 0\over{D}$ is a differential operator in $t_1,t_2$,
the radial  part of $D$ with respect to $N$. It turns out to be a
differential operator with constant coefficients. Introduce 
operators
\begin{eqnarray}
D_1&=&\left[\frac\partial{\partial t_1}+\frac\partial{\partial
t_2}+n-3\right]^2-(2n-7)\nonumber\\
D_2&=&\left[\frac\partial{\partial t_1}-\frac\partial{\partial
t_2}+1\right]^2-(2n-7).\nonumber
\end{eqnarray}

\begin{theorem}
We have that
\begin{eqnarray}
\buildrel 0\over{\Delta_2}&=&\frac12\left\{{D_1+D_2-(n-4)(n-6)}\right\},\nonumber \\
\buildrel 0\over{\Delta_4}&=&D_1D_2+2(n-4)^3.\nonumber
\end{eqnarray}
\end{theorem}

\section{Polynomial quantization on $G/H$}

We follow the scheme from [2]. The role of supercomplete system
is played by the kernel ${\Phi (\xi ,\eta )= \Phi _{\sigma
,\varepsilon }(\xi ,\eta )=N(\xi,\eta)^{\sigma,\,\varepsilon}}$ of
the intertwining operator ${A_{2-n-\sigma,\,\varepsilon}}$. As an
analogue of the Fock space, we 
%Old
%regard
%New
take
 the space of functions
$\varphi (\xi )$. We start from the algebra of operators
${D=T_{\sigma,\,\varepsilon}(X)}$, where $X$ belongs to the
universal enveloping algebra ${{\rm Env}(\got g)}$ for $\got g$.
{\it The covariant symbol} ${F(\xi,\eta)}$ of the operator $D$ is
defined by:
$$
F(\xi ,\eta )=\frac 1{\Phi (\xi ,\eta )}\,D_{\xi}\Phi (\xi ,\eta ),
$$
where ${D_{\xi}}$ means that the operator $D$ acts on $\Phi (\xi
,\eta )$ as on a function of $\xi $. These covariant symbols are
independent of $\varepsilon $. They are functions on $G/H$.
Moreover, they are {\it polynomials} on $G/H$ (i.e. restrictions on
$G/H$ of polynomials on the space of matrices $z$,
see~(\ref{matr-z})).

For generic $\sigma$ the space ${{\cA}_\sigma}$ of covariant
symbols is the space $S(G/H)$ of all polynomials on $G/H$.

The map ${D\mapsto F}$, which assigns to an operator its covariant
symbol, is $\got g$--equivariant. For an arbitrary $\sigma$ the
operator $D$ is reconstructed from its covariant symbol $F$:
\begin{equation}
\label{cov}
(D\varphi )(\xi )=c(\sigma,\varepsilon)\int F(\xi ,v)\,\frac {\Phi (\xi ,v)}{\Phi (u,v)}\, \varphi
(u)\,dx(u,v),
\end{equation}
where ${c(\sigma,\varepsilon)=\omega_0(\sigma,\varepsilon)^{-1}}$.

The multiplication of operators gives rise to a multiplication
(denote it by $*$) of covariant symbols. Let $F_1,\ F_2$ be the
covariant symbols of operators $D_1,\ D_2$ respectively. We have that
$$
F_1*F_2=\frac1{\Phi}(D_1)_{\xi}(\Phi F_2).
$$
This multiplication is given by 
%Old
%an integral:
$$ (F_1*F_2)(\xi ,\eta
)=\int F_1(\xi ,v)F_2(u,\eta ){\cB}(\xi ,\eta ; u,v)\,dx(u,v),
$$
where $dx(u,v)$ is an invariant measure on $G/H$, and
$$
{\cB}(\xi ,\eta ;u,v)=c\,\frac {\Phi (\xi ,v) \Phi (u,\eta
)}{\Phi(\xi ,\eta )\Phi (u,v)}.
$$
Let us call this kernel ${\cB}$ the {\it Berezin kernel}.

Thus, the spaces ${{\cA}_\sigma}$ turn out to be associative
algebras with unit (with respect to $*$).

On the other hand, we can define {\it contravariant symbols} of the
operators. 
%Old
%A function $F(\xi ,\eta )$ is the contravariant symbol
%for the following operator $A$ (acting on functions $\varphi (\xi
%)$):
%New
A function $F(\xi,\eta)$ can be viewed as the contravariant symbol
for the following operator $A$ (acting on functions $\varphi (\xi
)$):
$$
(A\varphi )(\xi )=c(\sigma,\varepsilon )\int
F(u,v)\,\frac {\Phi (\xi ,v)}{\Phi (u,v)} \varphi (u)\,dx(u,v).
$$
Notice that this expression differs from (\ref{cov}) only by the first argument of
function $F$. A contravariant symbol can be reconstructed from the
corresponding operator.

Thus we obtain two maps $D\mapsto F$ ("co") and $F\mapsto A$
("contra"), connecting operators $D$ and $A$ with polynomials $F$ on
$G/H$.

The passage from the contravariant symbol of an operator to its
covariant symbol is an integral operator with the Berezin kernal.
Let us call ${\cB}$ the {\it Berezin transform}.

\begin{theorem} The Berezin transform can be expressed in terms of Laplace operators:
$$
{\cB}=\frac{\Gamma (\sigma +n-2+\frac{a+b}2)\Gamma (\sigma
+1-\frac{a+b}2)\Gamma (\sigma +\frac{n}2+\frac{a-b}2)\Gamma (\sigma
+\frac{n}2-1-\frac{a-b}2)}{\Gamma (\sigma +n-2)\Gamma (\sigma
+1)\Gamma (\sigma +\frac{n}2)\Gamma (\sigma +\frac{n}2-1)}
$$
%Overfull 5pt!
%????????????
where ${a,b}$ are some variables and one has to consider
$$
D_1=(a+b)^2+2(n-3)(a+b)+(n-4)^2,
$$
$$
D_2=(a-b)^2+2(a-b)-2(n-4).
$$
%????????????
\end{theorem}

Note that on finite-dimensional subspaces in ${S(G/H)}$ the Berezin transform
is a differential operator.

Now let ${\sigma \to -\infty}$. The first two terms of the asymptotic
expansion of ${\cB}$ are given by:
\markboth{Cormac Walsh}{Polynomial quantization with pseudo-orthogonal group}

\begin{equation}
\label{asim}
{\cB}\sim 1-\frac 1\sigma\,\Delta_2.
\end{equation}
The relation (\ref{asim}) implies the following {\it correspondence principle}
(as the "Planck constant"  one has to take $h=-1/\sigma$):
\begin{equation}
\label{umn}
F_1*F_2\longrightarrow F_1F_2,
\end{equation}
\begin{equation}
\label{skob}
-\sigma \,(F_1*F_2-F_2*F_1)\longrightarrow
\{F_1,\,F_2\},
\end{equation}
as $\sigma \to -\infty $,
%Old
%In the right-hand sides of (\ref{umn})
%and (\ref{skob}) the pointwise multiplication and the Poisson
%bracket stand respectively.
%New
In (\ref{umn}) and (\ref{skob}),  $F_1 F_2$ denotes the pointwise multiplication of $F_1$ and $F_2$, and $\{F_1, F_2\}$ stands
for the Poisson bracket of $F_1$ and $F_2$.

\setcounter{section}{0}
\setcounter{footnote}{0}
\nachaloe{Cormac Walsh}{The horofunction boundary}{This work was funded in part by grant RFBR/CNRS 05-01-02807.}
\label{wal-abs}
The horofunction boundary (also known as the `metric' or `Busemann' boundary)
is a means of compactifying metric spaces. Its definition goes back to
Gromov~\cite{gromov:hyperbolicmanifolds}
in the 1970s but it seems not to have received much study until recently,
when it has appeared in several different
domains~\cite{ishii_mitake,rieffel_group,newman,AGW-m,karlsson}.
To define this boundary for a metric space $(X,d)$, one assigns
to each point $z\in X$ the function $\phi_z:X\to \R$,
\begin{equation*}
\phi_z(x) := d(x,z)-d(b,z),
\end{equation*}
where $b$ is some basepoint.
If $X$ is proper, then the map
$\phi:X\to C(X),\, z\mapsto \phi_z$ defines an embedding of $X$ into $C(X)$,
the space of continuous real-valued functions on $X$ endowed
with the topology of uniform convergence on compacts.
The horofunction boundary is defined to be
$X(\infty):=\closure\{\phi_z\mid z\in X\}\backslash\{\phi_z\mid z\in X\}$,
and its elements are called horofunctions.

This boundary is not the same as the better known Gromov boundary of
a $\delta$-hyperbolic space. For these spaces, it has been
shown~\cite{coornaert_papadopoulos_horofunctions,winweb_hyperbolic,
storm_barycenter} that the horoboundary is finer than the Gromov boundary
in the sense that there exists a continuous surjection from
the former to the latter.

Of particular interest are those horofunctions that are the limits of
almost-geodesics.
An almost-geodesic, as defined by Rieffel~\cite{rieffel_group},
is a map $\gamma$ from an unbounded set $T\subset \R_+$ containing 0 to $X$,
such that for any $\epsilon>0$,
\begin{equation*}
|d(\gamma(t),\gamma(s))+d(\gamma(s),\gamma(0))-t| < \epsilon
\end{equation*}
for all $t\in T$ and $s\in T$ large enough with $t\ge s$.
Rieffel called the limits of such paths Busemann points.
See~\cite{AGW-m} for a slightly different definition of almost-geodesic
which nevertheless gives rise to the same set of Busemann points.

As noted by Ballmann~\cite{ballmann:spaces},\markboth{Cormac Walsh}{The horofunction boundary}
the construction above is an additive analogue of the way the Martin
boundary is constructed in Probabilistic Potential Theory.
One may pursue the analogy further
in the framework of max-plus algebra, where one replaces the usual
operations of addition and multiplication by those of maximum and addition.
Indeed, this approach has already provided inspiration for many results
about the horofunction boundary~\cite{AGW-m,walsh}.
We mention, for example, the characterisation
of Busemann points as the functions in the horoboundary that are extremal
generators in the max-plus sense of the set of 1-Lipschitz functions.
So the set of Busemann points is seen to be an
analogue of the \emph{minimal} Martin boundary.
There is also a representation of 1-Lipschitz functions in terms of
horofunctions analogous to the Martin representation theorem.

%\providecommand{\doi}[1]{}
%\def\eqref#1{(\ref{#1})}
%\newcommand{\dd}{\,\mrm{d}}
%\newcommand{\mbuu}{\sM_{\mrm{bu}}}
%\newcommand{\mrm}[1]{\text{\rm #1}}
%\newcommand{\cR}{\mathcal{R}} 
%\newcommand{\sB}{\mathscr{B}}
%\newcommand{\sK}{\mathscr{K}}
%\newcommand{\sH}{\mathscr{H}}
%\newcommand{\sM}{\mathscr{M}}
%\newcommand{\sD}{\mathscr{D}}
%\newcommand{\sS}{\mathscr{S}}
%\newcommand{\sF}{\mathscr{F}}
%\newcommand{\sFb}{\bar{\sF}}
%\newcommand{\sC}{\mathscr{C}}
%\newcommand{\sMin}{\sM^{m}}
%\newcommand{\sMinp}{\sM'{}^{m}}
%\newcommand{\set}[2]{\{#1\mid\,#2\}}
%\newcommand{\Zb}{\overline{\Z}}
%\newcommand{\rmax}{\R_{\max}}
%\newcommand{\rmaxb}{\overline{\R}_{\max}}
%\newcommand{\rmaxbs}{\overline{\R}_{\max}^{_{\scriptstyle S}}}
%\newcommand{\rmaxbss}{\overline{\R}_{\max}^{_{\scriptstyle S\times S}}}
%\newcommand{\new}[1]{{\em #1}\index{#1}}
%\newcommand{\ind}[1]{\chi_{#1}}
%\newcommand{\tr}{\mrm{tr}\,}
%\newcommand{\pl}{\flat}
%\newcommand{\access}{\stackrel{*}{\to}}
%% defining 0,1 font 
%\newcommand{\bbfamily}{\fontencoding{U}\fontfamily{bbold}\selectfont}
%\DeclareMathAlphabet{\mathbbold}{U}{bbold}{m}{n}
%\newcommand{\zero}{\mathbbold{0}}
%\newcommand{\unit}{\mathbbold{1}}
%\newcommand{\abar}{\bar A}

There are few examples of metric spaces where the horofunction
boundary or Busemann points are explicitly known.
The first cases to be investigated
were those of Hadamard manifolds~\cite{ballmann:manifolds}
and Hadamard spaces~\cite{ballmann:spaces}, where the horofunction
boundary turns out to be homeomorphic to the ray boundary and all
horofunctions are Busemann points.
The case of finite-dimensional normed spaces has also received attention.
Karlsson \textit{et.~al.}~determined the horofunction boundary in the case when
the norm is polyhedral~\cite{karl_metz_nosk_horoballs}.
Other examples of metric spaces where the horofunction
boundary has been studied include the Cayley graphs of
finitely-generated abelian groups, studied by
Develin~\cite{develin_cayley},
and Finsler $p$--metrics on $\text{GL}(n,\mathbb{C})/\text{U}_n$,
where explicit expressions for the horofunctions were found by
Friedland and Freitas~\cite{friedland_freitas_pmetrics,
friedland_freitas_revisiting1}.
Webster and Winchester have some general results on when all horofunctions
are Busemann points~\cite{winweb_busemann},~\cite{winweb_metric}.

In the following sections, we describe our recent work elucidating
the hor\-obo\-und\-ary of some particular metric spaces.

\section{Normed spaces}

%In~\cite{gromov:hyperbolicmanifolds}, Gromov defines a boundary of a metric

Rieffel comments that it is an interesting question as to when
all boundary points of a metric space are Busemann points and asks whether
this is the case for general finite-dimensional normed spaces.
In~\cite{walsh:normed}, we answer this question in the negative and give
a necessary and sufficient criterion for it to be the case.

Let $V$ be an arbitrary finite-dimensional normed space with unit ball $B$.
Recall that a convex subset $E$ of a convex set $D$ is said to be an
\emph{extreme set} if the endpoints of any line segment in $D$ are
contained in $E$ whenever any interior point of the line segment is.
For any extreme set $E$ of the dual unit ball $\dualball$
and point $p$ of $V$, define the function
$\hnorm_{E,p}$ from the dual space $\dualspace$ to $[0,\infty]$ by
\begin{align*}
\hnorm_{E,p}(q):= \indicator_E(q) + \dotprod{q}{p} - \inf_{y\in E}\dotprod{y}{p}
\qquad\text{for all $q\in V^*$}.
\end{align*}
Here $\indicator_E$ is the indicator function, taking value $0$ on $E$
and $+\infty$ everywhere else.

Our first theorem characterises the Busemann points of $V$
as the Legendre-Fenchel transforms of these functions.
\begin{theorem}
\label{thm:theorem1}
The set of Busemann points of a finite-dimensional normed space $(V,||\cdot||)$
is
\begin{equation*}
\{ \hnorm^*_{E,p}  \mid \text{$E$ is a proper extreme set of $\dualball$
and $p\in V$} \}.
\end{equation*}
\end{theorem}

We use this knowledge to characterise those norms for which all horofunctions
are Busemann points.
\begin{theorem}
\label{thm:theorem2}
A necessary and sufficient condition for every horofunction
of a finite-dimensional normed space to be a Busemann point
is that the set of extreme sets of the dual unit ball be closed in the
Painlev\'e--Kuratowski topology.
\end{theorem}

\section{The Hilbert metric}

Let $x$ and $y$ be distinct points in a bounded open convex subset
$D$ of $\R^N$, with $N\ge 1$.
Define $w$ and $z$ to be the points in the Euclidean boundary of $D$
such that $w$, $x$, $y$, and $z$ are collinear and arranged in this
order along the line in which they lie.
The Hilbert distance between $x$ and $y$ is
defined to be the logarithm of the cross ratio of these four points:
\begin{equation*}
\hil{x}{y}{}:= \log \frac{|zx|\,|wy|}{|zy|\,|wx|}.
\end{equation*}
If $D$ is the open unit disk, then the Hilbert metric is exactly
the Klein model of the hyperbolic plane.

As pointed out by Busemann~\cite[p105]{busemann:geometry},
the Hilbert geometry is related to hyperbolic geometry in much the same way
that normed space geometry is related to Euclidean geometry.
It will not be surprising therefore that there are similarities between
the horofunction boundaries of Hilbert geometries and of normed spaces.

Define the function
\begin{equation}
\label{eqn:funkratio}
\funk{x}{y}{}:= \log \frac{|zx|}{|zy|},
\qquad\text{for all $x$ and $y$ in $D$.}
\end{equation}
This function satisfies the usual metric space axioms,
apart from that of symmetry.

Hilbert's metric can now be written
\[
\hil{\poynta}{\poyntb}{}
   := \funk{\poynta}{\poyntb}{}+\funk{\poyntb}{\poynta}{},
\qquad\text{for all $\poynta$ and $\poyntb$ in $D$}.
\]
This expression of the Hilbert metric as the symmetrisation
of the Funk metric plays a crucial role.
It turns out that every Hilbert horofunction is the sum of a
horofunction in the Funk geometry and a horofunction in the reverse Funk
geometry, where the metric in the latter is given by
\begin{equation*}
\rev{x}{y}{}:=\funk{y}{x}{}.
\end{equation*}
This allows us to simplify the problem by investigating separately
the horofunction boundaries of these two geometries and then combining the
results. Determining the boundary of the Funk geometry turns out to be very
similar to determining that of a normed space,
which was done in~\cite{walsh:normed}.

In~\cite{walsh_hilbertboundary},
we characterise those Hilbert geometries for which all
horofunctions are Busemann points.
\begin{theorem}
\label{theorem3}
A necessary and sufficient condition for every horofunction on a
bounded convex open subset of $\R^N$ containing the origin
to be a Busemann point in the Hilbert geometry is that
the set of extreme sets of its polar
be closed in the Painlev\'e--Kuratowski topology.
\end{theorem}

It had previously been shown~\cite{karl_metz_nosk_horoballs}
that all horofunctions of the Hilbert geometry on a polytope
are Busemann points.

\begin{theorem}
\label{theorem2}
Let $D$ be a bounded convex open subset of $\R^N$.
If a sequence in $D$ converges to a point in the horofunction boundary of the
Hilbert geometry, then the sequence converges in the usual sense to a
point in the Euclidean boundary~$\partial D$.
\end{theorem}

\section{Finitely generated groups}

An interesting class of metric spaces are the Cayley graphs of finitely
generated groups with their word metric.
Here one may hope to have a combinatorial description of the horoboundary.

The first to consider the horoboundary in this setting was
Rieffel~\cite{rieffel_group} who studied the horoboundary of $\Z^n$ with an
arbitrary finite generating set in connection with his work on
non-commutative geometry.

In~\cite{walsh_artin}, we investigate the horofunction boundary of
Artin groups of dihedral type.
Let $\prodd(s,t;n):= ststs\cdots$, with $n$ factors in the product.
The Artin groups of dihedral type have the following presentation:
\begin{align*}
\dihed = \langle a,b \mid \prodd(a,b;k)=\prodd(b,a;k) \rangle,
\qquad\text{with $k\ge3$.}
\end{align*}
Observe that $A_3$ is the braid group on three strands.
The generators traditionally considered are the Artin generators
$\gens:=\{a,b,a^{-1},b^{-1}\}$.

In what follows, we will have need of the Garside normal form
for elements of $\dihed$.
The element $\Delta:=\prodd(a,b;k)=\prodd(b,a;k)$ is called the
Garside element.
Let
\begin{align*}
M^+:= \{a,b,ab,ba,\dots,\prodd(a,b;k-1),\prodd(b,a;k-1)\}.
\end{align*}
It can be shown~\cite{word_processing} that $w\in \dihed$ can be written
\begin{align*}
w = w_1\cdots w_n \Delta^r
\end{align*}
for some $r\in\Z$ and $w_1,\dots,w_n\in M^+$.
This decomposition is unique if $n$ is required to be minimal.
We call it the right normal form of $w$.
The factors $w_1,\dots,w_n$ are called the canonical factors of $w$.

An algorithm was given in~\cite{mairesse_matheus_growth} for
finding a geodesic word representing any given element of $A_k; k\ge 3$.
We use this algorithm to find a simple formula for the word length metric.
\begin{proposition}
Let $x=z_1\cdots z_m \Delta^{r}$ be an element of $\dihed$ written in right
normal form. Let $(p_0,\dots,p_{k-1})\in\N^{k}$ be such that
$p_0:=r$ and, for each $i\in\{1,\dots,k-1\}$, $p_i-p_{i-1}= m_{k-i}$,
where $m_i$ is the number of canonical factors of $x$ of length $i$.
Then the distance from the identity $e$ to $x$ in the Artin-generator
word-length metric is
\begin{align*}
\dist(e,x)=\sum_{i=0}^{k-1} |p_i|.
\end{align*}
\end{proposition}
Since $d$ is invariant under left multiplication, that is,
$d(y,x)=d(e,y^{-1}x)$, we can use this formula to calculate the distance
between any pair of elements $y$ and $x$ of $\dihed$.
With this knowledge we can find the following description of the horofunction
compactification.
{\sloppy

}
\newcommand\emme{m}

Let $Z$ be the set of possibly infinite words of positive generators having
no product of consecutive letters equal to $\Delta$.
We can write each element $z$ of $Z$ as a concatenation of substrings in
such a way that the products of the letters in every substring equals
an element of $M^+$ and the combined product of letters in each consecutive
pair of substrings is not in $M^+$. Because $z$ does not contain $\Delta$,
this decomposition is unique. Let $\emme_i(z)$ denote the number of substrings
of length $i$. Note that if $z$ is an infinite word, then this number will
be infinite for some $i$.

Let $\Omega'$ denote the set of $(p,z)$ in
$(\Z\union\{-\infty,+\infty\})^k\times Z$
satisfying the following:
\begin{itemize}
\item
$p_i-p_{i-1}\ge \emme_{k-i}(z)$ for all $i\in\{1,\dots,k-1\}$ such that
$p_i$ and $p_{i-1}$ are not both $-\infty$ nor both $+\infty$;
\item
if $z$ is finite, then $p_i-p_{i-1} = \emme_{k-i}(z)$
for all $i\in\{1,\dots,k-1\}$ such that
$p_i$ and $p_{i-1}$ are not both $-\infty$ nor both $+\infty$.
\end{itemize}
We take the product topology on $\Omega'$.

We now define $\Omega$ to be the quotient topological space of $\Omega'$
where the elements of $(+\infty,\dots,+\infty)\times Z$ are considered
equivalent and so also are those in $(-\infty,\dots,-\infty)\times Z$.
We denote these two equivalence classes by $\plusclass$ and $\minusclass$,
respectively.

We let $\mathcal{M}$ denote the horofunction compactification of $\dihed$
with the Artin-generator word metric.
The basepoint is taken to be the identity.

\begin{theorem}
The sets $\Omega$ and $\mathcal{M}$ are homeomorphic.
\end{theorem}
Let $Z_0$ be the set of elements of $Z$ that are finite words.
Let $\Omega_0$ denote the set of $(p,z)$ in $\Z^k\times Z_0$
such that $p_i-p_{i-1}= \emme_{k-i}(z)$ for all $i\in\{1,\dots,k-1\}$.
One can show that the elements of $\Omega_0$ are exactly the elements of
$\Omega$ corresponding to functions of the form $d(\cdot,z)-d(e,z)$
in $\mathcal{M}$.

In the present context, since the metric takes only integer values,
the Busemann points are exactly the limits of geodesics
(see~\cite{winweb_busemann}).
Develin~\cite{develin_cayley}, investigated the horoboundary of finitely
generated abelian groups with their word metrics and showed that all their
horofunctions are Busemann.

We have the following characterisation of the Busemann points of $\dihed$.
\begin{theorem}
A function in $\mathcal{M}$ is a Busemann point if and only if the
corresponding element $(p,z)$ of $\Omega$ is in $\Omega\backslash\Omega_0$
and satisfies the following:
$p_i-p_{i-1} = \emme_{k-i}(z)$ for every $i\in\{1,\dots,k-1\}$ such that
$p_i$ and $p_{i-1}$ are not both $-\infty$ nor both $+\infty$.
\end{theorem}

The group $A_k$ also has a dual presentation:
\begin{align*}
A_k = \langle \sigma_1,\dots,\sigma_k
 \mid \sigma_1\sigma_2 = \sigma_2\sigma_3 = \cdots = \sigma_k\sigma_1 \rangle,
\qquad\text{with $k\ge3$.}
\end{align*}
The set of dual generators is
$\dualgens := \{\sigma_1,\dots,\sigma_k,\sigma_1^{-1},\dots,\sigma_k^{-1}\}$.

Again, one can find a formula for the word length metric and
use it to determine the horoboundary. This time however,
it turns out that there are no non-Busemann points.
\begin{theorem}
In the horoboundary of $\dihed$ with the dual-generator word metric,
all horofunctions are Busemann points.
\end{theorem}

In general, one would expect the properties of the horoboundary
of a group with its word length metric to depend strongly on the generating
set. It would be interesting to know for which groups and for which properties
there is not this dependence.
As already mentioned, all boundary points of abelian groups are
Busemann no matter what the generating set~\cite{develin_cayley}.
On the other hand, the above results show that
for Artin groups of dihedral type
the existence of non-Busemann points depends on the generating set.
\markleft{Evgeny M.~Beniaminov}%{The horofunction boundary}

\selectlanguage{russian}
%\newpage
\setcounter{section}{0}
\setcounter{footnote}{0}
\nachalo{Е.М.~Бениаминов}{Квантование как приближенное
описание некоторого диффузионного процесса}
%\markboth{Е.М.~Бениаминов}{Квантование как приближенное описание}
\label{ben-rus-abs}
\section{Описание и некоторые свойства модели}
Рассматривается некоторая математическая модель процесса, состояние которого 
в каждый момент времени задается волновой функцией -- комплекснозначной функцией 
$\varphi (x, p)$, где  $(x,p) \in R^{2n},$ 
и $n$  ---  размерность конфигурационного пространства. 
 В отличие от квантовой механики, где волновая функция зависит только от координат или 
только от импульсов, в нашем случае волновая функция зависит и от координат и и от импульсов. 
Так же, как в квантовой механике, предполагается, что для волновых функций выполняется 
принцип суперпозиции, и плотность вероятности 
$\rho (x,p)$ на фазовом пространстве, соответствующая  волновой функции 
$\varphi(x, p),$ задается стандартной формулой
\begin{equation}\label{rho4668}
\rho (x,p)= \varphi^*(x, p) \varphi(x, p)=|\varphi(x, p)|^2. 
\end{equation}

В работе рассматривается   классическая  модель диффузионного процесса для волновой 
функции $\varphi (x, p)$ на  фазовом пространстве. Предполагается, что каждый комплексный вектор 
волновой функции одновременно находится в 4-х движениях: 
\markboth{Evgeny M.~Beniaminov}{Quantization as approximate description
(in Russian)}
{\tolerance=500

}
точка приложения вектора движется по классической траектории, заданной функцией Гамильтона
$H(x, p)$;

точка приложения вектора перемещается случайно по координатам и импульсам, находясь 
в диффузионном процессе с постоянными коэффициентами диффузий $a^2$ и $b^2$ 
по координатам и  импульсам, соответственно;

точка приложения каждого вектора движется по случайной траектории в результате 
движений, описанных в двух предыдущих пунктах, а сам вектор вращается с постоянной
угловой скоростью $\omega ={mc^2}/{\hbar}$ в  системе координат, связанной с этой точкой, где
$m$ -- масса частицы, $c$ -- скорость света, ${\hbar}$ -- постоянная Планка;  

длина всех комплексных векторов волновой функции в момент времени $t$ умножается на
$\exp (abnt/\hbar)$ (это чисто техническое требование, которое не сказывается на относительных 
вероятностях нахождения частицы в фазовом пространстве). 
 
Предполагается, что волновой вектор $\varphi (x, p, t)$
в точке $(x, p)$ в момент времени $t$ по принципу суперпозиции равен сумме волновых векторов, 
заданных распределением 
векторов $\varphi^0 (x, p)$ в начальный момент времени  и попавших в результате описанных выше 
движений в точку $(x, p)$ в момент времени $t$.

Процесс описывается  дифференциальным уравнением диффузионного типа. 
 Анализ дифференциального
 уравнения модели показывает, что движение в модели раскладывается на быстрое и медленное.
В результате быстрого движения система, начиная с 
произвольной волновой функции на фазовом пространстве, переходит к функции, 
принадлежащей некоторому особому подпространству. Элементы 
этого подпространства параметризуются  волновыми функциями, зависящими только от 
координат. Медленное движение по подпространству  описывается уравнением Шредингера.

Исходя из предположений о тепловой причине диффузий и соответствии 
следствий модели известным физическим экспериментам  Лэмба - Резерфорда \cite{lamb}
(сдвиг Лэмба в спектре атома водорода), в работе делается оценка коэффициентов 
диффузий и времени переходного процесса от классического описания процесса, в котором принцип 
неопределенности Гейзенберга может не выполняться, к квантовому, в котором принцип Гейзенберга
 уже выполняется. Время переходного процесса имеет порядок 
$1/T \cdot 10^{-11} \text{с}$, где 
$T$ --- температура среды.

\section{Основные результаты}
Рассмотрим диффузионный процесс  на  фазовом пространстве, в котором волновая функция
 $\varphi (x, p, t)$ в момент времени $t$  удовлетворяет 
дифференциальному  уравнению
\begin{equation}\label{eq_diff}
\frac{\partial\varphi}{\partial{t}}=\sum_{k=1}^{n}
\biggl(
\frac{\partial H}{\partial x_k} \frac{\partial\varphi}{\partial p_k}-
\frac{\partial H}{\partial p_k} \frac{\partial\varphi}{\partial x_k}
\biggr)
-\frac{i}{\hbar}
\biggl(H-\sum_{k=1}^{n}\frac{\partial H}{\partial p_k}p_k\biggr)\varphi
+\Delta_{a,b}{\varphi},
\end{equation}
\begin{equation}\label{delta}
\mbox{где }\ \ \ \ \ \Delta_{a,b}{\varphi}=
a^2\sum_{k=1}^{n}\biggl(\frac{\partial}{\partial{x_k}}-
     \frac{ip_k}{\hbar}\biggr)^{2}\varphi
+b^2\sum_{k=1}^{n}\frac{\partial^2 }{\partial{p^2_k}}\varphi
     +\frac {abn}{\hbar}{\varphi},
\end{equation}
где
$H(x, p)$ --- функция Гамильтона;
$a^2$ и $b^2$ --- коэффициенты диффузий по координатам и  импульсам, соответственно.

Если в уравнении~(\ref{eq_diff}) отбросить последнее слагаемое, то получим 
дифференциальное уравнение в частных производных 
первого порядка. Эта  часть
уравнения~(\ref{eq_diff}) описывает детерминированную
составляющую 
движения
комплексных векторов $\varphi(x, p,t)$.  Согласно уравнению, в этом движении точка
приложения каждого вектора движется по классической траектории, заданной 
гамильтонианом $H(x, p)$, а сам вектор при этом вращается в каждой точке траектории
с угловой скоростью
{\tolerance=500

}
\begin{equation}
\omega'= \frac{1}{\hbar}\biggl(H-\sum_{k=1}^{n}
\frac{\partial H}{\partial p_k}p_k\biggr).
\end{equation}
Заметим, что в случае, когда конфигурационное пространство  трехмерно и 
$H= c\sqrt {m^2c^2 +p^2}$,  то
$\omega' dt  =
\frac{mc^2}{\hbar}\frac{mc^2dt}{H}= \frac{mc^2}{\hbar} d\tau ,$
где $\tau$ --- собственное время в системе координат, связанной с частицей, движущейся
с импульсом $p$.  То есть в этом случае, вектор, точка приложения которого движется
по классической траектории, вращается с постоянной угловой скоростью 
$\omega ={mc^2}/{\hbar}$ в  системе координат, связанной с этой точкой.

Наоборот, если в правой части уравнения~(\ref{eq_diff}) оставить только последнее слагаемое
вида~(\ref{delta}), то получим уравнение 
\begin{equation}\label{eq_delta}
\frac{\partial\varphi}{\partial{t}}=
a^2\sum_{k=1}^{n}\biggl(\frac{\partial}{\partial{x_k}}-
     \frac{ip_k}{\hbar}\biggr)^{2}\varphi
+b^2\sum_{k=1}^{n}\frac{\partial^2 }{\partial{p^2_k}}\varphi
     +\frac {abn}{\hbar}{\varphi}.
\end{equation}
Это уравнение описывает диффузионную составляющую движения 
векторов $\varphi(x, p,t)$ на фазовом пространстве. В этом движении точки приложения 
векторов перемещаются в соответствии с классическим однородным диффузионным 
процессом с коэффициентами диффузий по координатам и  импульсам равными $a^2$ и $b^2$,
соответственно. При этом сам вектор при малых случайных перемещениях из точки $(x, p)$  
в точку $(x+dx, p+dp)$ переносится параллельно, 
а его длина  в момент времени $t$ умножается на $\exp (abnt/\hbar)$. Заметим, что 
параллельный перенос векторов на фазовом пространстве задается связностью, которая
выражается формулой:
$
L_{(dx,dp)}\varphi(x, p) - \varphi(x, p)\approx  -({i}/{\hbar})\varphi(x, p) p dq,
$
где $L_{(dx,dp)}\varphi(x, p)$ --- параллельный перенос вектора $\varphi(x, p)$ из точки
$(x, p)$  по бесконечно малому вектору $(dx, dp)$.
В частном случае, когда конфигурационное пространство трехмерно, такая связность на 
фазовом пространстве вызвана синхронизацией движущихся часов в точках фазового
пространства.

Правая часть уравненния~(\ref{eq_delta}) --- самосопряженный оператор. Задача на 
собственные  значения для этого оператора преобразованием Фурье по координатам 
сводится к стационарному уравнению Шредингера  для гармонических колебаний. 
Отсюда показывается, что собственные значения оператора 
уравнения~(\ref{eq_delta}) неположительны, и верна  следующая теорема.

\begin{teorema}
Пусть $\varphi(x,p,0)$ --- произвольная функция, преобразование
Фурье которой по $p$ стремится к нулю при $x \rightarrow  \infty $. 
Тогда решение $\varphi(x,p,t)$ диффузионного уравнения~(\ref{eq_delta}) 
 экспоненциально  по времени (с показателем равным $-abt/{\hbar}$) стремится к 
стационарному решению вида: 
\begin{align}
\varphi(x,p)&=\lim_{t \to \infty} \varphi(x,p,t)=
\frac{1}{(2\pi{\hbar})^{n/2}}
\!\int\limits_{R^n}\!\!\psi(y)\chi(x,y)
e^{-{{i  (y-x)p}/{\hbar}}}
dy,\label{view_varphi}\\
\intertext{где}
\psi(y) &=\frac{1}{(2\pi{\hbar})^{n/2}}
\int\limits_{R^{2n}}\!\!\varphi(x,p,0)
e^{{{i  (y-x)p}/{\hbar}}}
\chi(x,y)dpdx,\label{psi}\\
\chi(x,y)&=\left(\frac{b}{a\pi\hbar}\right)^{n/4}e^{-{{b}(x-y)^2}/{(2a\hbar)}}.\label{chi}
\end{align}
\end{teorema}

Заметим, что   $\chi^2(x,y)$   представляет   собой   плотность вероятностей
нормального   распределения   по   $x$  с  математическим
ожиданием $y$, и дисперсией  $a\hbar/(2b)$. Если величина
$a\hbar/(2b)$ мала, то функция $\chi^2(x,y)$ близка к дельта-функции от
$x-y$.

Композиция выражений~(\ref{psi}) и (\ref{view_varphi}) строит 
проектор из пространства всех волновых функций, заданных на фазовом пространстве, 
на некоторое подпространство.
Элементы этого подпространства параметризуются функциями вида $\psi(y),$ где $y\in R^n$,  
т.~е. волновыми функциями на конфигурационном пространстве. 

Если же предполагать, что диффузия вызывается тепловыми воздействиями на электрон, то  
коэффициенты диффузий по координатам и импульсам
выражаются в статистической физике 
(см., например \cite{isihara}, гл.7, \S 4 и \S 9) 
через температуру $T$ по формулам:
$a^2= kT/(m\gamma) \ \  \text{      и     }\ \ b^2=\gamma kTm,$
 где $k$ --- постоянная Больцмана,  $m$ --- масса электрона, 
$\gamma$ --- 
коэффициент трения среды на единицу массы.
Отсюда, $a/b=(\gamma m)^{-1}$  и $ab=kT$. 
То есть, в этом случае, величина $a/b$, которая входит в 
выражение~(\ref{chi}), 
не зависит от температуры.
С другой стороны, $t$ --- время переходного процесса, определенное в теореме~1,
имет вид:
$t \sim {\hbar}/{(ab)}={\hbar}/{(kT)}= T^{-1}\cdot 7.638\cdot 10^{-12}\text{с}.$

С учетом этой оценки, будем считать в уравнении~(\ref{eq_diff})  величину ${\hbar}/{(ab)}$
малым параметром и предполагать, что координаты и импульсы мало меняются за
это время при классическом движении, определенном гамильтонианом $H(x, p)$.
{\tolerance=500

}
\begin{teorema}   
Движение, описываемое уравнением~(\ref{eq_diff}), асимптотически распадается при 
${\hbar}/{(ab)} \rightarrow   0$ на быстрое движение и медленное движение. 
В результате быстрого движения  произвольная волновая функция $\varphi(x, p, 0)$
переходит за время порядка ${\hbar}/{(ab)}$ к виду~(\ref{view_varphi}).  Волновые функции 
вида~(\ref{view_varphi}) образуют линейное подпространство.  
Элементы этого подпространства параметризуются  волновыми функциями  $\psi(y)$, зависящими 
только от координат $y\in R^n$. Медленное движение, начинающееся с ненулевой волновой функции из  
этого подпространства,   происходит по подпространству и параметризуется волновой функцией
$\psi(y,t)$, зависящей от времени. Функция $\psi(y,t)$ удовлетворяет уравнению Шредингера вида 
$ i\hbar {\partial \psi}/{\partial t} = \hat{H}\psi$,  где
{\tolerance=600

}
\begin{eqnarray}
\hat{H}\psi & = & \frac {1}{(2\pi \hbar)^n} \int \limits_{R^{3n}}  
\biggl (H(x, p)-\sum_{k=1}^{n}\biggl(\frac{ \partial H}{\partial x_k}
+\frac {i b}{a}\frac{ \partial H}{\partial p_k}\biggr)(x_k-y'_k)
\biggr)\times \nonumber\\
   &  & \times \chi(x, y) \chi (x, y')  e^{\frac{i}{\hbar}(y-y') p}  \psi(y',t) dy' dx dp, \nonumber  
\end{eqnarray}
и $\chi(x, y) $ задается формулой~(\ref{chi}).
\end{teorema}

\begin{teorema}
Если $\frac {a\hbar}{b}$ --- малая величина и 
$H(x, p) = \frac{p^2}{2m}+ V(x),$
то оператор $\hat H$ с точностью до членов порядка  $a\hbar/b$ имеет вид:
\begin{equation}\label{hatH}
\hat H \approx - \frac{\hbar^2}{2m}\biggl(\sum_{k=1}^{n}\frac{\partial^2 }{\partial{y^2_k}}\biggr)
+V(y)-\frac{a\hbar}{4b}\sum_{k=1}^{n}\frac{\partial^2 V}{\partial{y^2_k}} +\frac{3nb\hbar}{4ma}.
\end{equation}
\end{teorema}

Первые два слагаемые в формуле~(\ref{hatH}) дают стандартный оператор Гамильтона.
Последнее слагаемое ---  константа, и ею 
можно пренебречь. Предпоследнее слагаемое  рассмотрим (ввиду малости $a\hbar/b$) 
как возмущение к оператору Гамильтона. 
\markboth{A.M.~Gel'fand, B.Kh.~Kirshteyn}{Quantization as 
approximate description (in Russian)}

Считая, что отклонения в спектре атома водорода (сдвиг Лэмба), наблюдаемые в экспериментах 
Лэмба-Резерфода~\cite{lamb}, вызываются предпоследним слагаемым в формуле~(\ref{hatH}),
можно оценить величину $a/b$. Расчеты стандартным методом возмущений, аналогичные расчетам,
выполненным в \cite{welt}, дают следующую оценку: 
$a/b =3.41\cdot 10^4 \text{с/г}$. 
Отсюда, стандартное
отклонение для нормального распределения $\chi^2$, по
 которому производится сглаживание  
волновых функций, имеет вид 
$\sqrt{a\hbar/(2b)}=4.24 \cdot 10^{-12} \text{см}.$  
Эта величина существенно меньше радиуса атома водорода и близка к 
комптоновской длине волны электрона $\hbar/(mc)=3.86\cdot 10^{-11}
\text{см}$. 

Таким образом, расчеты показывают, что предложенная модель в виде дифференциального 
уравнения~(\ref{eq_diff}) достаточно адекватно
описывает физические процессы в стандартных случаях для стандартного гамильтониана. Но
эту модель можно применить и для расчетов процессов с нестандартным гамильтонианом или 
с гамильтонианом, быстро меняющимся во времени, как при внезапных возмущениях
или для периодически меняющегося потенциала с частотой порядка $ab/\hbar$, и сравнить с 
экспериментальными данными.

\renewcommand{\theequation}{\arabic{equation}}
\setcounter{equation}{0}
\setcounter{footnote}{0}
\nachaloe{А.М.~Гельфанд и Б.Х.~Кирштейн}{Идемпотентные 
системы нелинейных уравнений и задачи расчета 
электроэнергетических сетей}{Работа выполнена при частичной финансовой поддержке
гранта РФФИ 05-01-02807-НЦНИЛ\_а.}
%\markboth{А.М.~Гельфанд и Б.Х.~Кирштейн}{Идемпотентные системы
%уравнений в электроэнергетике}
\markboth{A.M.~Gel'fand, B.Kh.~Kirshteyn}{Idempotent equations in
electroenergetics (in Russian)}
\label{kir-rus-abs}
Электроэнергетическую сеть можно рассматривать как граф с $n$ вершинами, каждой вершине 
(узлу) которого сопоставлены два вещественных числа - активная $P_{k}$ и 
реактивная $Q_{k}$ составляющие инъекции мощности в узле $k$ $(k=1,...,n)$, и каждому ребру 
(линий электропередач), соединяющих $k$-ый и $j$-ый узлы - активная $Y_{kj}^{a}$ и реактивная 
$Y_{kj}^{r}$ составляющие проводимости.\\  
  Установившиеся режимы электроэнергетических сетей характеризуются значениями активных 
$U_{k}^{a}$ и реактивных $U_{k}^{r}$ составляющих напряжений в узлах, которые должны удовлетворять
системе $2n$ алгебраических уравнений \cite{E} - узловых уравнений балансов активной и реативных мощностей.\\  
  Такие уравнения можно рассматривать как вещественную и комплексную составляющую системы $n$ комплексных
уравнений вида
\begin{equation}\label{rho0}
 W_{k}=E_{k}(\sum_{j\in (k)} a_{kj}\overline{E}_{j}), 
\end{equation}
где\ 
$W_{k}= P_{k}+iQ_{k}$ , $E_{k}= U_{k}^{a}+iU_{k}^{r}$, $\overline{E}_{k}= U_{k}^{a}-iU_{k}^{r}$, 
суммирование идет по всем узлам $j$ связанных линией с узлом $k$, $a_{kj}$ комплексные числа, 
которые определяется по комплексным проводимостям $Z_{j}=U_{j}^{a}+iU_{j}^{r}$. Эти уравнения 
дают $2n$ вещественных алгебраических уравнений относительно $2n$ неизвестных 
$U_{k}^{a},  U_{k}^{r}, k=(1,...,n)$. \\
  Важной задачей анализа электроэнергетических сетей является анализ устойчивости (определение 
запаса устойчивости) установившегося режима. На практике, такой анализ сводится к 
проверке существовании вещественной деформации решения при деформациях коэффициентов уравнений, 
входящих в систему. Такие деформации отвечают изменениям инъекций мощности в узлах или измению
структуры графа сети (отключение линий). \\При этом в качестве критерия потери устойчивости при 
деформациях уравнений системы обычно принимается либо расходимость итерационного процесса 
метода Ньютона для нахождении решения, либо вырождение матрицы Якоби системы уравнений. Такой 
подход не всегда математически корректно отвечает поставленной залаче, но, главное, не позволяет 
быстро и наглядно получить границу области устойчивости установившегося режима.\\
Мы рассматриваем некоторую процедуру простого макетирования задачи анализа устойчивости электроэнергетических систем
с помощью анализа решений идемпотентной системы уравнений, полученной в результате деквантования по 
Маслову \cite {L}  исходной системы (1).\\
   Запишем систему уравнений установившегося режима более симметрично в виде системы $2n$
комплексных алгебраических уравнений \cite{M}, относительно $2n$ комплексных переменных. Для этого 
введем новые комплексные переменные $S_{k}$ вместо $\overline{E}_{k}$ и новые уравнения, полученные 
с помощью комплексного сопряженния уравнений (1) и последующей аналогичной заменной переменных. \\
Получим систему из $2n$ комплексных уравнений 
\begin{equation}\label{rho00}
W_{k}=E_{k}(\sum_{j\in (k)} a_{kj}S_{j}),\quad \overline{W}_{k}=S_{k}(\sum_{j\in (k)} \overline{a}_{kj}E_{j})
\end{equation}
%и
%\begin{equation}\label{rho000}
 %\overline{W}_{k}=S_{k}(\sum_{j\in (k)} \overline{a}_{kj}E_{j}) 
%\end{equation}
относительно $2n$ комплексных переменных $E_{k}, S_{k}, k=1,...,n$.   
Решение системы \eqref{rho00} удовлятворяет системе (1) тогда и только тогда когда
\begin{equation}\label{rho1}
 \overline{E}_{k}=S_{k}  
\end{equation}
для всех $k=1,...,n$.\\
В \cite{GKZ-Kirsh} для всякой алгебраической поверхности $f(z)=0$ в $\mathrm{C}^{2n}$ определена
ее амеба $\mathcal{A}_{f}$ в $\mathrm{R}^{2n}$, образ пересечения этой поверхности с комплексным 
тором относительно отображения
\[\mathrm{Log}:(z_{1},\ldots,z_{2n})\in{(\mathbf{C}\backslash \{0\})^{2n}}\to 
(\mathrm{log}\,|z_{1}|,\ldots,\mathrm{log}\,|z_{2n}|) \in \mathrm{R}^{2n} .\]
Условие \eqref{rho1} в этих терминах означает, что пересечение амеб всех уравнений системы \eqref{rho00} 
содержит точки, удовлетворяющие условиям
\begin{equation}\label{rho11}
 \mathrm{log}|E_{k}|=\mathrm{log}|S_{k}|  
\end{equation}
для всех $k=1,...,n$.\\
Известно \cite{GKZ-Kirsh}, что амеба $\mathcal{A}_{f}$ совпадает с дополненнием конечного числа открытых выпуклых 
подмножеств $\mathrm{E}_{\nu}$ в $\mathrm{R}^{2n}$: 
\[\mathrm{R}^{2n} \backslash \mathcal{A}_{f} = \cup\{\mathrm{E}_{\nu}\}.\]
В \cite{PR} определены наборы линейных функций на этих выпуклых множествах, нижняя грань которых
определяет кусочно-линейное подмножество $\mathrm{R}^{2n}$, которое называется спайном амебы,
лежит внутри $\mathcal{A}_{f}$ и является ее гомотопическим ретрактом.\\
Чтобы определить идемпотентную систему уравнений определим тропические аналоги операций сложения и 
умножения в $\mathrm{R}$ обычным образом: тропическое сложение как $x\oplus y = \max\{x,y\}$ и тропическое 
умножение как $x\otimes y = x+ y $.\\
Рассмотрим идемпотентную систему уравнений, полученную из уравнений \eqref{rho00} заменой
обычных операций на их тропические аналоги  и комплексных коэффициентов на логарифмы их модулей:
\begin{equation}\label{rho111}
 \mathrm{log}|W_{k}|=\mathrm{log}|E_{k}|\otimes(\bigoplus_{j\in (k)} \mathrm{log}|a_{kj}| \otimes \mathrm{log}|S_{j}|), 
\end{equation}
и
\begin{equation}\label{rho}
 \mathrm{log}|W_{k}|=\mathrm{log}|S_{k}|\otimes(\bigoplus_{j\in (k)} \mathrm{log}|a_{kj}| \otimes \mathrm{log}|E_{j}|), 
\end{equation}
Заметим, что многогранники Ньютона для уравнений \eqref{rho00}
совпадают с выпуклой обoлочкой подмножества вершин единичного куба в $\mathrm{R}^{2n}$ 
и не содержат поэтому внутри себя точек целочисленной решетки. \\
Как доказано в \cite{PT}, в этом случае уравнение 
спайна амебы  совпадают с уравнениями идемпотентной  системы, 
результата деквантования $f$. Отсюда легко получается следующая\\

\medskip
{\sc Теорема.} {\sl Решение идемпотентной системы уравнений \eqref{rho111}, \eqref{rho}
совпадает с множеством точек пересечения 
пределов амеб уравнений \eqref{rho00} при ретракции их на свои спайны.}

\medskip
Используя эту теорему можно следующим образом получить простую модель для анализа запаса устойчивости электроэнергетической системы.
\\Пусть 
\[E^{*} = (E_{1},...,E_{n})\in{(\mathbf{C}\backslash \{0\})^{n}}\] решение системы (1). Найдем решение идемпотентной системы уравнений (6), (7), ближайшее к вектору $\mathrm{Log}((E^{*},\overline{E}^{*}))$ в
$\mathrm{R}^{2n}$. Будем говорить, что область параметров деформации является областью притяжения (отталкивания), если
при соответствующей деформации решение идемпотентной системы приближается (удаляется) от подпространства, 
определяемого уравнениями \eqref{rho11} в 
$\mathrm{R}^{2n}$. 
\\Задача нахождения границы, разделяющих области 
притяжения и отталкивания в случае идемпотентой системы может рассматриваться, 
как естественный модельный аналог задачи нахождения запаса устойчивости. 
Такая задача, по сравнению с анализом зависимости от параметров решений 
вещественных многомерных систем алгебраических уравнений, 
решается существенно проще в идемпотентном анализе, где она, по существу,
сводится к анализу систем линейных уравнений зависящих от параметров.  

\renewcommand{\theequation}{\thesection.\arabic{equation}}
 
\newpage
\renewcommand{\theopred}{\arabic{opred}}
\renewcommand{\theteorema}{\arabic{teorema}}
\renewcommand{\theequation}{\arabic{equation}}
\setcounter{equation}{0}
\setcounter{teorema}{0}
\setcounter{footnote}{0}
\nachalo{О.В. Знаменская}{Классические и неархимедовы амебы в вопросах
расширения полей}
\label{zna-rus-abs}
Поле $\mathbb P$ называется алгебраическим расширением, 
или  расширением  Галуа поля $\mathbb K$, если существует 
алгебраическое уравнение 
\begin{equation}
	c_0+c_1x+c_2x^2+\dots+c_nx^n=0
\end{equation}
с коэффициентами в поле $\mathbb K$, такое, что поле $\mathbb P$ получается присоединением к $\mathbb K$ всех корней этого уравнения. Известно, что все такие расширения конечномерны.

\markboth{Oksana V.~Znamenskaya}{Classical and nonarchimedean amoebas (in Russian)}
Таким образом, в классической теории конечные расширения строятся при помощи присоединения к исходному полю нулей полиномов от одного переменного. Наша цель --- изучение бесконечных аналогов этих расширений, определяемых полиномами от нескольких переменных, для случая неархимедовых полей.
Более точно, наша задача состоит в описании бесконечных расширений подполей  неархимедова поля  $\mathbb K$ рядов Пюиз{о}.

Напомним, что в поле $\mathbb K$ с неархимедовым нормированием норма элемента $a\in \mathbb K$ может быть определена   через  показатель нормирования $\mbox{\rm val}(a)$ поля $\mathbb K$ при помощи соотношения $|a|=e^{-\mbox{\rm val}(a)}$.
Здесь $\mbox{\rm val}(a)$ есть отображение $\mathbb K\rightarrow \mathbb R\cup\{\infty\}$, определенное на элементах $\mathbb K$ и удовлетворяющее следующим условиям \cite{BoSh}:
\begin{enumerate}
	\item[a)] $\mbox{\rm val}(a)=\infty$ тогда и только тогда, когда $a=0$; 
	\item[b)] $\mbox{\rm val}(ab)=\mbox{\rm val}(a)+\mbox{\rm val}(b)$;
	\item[c)] $\mbox{\rm val}(a+b)\geqslant\min(\mbox{\rm val}(a),\mbox{\rm val}(b))$.
\end{enumerate}

Пусть $\mathbb K$ --- поле рядов Пюиз{о} с коэффициентами в произвольном поле $k$, т.е. рядов $a(t)$ вида 
$$
a(t)=\sum\limits_{q_j\in A_a} \xi_j {\displaystyle t}^{q_j}
$$
по дробным степеням $q_j$ переменного  $t$, где $A_a\subset \mathbb Q$ --- вполне упорядоченное множество.
Показатель неархимедова нормирования $\mbox{\rm val}$ в этом случае полагается  равным
$\min A_a$.

Бесконечные расширения поля $\mathbb K$ будем строить следующим образом.
Рассмотрим полином 
 \begin{equation}
 \label{znam1}
	f=\sum\limits_{\alpha\in A} a_{\alpha}(t) z^{\alpha}  
\end{equation}
из $\mathbb K[z_1,\dots,z_n]$ с коэффициентами 
$a_{\alpha}(t)\in L\subset \mathbb K$, где $L$ --- подполе $\mathbb K$.
Определим расширение $\mathbb P/L$ как множество всевозможных значений полиномов 
 \begin{equation*}
\sum\limits_{\beta\in B} b_{\beta} {z(t)}^{\beta},
\end{equation*}
где $z(t)=\big(z_1(t),\dots,z_n(t)\big)$ является решением уравнения $f=0$ для полинома $f$ вида \eqref{znam1}.

Отметим, что самый простой случай расширений Галуа --- циклические расширения, получаются присоединением к исходному полю всех корней из единицы, т.е. решений двучленного уравнения $x^m-a=0$.
 Очевидно, все корни полинома Галуа $f(x)=x^m-a$:
\begin{enumerate}
	\item[{\small$\bullet$}]  лежат на окружности;
	\item[{\small$\bullet$}] на окружности они равномерно распределены.
\end{enumerate}
 Многомерный  аналог первой из указанных геометрических характеристик решений полиномов Галуа может быть сформулирован на языке амеб.

Амебой $\mathcal A_f$ комплексной
гиперповерхности $V\in (\mathbb\setminus \{0\})^n$, задаваемой
полиномом $f$ (см. \cite{GKZ-Znam}), называется ее образ при отображении   
 $$
 \mbox{\rm Log }: (\mathbb C\setminus \{0\})^n\longrightarrow \mathbb R^n,
 $$
действующем по правилу   
$$
 (z_1,\dots,z_n)\longrightarrow (\log|z_1|,\dots,\log|z_n|).
 $$
Амебу, определяемую таким образом, назовем классической.

\begin{opred}[см. \cite{PasTh-Znam}]
Классическая амеба $\mathcal A_f$ называется
солидной, если число связных компонент дополнения к ней минимально. 
\end{opred}

Многомерное обобщение того факта,
что все корни $f(x)=x^m-a$ лежат на единичной окружности
на языке амеб выражается в том, что \textit{амеба $\mathcal A_f$ полинома $f$ солидна}. 

Заметим, что если $n=1$, то амеба произвольного
полинома от одного переменного есть конечное
 множество точек в $\mathbb R^1$. Солидными в этом случае
 будут только амебы, состоящие из одной точки и имеющие
 лишь две связные компоненты в дополнении, а  это и есть
 в точности амебы полиномов Галуа.  

Напомним, что многогранником Ньютона  $\mathcal N_f$ полинома  $f$ от $n$  переменных называется выпуклая оболочка показателей
его мономов в  $\mathbb R^n$. 
Конусом рецессии выпуклого множества $E\subset \mathbb R^n$  называется максимальный конус среди тех, 
которые сдвигом можно поместить в  $E$.

Согласно результатам М. Форсберга, М. Пассаре и А.К. Циха \cite{FPT-Znam}, 
справедлива
\begin{teorema} \label{teor}
Существует естественная инъективная функция порядка $\nu$  на 
множестве $\{E\}$  связных компонент дополнения 
$\mathbb R^n\setminus \mathcal A_f$   амебы гиперповерхности  $f=0$, 
сопоставляющая каждой компоненте  $E$ некоторую 
целочисленную точку $\nu(E)$  из многогранника 
Ньютона  $\mathcal N_f$. Конус рецессии
компоненты $E$  совпадает с конусом,
двойственным к  $\mathcal N_f$ в точке  $\nu(E)$.
\end{teorema}
Таким образом, многогранник Ньютона отражает структуру классической амебы. 
В частности, число связных компонент дополнения 
$\mathbb R^n\setminus \mathcal A_a$ не меньше числа вершин и не 
больше числа всех целых точек многогранника Ньютона  $\mathcal N_f$. 
Ясно, что классическая амеба солидна, если число 
компонент дополнения строго равно числу вершин $\mathcal N_f$.

С точки зрения многогранников Ньютона, многомерным аналогом полиномов Галуа являются, так называемые, \textit{максимально разреженные полиномы} \cite{PasTh-Znam}, т.е. 
полиномы вида:
\begin{equation*}
	f=\sum\limits_{\alpha\in A} a_{\alpha} z^{\alpha},
\end{equation*}
где с ненулевыми коэффициентами входят только мономы, соответствующие вершинам $\mathcal N_f$.
М.~Ниссе был заявлен результат, что классическая амеба любого максимально разреженного полинома солидна.

Далее нас будет интересовать вопрос солидности неархимедовых амеб нулевого множества максимально разреженных полиномов \eqref{znam1}, при помощи которых строятся бесконечные расширения $\mathbb P/L$ неархимедова поля рядов Пюиз{о}.

Определим  по аналогии свойство солидности для неархимедовых амеб.
Пусть $\mathbb K$ --- произвольное неархимедово поле и $\mbox{\rm val}(a)$ --- его показатель нормирования.
\begin{opred}[см. \cite{K-Znam, EKL}]
Амебой $\mathcal A(V)$ алгебраической гиперповерхности $V\subset{(\mathbb K^*)}^n$ называется замыкание образа $V$ при отображении 
$$
\mbox{\rm Log}: (z_1,\dots, z_n)\rightarrow (-\mbox{\rm val}(z_1),\dots, -\mbox{\rm val}(z_n)).
$$
\end{opred}

Из теоремы \ref{teor} следует, что конусы рецессии  всех компонент дополнения классической солидной амебы полномерны. Соответственно, для неархимедова случая можно дать следующее 
\begin{opred}Неархимедова амеба $\mathcal A(V)$ называется солидной, если любая связная компонента дополнения к ней имеет полномерный конус рецессии.
\end{opred}

Пусть $\mathbb K$ --- поле рядов Пюиз{о}.
\begin{teorema}
Неархимедова амеба максимально разреженного полинома, определяющего многомерное расширение Галуа поля $\mathbb K$, солидна.
\end{teorema}

Согласно \cite{EKL} существует двойственность между неархимедовой амебой 
$\mathcal A(V)$ и подразбиением многогранника Ньютона полинома, определяющего $V$. С учетом этого справедлива
\begin{teorema}
В случае $n=2$ если неархимедова амеба солидна, то, как граф, она не имеет циклов.
\end{teorema}

 Пусть  $f(z)=\sum\limits_{\alpha\in A} c_\alpha z^\alpha$, $z=(z_1,\dots,z_n)\in \mathbb C^n$ --- полином, все  коэффициенты которого  имеют рациональные модули, т.е. все $|c_\alpha|\in\mathbb Q$. Определим полином $F$ следующим образом: 
 $$
 F=\sum\limits_{\alpha\in A} c_\alpha(t) \xi^\alpha,
 $$
 где $c_\alpha(t)$ таковы, что $|c_\alpha(t)|=e^{-|c_\alpha|}$. 
 В указанных предположениях справедлива 
 
  \begin{teorema}
 Если классическая амеба $\mathcal A_f$ солидна, то и неархимедова амеба $\mathcal A(V)$ солидна.
 \end{teorema}

В доказательстве теоремы используются понятия хребта амебы и тропического многообразия, определяемого при помощи тропикализации полинома $F$.
\markboth{G.V.~Koval, V.P.~Maslov}{Classical and nonarchimedean amoebas (in Russian)} 

\renewcommand{\theopred}{\thesection.\arabic{opred}}
\renewcommand{\theteorema}{\thesection.\arabic{teorema}}
\renewcommand{\theequation}{\thesection.\arabic{equation}}
\renewcommand{\theequation}{\arabic{equation}}
\renewcommand{\thepredl}{\arabic{predl}}
\setcounter{predl}{0}
\setcounter{equation}{0}
\setcounter{footnote}{0}
\nachaloe{Г.В.~Коваль и В.П.~Маслов}{Обобщение ультравторичного квантования
для фермионов при ненулевой температуре}%
{Работа выполнена при поддержке грантов РФФИ  05-01-00824 и 
05-01-02807-НЦНИЛ\_а.}
%\markboth{Г.В.~Коваль и В.П.~Маслов}{Обобщение ультравторичного
%квантования}

\label{kov-rus-abs}

В работах В.П. Маслова развит метод ультравторичного квантования и
концепция истинного символа\markboth{G.V.~Koval, V.P.~Maslov}{Generalization of ultrasecond
quantization (in Russian)}
\cite{Ultravtor_kvant,Maslov_TMF_2002,Maslov_RJMP_2002}. Этот
метод позволяет находить асимптотические серии систем большого
числа частиц при нулевой температуре. В частности, некоторые серии
определяются периодическими решениями системы уравнений
Гамильтона, соответствующей истинному
символу~\cite{Maslov_TMF_2002,Maslov_RJMP_2002} рассматриваемой
физической системы. В данной работе найдено соответствие между
уравнениями метода ультравторичного квантования по парам для
фермионов и уравнениями вариационного метода Боголюбова. Так как
вариационный метод Боголюбова применим для случая ненулевой
температуры, из принципа соответствия получено обобщение уравнений
метода ультравторичного квантования фермионов на температурный
случай.

В статье~\cite{Maslov_RJMP_2002} показано, что асимптотика серий
собственных значений системы $N$ тождественных фермионов в пределе
при $N\to\infty$ определяется решениями следующей системы
уравнений
\begin{equation}
\begin{split}
\Omega\Phi(x,y)&= \left(-\frac{\hbar^2}{2m}
\left(\Delta_x+\Delta_y\right)+U(x)+U(y)+V(x,y)\right) \Phi(x,y)+\\
&+2\iint dzdw\ \left(V(x,y)+V(z,w)\right)
\Phi^+(z,w)\Phi(x,z)\Phi(w,y),\\[1.5ex]
\Omega\Phi^+(x,y)&= \left(-\frac{\hbar^2}{2m}
\left(\Delta_x+\Delta_y\right)+U(x)+U(y)+V(x,y)\right) \Phi^+(x,y)
+\\
&+2\iint dzdw\ \left(V(x,z)+V(y,w)\right)
\Phi^+(x,z)\Phi^+(w,y)\Phi(z,w), \label{koval1}
\end{split}
\end{equation}
где $x,y\in\cM$ --- координаты частиц, пространство $\cM$
определяется задачей, например, это может быть $\bR^3$, или
трехмерный тор, $\Delta_x$, $\Delta_y$ --- операторы Лапласа,
действующий по соответствующей переменной $x$ или $y$, $U(x)$
--- потенциал внешнего поля, $V(x,y)$ --- потенциал
взаимодействия, симметричный относительно перестановки переменных
$x$ и $y$, $m$
--- масса частиц, $\hbar$ --- постоянная Планка, $\Omega$ ---
действительное число. Функции $\Phi(x,y),\Phi^+(x,y)\in L_2(\cM)$
--- антисимметричны относительно перестановок переменных $x$ и $y$
и удовлетворяют условию
\begin{equation}
\iint dxdy\ \Phi^+(x,y)\Phi(x,y) = \frac{N}2. \label{koval2}
\end{equation}

Уравнения~(\ref{koval1}) получаются при ультравторичном квантовании по
парам. В таком квантовании рассматриваемой системе фермионов
отвечает истинный символ вида
\begin{equation}
\begin{split}
&A[\Phi^+,\Phi]=\\
&\int\!\int dxdy \Phi^+(x,y) \left(
-\frac{\hbar^2}{2m} \left( \Delta_x + \Delta_y \right) + U(x) +
U(y) + V(x,y) \right)
\Phi(x,y) +\\
&+2\int\!\int dxdydzdw V(x,y)\Phi^+(x,y)\Phi^+(z,w)
\Phi(x,z)\Phi(w,y), \label{koval000}
\end{split}
\end{equation}
который является функционалом от двух 
антисимметричных функций
$\Phi^+(x,y)$ и $\Phi(x,y)$ из $L_2(\cM^2)$. 
Этому символу соответствует
система уравнений Гамильтона
{\tolerance=500

}
\begin{equation}
i\frac{\pa\Phi}{\pa t}(x,y,t) = \frac{\delta A}
{\delta\Phi^+(x,y,t)}, \qquad -i\frac{\pa\Phi^+}{\pa t}(x,y,t) =
\frac{\delta A}{\delta\Phi(x,y,t)}, \label{koval001}
\end{equation}
где в правой части уравнений стоит вариационная производная
функционала~(\ref{koval000}). Система уравнений~(\ref{koval001}) имеет
интеграл движения, вид которого совпадает с выражением в левой
части равенства~(\ref{koval2}). Уравнения~(\ref{koval1}) получаются
из~(\ref{koval001}) в частном случае, когда
\begin{equation}
\Phi(x,y,t)=\Phi(x,y) e^{-i\Omega t},\qquad
\Phi^+(x,y,t)=\Phi^+(x,y) e^{i\Omega t}. \label{koval015}
\end{equation}

Запишем~(\ref{koval1}) в другом виде. Введем функции $G(x,y)$, $R(x,y)$
и $\wt{R}(x,y)$
\begin{equation}
\begin{split}
&\wt{R}(x,y)=\Phi^+(x,y),\qquad G(x,y)=2\int dz\
\Phi^+(x,z)\Phi(y,z),\\
&R(x,y)=2\left(\Phi(x,y)-\int dz\ \Phi(x,z)G(z,y)\right).
\label{koval3}
\end{split}
\end{equation}
В силу антисимметрии функций $\Phi^+(x,y)$, $\Phi(x,y)$ для
функций~(\ref{koval3}) выполняются равенства
\begin{align}
&R(x,y)=-R(y,x),\qquad \wt{R}(x,y)=-\wt{R}(y,x),\label{koval4} \\
&G(x,y)=\int dz\ G(x,z)G(z,y) + \int dz\ \wt{R}(z,x) R(z,y). \label{koval5}
\end{align} 
Из~(\ref{koval2}) следует
\begin{equation}
\int dx\ G(x,x) = N. \label{koval6}
\end{equation}

%\begin{theorem}
\begin{predl}
Функции $G(x,y)$, $\wt{R}(x,y)$, $R(x,y)$
удовлетворяют системе уравнений
\begin{equation}
\begin{split}
&\left( -\frac{\hbar^2}{2m} \left(\Delta_x-\Delta_y\right) + U(x)- U(y) \right) G(x,y) -\\
&- \int dz\ \left(V(x,z)-V(y,z)\right)\wt{R}(x,z) R(z,y)=0,\\[1.5ex]
&\left( -\frac{\hbar^2}{2m} \left(\Delta_x+\Delta_y\right) + U(x)+ U(y) +V(x,y) \right) R(x,y) - \\
&- \int dz\ \left(V(x,z)G(z,y)R(x,z) +V(y,z)G(z,x)R(z,y)\right) = \Omega R(x,y), \\[1.5ex]
&\left( -\frac{\hbar^2}{2m} \left(\Delta_x+\Delta_y\right) + U(x)+ U(y) +V(x,y) \right) \wt{R}(x,y) - \\
&- \int dz\ \left(V(x,z)G(y,z)\wt{R}(x,z) +V(y,z)G(x,z)\wt{R}(z,y)\right) = \Omega \wt{R}(x,y).\label{koval7}
\end{split}
\end{equation}
\end{predl}
%\end{theorem}

Действительно, непосредственной проверкой удостоверяется, что
функции~(\ref{koval3}) удовлетворяют уравнениям~(\ref{koval7}), если функции
$\Phi^+(x,y)$, $\Phi(x,y)$ удовлетворяют~(\ref{koval1}).

Чтобы обобщить уравнения~(\ref{koval1}) на случай ненулевой
температуры, применим принцип соответствия между этими уравнениями
и уравнениями вариационного принципа Боголюбова при ненулевой
температуре. То есть исходя из температурных уравнений
вариационного принципа Боголюбова, по принципу соответствия найдем
обобщение уравнений~(\ref{koval1}) на температурный случай.

Рассмотрим уравнения вариационного метода Боголюбова. При
температуре $\theta\ge0$ для рассматриваемой системы фермионов из
вариационного метода Боголюбова~\cite{Bogol-Koval} получаются уравнения
%\newpage
\begin{equation}
\begin{split}
\la_\alpha u_\alpha(x)&= \left(-\frac{\hbar^2}{2m}\Delta +
U(x)-\mu\right) u_\alpha(x) + \int dy\ V(x,y) R_B(x,y)v_\alpha^*(y) +\\
&+\int dy\ V(x,y)\left(G_B(y,y)u_\alpha(x)-G_B(y,x)u_\alpha(y)
\right),\\[1.5ex]
-\la_\alpha v_\alpha(x)&= \left(-\frac{\hbar^2}{2m}\Delta + U(x)
-\mu\right) v_\alpha(x)+\int dy\ V(x,y) R_B(x,y)u_\alpha^*(y)+\\
&+\int dy\ V(x,y)\left(G_B(y,y)v_\alpha(x)-G_B(y,x)v_\alpha(y)
\right), \label{koval8}
\end{split}
\end{equation}
где $\alpha=1,2,\dots$, функции $u_\alpha(x)$ и $u_\alpha^*(x)$, а
также $v_\alpha(x)$ и $v^*_\alpha(x)$, комплексно сопряжены друг
другу и удовлетворяют условиям
\begin{equation}
\begin{split}
\int dx\ \left(u^*_\alpha(x) v_\beta(x) + v_\alpha(x)
u^*_\beta(x)\right)&=\int dx\ \left(u_\alpha(x) v^*_\beta(x) +
v^*_\alpha(x) u_\beta(x)\right) = 0,  \\
\int dx\ \left(u^*_\alpha(x) u_\beta(x) + v_\alpha(x)
v^*_\beta(x)\right)&=\delta_{\alpha\beta},\ 
\forall\alpha,\beta=1,2,\dots, \label{koval9}
\end{split}
\end{equation}
где $\delta_{\alpha\beta}$ --- символ Кронекера. Кроме того в
уравнениях~(\ref{koval9}) функции $R_B(x,y)$ и $G_B(x,y)$ имеют вид
\begin{equation}
\begin{split}
R_B(x,y)&=\sum_{\alpha=1}^\infty\left( \frac12 - n_\alpha
\right)\left(v_\alpha(x)
u_\alpha(y)- v_\alpha(y)u_\alpha(x)\right) ,\\
G_B(x,y)&=\sum_{\alpha=1}^\infty \left( v^*_\alpha(x)
v_\alpha(y) \left(1-n_\alpha\right) + u^*_\alpha(x) u_\alpha(y)
n_\alpha \right), \label{koval10}
\end{split}
\end{equation}
где
\begin{equation}
n_\alpha= \frac1{\exp(\la_\alpha/\theta)+1}, \label{koval11}
\end{equation}
а $\mu$ определяется из условия, что функция $G_B(x,y)$~(\ref{koval10})
удовлетворяет равенству
\begin{equation}
\int dx\ G_B(x,x)=N. \label{koval12}
\end{equation}

Функции~(\ref{koval10}) и комплексно сопряженная к 
$R_B(x,y)$ функция
$R_B^*(x,y)$ при любой температуре $\theta$ 
удовлетворяют
равенствам
{\sloppy

}
\begin{eqnarray}
&&R_B(x,y)=-R_B(y,x),\qquad R_B^*(x,y)=-R_B^*(y,x),\label{koval13}\\
&&G_B(x,y)=G_B^*(y,x), \label{koval14}
\end{eqnarray}
а из уравнений~(\ref{koval8}) следует, что эти функции также
удовлетворяют системе уравнений
\begin{equation}
\begin{split}
&\left(-\frac{\hbar^2}{2m} \left(\Delta_x-\Delta_y\right) + U(x)
- U(y) \right) G_B(x,y) +\\
& + \int dz\left(V(x,z)-V(y,z)\right)
R_B^*(x,z) R_B(z,y)+\\
& +\int dz\left(V(x,z)-V(y,z)\right)\left( G_B(z,z)G_B(x,y) -
G_B(x,z)G_B(z,y) \right)=0,\\[1.5ex]
&\left(-\frac{\hbar^2}{2m} \left(\Delta_x+\Delta_y\right) + U(x)
+ U(y) +V(x,y) \right) R_B(x,y) - \\
& - \int dz\left(V(x,z)G_B(z,y)R_B(x,z) +
V(y,z)G_B(z,x)R_B(z,y)\right)+\\
&+\int dz V(x,z)\left(G_B(z,z)R_B(x,y)-G_B(z,x)R_B(z,y)\right)+\\
&+\int dz V(y,z)\left(G_B(z,z)R_B(x,y)-G_B(z,y)R_B(x,z)\right)
=2\mu R_B(x,y),\label{koval15}
\end{split}
\end{equation}
где дополнительное уравнение получается комплексным сопряжением
второго уравнения формулы ~(\ref{koval15}).

Если $\theta=0$, то из~(\ref{koval11}) следует, что $n_\alpha$
принимает значение $0$ или $1$ для всех $\alpha=1,2,\dots$. Тогда
из~(\ref{koval9}) следует, что функции~(\ref{koval10}) при нулевой
температуре удовлетворяют условию
\begin{equation}
G_B(x,y)=\int dz\ G_B(x,z)G_B(z,y) + \int dz\ R_B^*(z,x) R_B(z,y).
\label{koval16}
\end{equation}
Равенства~(\ref{koval16}), (\ref{koval13}) и (\ref{koval12}) совпадают
соответственно с равенствами~(\ref{koval5}), (\ref{koval4}) и (\ref{koval6}),
если по принципу соответствия заменить $R_B^*(x,y)$ на
$\wt{R}(x,y)$, $R_B(x,y)$ на $R(x,y)$, $G_B(x,y)$ на $G(x,y)$.
Уравнения~(\ref{koval15}) при такой замене не переходят в
уравнения~(\ref{koval7}), однако, если $2\mu$ в~(\ref{koval15}) заменить на
$\Omega$, то очевидно соответствие между одной системой и другой,
соответствующие друг другу уравнения отличаются несколькими
слагаемыми в левой части. Кроме того, для частного вида
взаимодействия, например, такого как в модели БКШ~\cite{Shriff},
уравнения~(\ref{koval15}) после замены совпадают с~(\ref{koval7}).
{\tolerance=1000

}
Отметим, что функции $G_B(x,y)$, $R_B(x,y)$, $R^*_B(x,y)$
удовлетворяют большему числу условий, чем функции $G(x,y)$,
$R(x,y)$, $\wt{R}(x,y)$. Функции $\Phi(x,y)$, $\Phi^+(x,y)$,
удовлетворяющие уравнениям~(\ref{koval1}), в общем случае не являются
комплексно сопряженными друг другу~\cite{Maslov_TMF_2002}. Поэтому
из формул~(\ref{koval3}) следует, что функция $\wt{R}(x,y)$ не должна
быть комплексно сопряженной к $R(x,y)$, а функция $G(x,y)$ не
должна удовлетворять условию~(\ref{koval14}).

Из соответствия между функциями~(\ref{koval3}) и~(\ref{koval10}), получим,
что температурным аналогом системы уравнений~(\ref{koval1}) является
следующая система уравнений:
\begin{equation}
\begin{split}
\la_\alpha u_\alpha(x) &= \left(-\frac{\hbar^2}{2m}\Delta +
U(x)-\mu\right) u_\alpha(x) + \int dy V(x,y) R(x,y)\wt{v}_\alpha(y),\\
-\la_\alpha v_\alpha(x) &= \left(-\frac{\hbar^2}{2m}\Delta + U(x)
-\mu\right) v_\alpha(x)+\int dy V(x,y) R(x,y)\wt{u}_\alpha(y),\\
\la_\alpha \wt{u}_\alpha(x) &= \left(-\frac{\hbar^2}{2m}\Delta +
U(x)-\mu\right) \wt{u}_\alpha(x)+\int dy V(x,y) \wt{R}(x,y)v_\alpha(y),\\
-\la_\alpha \wt{v}_\alpha(x) &= \left(-\frac{\hbar^2}{2m}\Delta +
U(x)-\mu\right) \wt{v}_\alpha(x)+\int dy V(x,y) \wt{R}(x,y)
u_\alpha(y), \label{17}
\end{split}
\end{equation}
где $G(x,y)$, $R(x,y)$ и $\wt{R}(x,y)$ выражаются следующим
образом:
\begin{equation}
\begin{split}
G(x,y)&=\sum_{\alpha=1}^\infty \left( \wt{v}_\alpha(x)
v_\alpha(y) \left(1-n_\alpha\right) + \wt{u}_\alpha(x) u_\alpha(y)
n_\alpha\right),\\
R(x,y)&=\sum_{\alpha=1}^\infty\left( \frac12 - n_\alpha \right)
\left(v_\alpha(x) u_\alpha(y) - v_\alpha(y)u_\alpha(x)\right) ,\\
\wt{R}(x,y)&=\sum_{\alpha=1}^\infty\left( \frac12 - n_\alpha
\right)\left( \wt{v}_\alpha(x) \wt{u}_\alpha(y)- \wt{v}_\alpha(y)
\wt{u}_\alpha(x)\right), \label{18}
\end{split}
\end{equation}
$n_\alpha$ выражается через $\lambda_\alpha$ и $\theta$
формулой~(\ref{koval11}), а функции $u_\alpha(x)$, $v_\alpha(x)$,
$\wt{u}_\alpha(x)$, $\wt{v}_\alpha(x)$, $\alpha=1,2,\dots$ кроме
уравнений~(\ref{17}) еще удовлетворяют условиям
\begin{equation}
\begin{split}
&\int dx\left(\wt{u}_\alpha(x) v_\beta(x) + v_\alpha(x)
\wt{u}_\beta(x)\right)=\int dx\left(u_\alpha(x) \wt{v}_\beta(x) +
\wt{v}_\alpha(x) u_\beta(x)\right) = 0,\\
&\int dx\left(\wt{u}_\alpha(x) u_\beta(x) + v_\alpha(x)
\wt{v}_\beta(x)\right)=\delta_{\alpha\beta},\quad
\forall\alpha,\beta=1,2,\dots. \label{19}
\end{split}
\end{equation}
Параметр $\mu$ в уравнениях~(\ref{17}) определяется из условия,
что функция $G(x,y)$ из (\ref{18}) удовлетворяет условию~(\ref{koval6}).
%\begin{theorem}
\begin{predl}
Если функции $u_\alpha(x)$,
$\wt{u}_\alpha(x)$, $v_\alpha(x)$, $\wt{v}_\alpha(x)$,
$\alpha=1,2,\dots$ удовлетворяют системе уравнений~(\ref{17}) и
условиям~(\ref{19}), то функции~(\ref{18}) при 
любом $\theta\ge0$
удовлетворяют системе уравнений~(\ref{koval7}) 
с $\Omega=2\mu$ и
условиям~(\ref{koval4},\ref{koval6}), 
а при $\theta=0$ еще удовлетворяют
условию~(\ref{koval5}).
{\tolerance=500

}
\end{predl}

%\end{theorem}

%{\bf Доказательство.}
В силу уравнений~(\ref{17}) и условий~(\ref{19}) это утверждение
доказывается прямой подстановкой функций~(\ref{18}) в
формулы~(\ref{koval4}-\ref{koval7}).

Множество решений уравнений~(\ref{koval15}) шире, чем множество
функций~(\ref{koval10}), выраженных через решения системы
уравнений~(\ref{koval8}). В~\cite{Maslov_TMF_2002} показано, что
система уравнений~(\ref{koval15}) может быть записана в виде пары:
\begin{equation}
\left[\wh{A},\wh{L}\right]=0, \label{20}
\end{equation}
а множеству решений температурных уравнений~(\ref{koval8})
соответствует такое решение уравнения~(\ref{20}), для которого
\begin{equation}
\wh{A}=f_\theta(\wh{L}),\label{21}
\end{equation}
где
\begin{equation}
f_\theta(\xi)=\frac1{\exp(\xi/\theta)+1}-\frac12. \label{22}
\end{equation}
Для уравнений~(\ref{17}) и~(\ref{17}) справедливо аналогичное
утверждение. Рассмотрим матрицы $\wh{A}$ и $\wh{L}$ вида
\begin{equation}
\wh{A}=\left(
\begin{array}{cc}
\wh{G}-\frac12&-\wh{\wt{R}}\\
\wh{R}&\frac12-\wh{G}^t\\
\end{array}\right),\qquad
\wh{L}=\left(
\begin{array}{cc}
\wh{T}&-\wh{\wt{B}}\\
\wh{B}&-\wh{T}\\
\end{array}\right),\label{koval23}
\end{equation}
где $\wh{G}$, $\wh{R}$, $\wh{\wt{R}}$ --- операторы в пространстве
$L_2(\cM)$, задаваемые интегральными ядрами $G(x,y)$, $R(x,y)$,
$\wt{R}(x,y)$ соответственно, $\wh{G}^t$ --- оператор, задаваемый
в $L_2(\cM)$ ядром $G^t(x,y)=G(y,x)$, $\wh{B}$ и $\wh{\wt{B}}$
--- ядрами $B(x,y)=V(x,y)R(x,y)$ и $\wt{B}(x,y)=V(x,y)\wt{R}(x,y)$
соответственно, а оператор $\wh{T}$ --- оператор Гамильтона для
одной частицы, то есть оператор вида
$$
\wh{T}=-\frac{\hbar^2}{2m}\Delta+U(x).
$$

Подстановка~(\ref{koval23}) в~(\ref{20}) приводит к четырем уравнениям,
из которых два совпадают, а три независимых приводятся к
виду~(\ref{koval7}). Поэтому справедливо следующее утверждение.
\markboth{Alexey G.~Kushner}{Generalization of ultrasecond 
quantization (in Russian)}
%\begin{theorem}
\begin{predl}
Система уравнений~(\ref{7}) может быть
записана в виде~(\ref{20}), где $\wh{A}$ и $\wh{L}$ имеют
вид~(\ref{koval23}), а $2\mu=\Omega$.
\end{predl}
%\end{theorem}

%{\bf Доказательство.}

Кроме того, решениям системы уравнений~(\ref{koval7}), вида~(\ref{18}),
которые получены из решений уравнений~(\ref{17},\ref{19}),
соответствуют такие $\wh{A}$ и $\wh{L}$, что для них справедливо
равенство~(\ref{21}). Это является следствием
уравнений~(\ref{17},\ref{19}).

В заключение отметим, что уравнения~(\ref{koval7}) в температурном
случае, полученные здесь из принципа соответствия, могут быть
строго получены из истинного символа
\cite{Maslov_TMF_2002,Maslov_RJMP_2002} для ультравторично
квантованного уравнения, отвечающего матрице плотности
\cite{Ultravtor_kvant}.

% Это позволяет строго перенести результаты \cite{Ultravtor_kvant,
%Maslov_RJMP_2005} на температурный случай.

\renewcommand{\theequation}{\thesection.\arabic{equation}}
\renewcommand{\thepredl}{\thesection.\arabic{predl}}
%\newpage
\setcounter{section}{0}
\setcounter{footnote}{0}
\nachalo{А.Г.~Кушнер}{Контактная классификация уравнений Монжа-Ампера}
%\markboth{В.П.~Маслов}{Ультравторичное квантование} 
\markboth{Alexey G.~Kushner}{Contact classification of Monge-Amp\`ere equations (in Russian)}
\label{kus-rus-abs}

\section{Геометрические структуры, ассоциированные с уравнениями Монжа-Ампера}

Класс уравнений Монжа-Ампера выделяется из многообразия уравнений
второго порядка тем, что он замкнут относительно контактных
преобразований. Это обстоятельство было известно еще Софусу Ли,
изучавшему уравнения Монжа-Ампера методами созданной им контактной
геометрии. В 1870-х и 1880-х он поставил проблемы классификации
уравнений Монжа-Ампера относительно (псевдо)группы контактных
преобразований, в частности, о приведении уравнений Мон\-жа-Ам\-пе\-ра к
квазилинейной форме и наиболее простом координатном представлении
таких уравнений \cite{Lie1872}.

В 1979 г. в работе \cite{Lch1979} В.\,В.~Лычагин показал,что
уравнения Монжа-Ампера  допускают эффективное описание в терминах
дифференциальных форм на многообразии $1$-джетов гладких функций.
Отправной точкой является следующее наблюдение.

Пусть $M$ --- $n$-мерное гладкое многообразие, $J^1M$ ---
многообразие $1$-джетов гладких функций на $M$. На $J^1M$
естественным образом определена контактная структура ---
распределение Картана $C$, в канонических локальных координатах
Дарбу $(q,u,p)=(q_1,\dots ,q_n,u,p_1,\dots , p_n)$ задаваемое
дифференциальной $1$-формой Картана $U=du-pdq$. Ограничение
дифференциала формы Картана на подпространство Картана не вырождено
на нем и определяет симплектическую структуру $ \Omega_{a}=
dU|_{C(a)}\in\Lambda^{2}\left( C^{\ast}(a)\right)$.

Со всякой дифференциальной $n$-формой $\omega\in\Omega^n(J^1M)$
свяжем нелинейный дифференциальный оператор
$\Delta_{\omega}:C^{\infty}(M)\rightarrow\Omega^n(M),$ действующий
на гладкую функцию $v$ следующим образом:
\begin{equation}
\Delta_{\omega}(v)=j_{1}(v)^{\ast}(\omega). \label{Intro_Delta}
\end{equation}
Здесь $j_{1}(v):M\rightarrow J^{1}M$ --- $1$-джет функции $v\in
C^{\infty}(M)$.

Оператры $\Delta_{\omega}$ называются операторами
\emph{Монжа-Ампера}, а уравнение $ E_{\omega}=\left\{ \Delta_{\omega
}(v)=0\right\}  \subset J^{2}M$
--- {\it уравнением Монжа-Ампера}. Следующее обстоятельство
оправдывает эти названия: будучи записанным в локальных канонических
координатах на $J^{1}M$, оператор $\Delta_{\omega}$ имеет тот же
самый тип нелинейности по производным второго порядка, что и
классические опреаторы Монжа-Ампера, а именно, нелинейности типа
определителя матрицы Гессе и ее миноров. При $n=2$ мы получаем
классическое уравнение Монжа-Ампера:
\begin{equation}
\label{MAE}
    Av_{xx}+2Bv_{xy}+Cv_{yy}+D(v_{xx}v_{yy}-v_{xy}^{2})+E=0,
\end{equation}
\noindent где $A,B,C,D,E$ --- функции от независимых переменных
$x,y$, функции $v=v(x,y)$ и ее первых производных $v_{x},v_{y}$.

Преимуществом такого подхода перед классическим является редукция
порядка пространства джетов: мы используем более простое
пространство $1$-джетов $J^{1}M$ вместо пространства $2$-джетов
$J^{2}M$, в котором, будучи уравнениями второго порядка, \emph{ad
hoc} должны лежать уравнения Монжа-Ампера.

В случае, когда коэффициенты уравнения (\ref{MAE}) не зависят явно
от функции $v$ ситуация еще более упрощается: в определении
оператора (\ref{Intro_Delta}) вместо пространства $1$-джетов можно
рассматривать кокасательное расслоение $T^\ast M$ многообразия $M$,
а вместо контактной геометрии --- симплектическую. Такие уравнения
Монжа-Ампера будем называть \emph{симплектическими}.

Заметим, что соответствие между дифференциальными $n$-фор\-ма\-ми на
$J^{1}M$ и операторами Монжа-Ампера не является вза\-им\-но-од\-но\-знач\-ным.
Дифференциальные формы, аннулирующимся на любом интегральном
многообразии распределении Картана, образуют идеал $\mathcal{C}$ во
внешней алгебре $\Omega^{\ast}(J^{1}M)$, который называется
\emph{идеалом Картана}. Им отвечает нулевой дифференциальный
оператор. Элементы фак\-тор-ал\-геб\-ры
$\Omega^{\ast}(J^{1}M)/\mathcal{C}$ по этому идеалу называются
\emph{эффективными} формами.

Далее мы будем рассматиривать случай когда $M$ --- двумерное гладкое
многообразие. В терминах эффективных форм можно определить тип
уравнения --- эллиптический, параболический, гиперболический или
переменный. Функция $\operatorname{Pf}\left( \omega\right) \in
C^{\infty}\left( J^{1}M\right)  $, определяемая поточечно равенством
$ \operatorname{Pf}(\omega_{a})
\Omega_{a}\wedge\Omega_{a}=\omega_{a}\wedge\omega_{a}$, называется
{\it пфаффианом} формы $\omega$. Уравнение $E_{\omega}$ называется
{\it гиперболическим, параболическим} или {\it эллиптическим} в
точке $a \in J^{1}M$, если пфаффиан $\operatorname{Pf}(\omega)$
отрицательный, нулевой или положительный в этой точке. Если пфаффиан
не равен нулю в точке, то уравнение называется \emph{невырожденным}.
Очевидным образом понятие типа распространяется на область.

В силу невырожденности симплектической структуры на подпространстве
Картана $C(a)$, формула
$
X_{a}~\rfloor\,\omega_{a}=A_{\omega_{a}}X_{a}\rfloor\,~\Omega_{a},
$
определяет на $C(a)$ ассоциированный с эффективной дифференциальной
$2$-формой $\omega$ линейный оператор $A_{\omega_{a}}$. Здесь
$X_{a}\in C(a)$. Этот оператор симметричен относительно
симплектической структуры, а его квадрат скалярен и
$A_\omega^2+\operatorname*{Pf}(\omega) =0$. Заметим, что операторы
$A_{\omega_{a}}$ не образуют поля эндоморфизмов на $J^1M$, ибо они
определены только на подпространствах Картана.

Если в точке $a\in J^1M$ пфаффиан формы $\omega$ не обращается в
нуль, то в некоторой ее окрестости этой точки форму $\omega$ можно
нормировать так чтобы $\operatorname{Pf}(\omega)=\pm 1$. В этом
случае на подпространсве Картана определена либо структура почти
произведения (для гиперболических уравнений), либо комплексная
структура (для эллиптических уравнений). В первом случае мы получаем
два вещественных, а во втором --- два комплексных $2$-мерных
распределения  на $J^1M$, которые будем называть
\emph{характеристическими} и обозначить через $C_+$ и $C_-$. Эти
распределения косоортогональны друг другу и на каждой из плоскостей
$C_\pm (a)$ $2$-форма $\Omega_a$ не вырождена.

Характеристические распределения порождают еще одно распределение
--- вещественное одномерное распределение $$l=[C_+,C_+]\bigcap
[C_-,C_-],$$ трансверсальное распределению Картана \cite{Lch1993}.

\section{Невырожденные уравнения и инварианты Лапласа}

Пусть $\omega$ --- невырожденная нормированная эффективная
дифференциальная $2$-форма на $J^1M$. В каждой точке $a\in J^1M$
комплексификация касательного пространства $T_a(J^1M)$ распадается в
прямую сумму
\begin{equation}
\label{decomp} T_{a}(J^{1}M)^\mathbb{C}=C_{+}(a)\oplus l(a)\oplus
C_{-}(a).
\end{equation}
Обозначим распределения $C_+,l,C_-$ через $P_1,P_2$ и $P_3$
соответственно. Формула (\ref{decomp}) порождает разложение в прямую
сумму комплекса де Рама многообразия $J^1M$, что позволяет найти
дифференциальные инварианты уравнения \cite{Ksh2006a}.

Определим тензорные поля $q_{j,k}^{s}:D(J^1M)^\mathbb{C}\times
D(J^1M)^\mathbb{C} \rightarrow D(J^1M)^\mathbb{C}$ на $J^1M$:
\begin{equation}
\label{q}
q_{j,k}^{s}(X,Y)=-\mathbf{P}_s\left[\mathbf{P}_jX,\mathbf{P}_kY\right].
\end{equation}
Здесь $\mathbf{P}_j$ --- проектор на распределение $P_j$, $D(J^1M)$
--- модуль векторных полей на $J^1M$.

Мы получаем всего $4$ нетривиальных тензорных поля: $q_{1,2}^3$,
$q_{2,3}^1,$ $q_{1,1}^2,$ $q_{3,3}^2$. Остальные тензоры (\ref{q}) равны
нулю. Определим две дифференциальные $2$-формы как свертки тензоров:
\begin{equation*}
\xi_+=\left\langle q_{1,1}^2,q_{3,2}^1\right\rangle, \quad
\xi_-=\left\langle q_{3,3}^2, q_{1,2}^3\right\rangle.
\end{equation*}
Эти формы мы будем называть \emph{формами Лапласа}, поскольку они
являются обобщением инвариантов Лапласа  для случая линейных
уравнений \cite{Frst1906}. Заметим, что классические инварианты
Лапласа определены только для гиперболических уравнений.

Формы Лапласа играют важную роль при решении вопроса о контактной
линеаризации уравнений Монжа-Ампера. Так, например, если обе формы
Лапласа нулевые, то уравнение Монжа-Ампера локально контактно
эквивалентно либо волновому уравнению $ v_{xx}-v_{yy}=0$, либо
уравнению Лапласа $v_{xx}+v_{yy}=0$ (см. также \cite{Tn1996b}).

В терминах форм Лапласа формулируется решение проблемы
эквивалентности уравнений Монжа-Ампера линейным уравнениям вида
\begin{equation*}
av_{xx}+2bv_{xy}+cv_{yy}+rv_{x}+sv_{y}+kv+w=0,
\end{equation*}
где $a,b,c,r,s,k,w$ --- функции только от независимых переменных
$x,y$ \cite{KLR2007}. В частности, для таких уравнений формы Лапласа
замкнуты.
%\begin{example}
Например, для уравнения Хантора-Сакстона
\begin{equation*}
v_{tx}=vv_{xx}+\kappa u_{x}^{2},
\end{equation*}
возникающего в теории жидких кристаллов, формы Лапласа имеют вид:
$$
\xi_{1}  =-dq_{2}\wedge dp_{1}, \quad  \xi_{2}
=2\left(1-\kappa\right) dq_{2}\wedge dp_{1}.
$$
Это уравнение контактно эквивалентно линейному уравнению
Эйлера-Пуассона \cite{Mrz2004}
\begin{equation*}
v_{tx}=\dfrac{1}{\kappa\left( t+x\right) }v_{t}+\dfrac{2\left(
1-\kappa\right) }{\kappa\left( t+x\right) }v_{x}-\dfrac{2\left(
1-\kappa\right) }{\left( \kappa\left( t+x\right) \right) ^{2}}u.
\end{equation*}
%\end{example}

Если выполняется условие $\xi_-\wedge\xi_-=\xi_+\wedge\xi_+=0$, то
формы Лапласа разложимы: $\xi_\pm=\eta_\pm\wedge\vartheta_\mp$ для
некоторых дифференциальных $1$-форм
$\eta_\pm,\vartheta_\pm\in\Omega^1(C_\pm)$. Рассмотрим следующие
$1$-мерные подраспределения распределения Картана:
$X_{\eta_\pm}=C_\pm\cap\ker\eta_\pm$ и
$X_{\vartheta_\pm}=C_\pm\cap\ker\vartheta_\pm$.

\markboth{Evgeny N.~Mikhalkin}{Contact classification of Monge-Amp\`ere equations (in Russian)}
Для уравнений общего положения эти распределения различны. Это
позволяет построить $e$-структуру для таких уравнений и найти полную
систему их скалярных дифференциальных инвариантов.

Заметим, что для уравнений Монжа-Ампера, контактно эквивалентных
уравнению, линейному относительно первых производных (т.е. уравнению
вида
\begin{equation*}
av_{xx}+2bv_{xy}+cv_{yy}+rv_{x}+sv_{y}+w=0,
\end{equation*}
где $a,b,c,r,s,w$ --- функции от $x,y,v$), посторенные $1$-мерные
распределения попарно совпадают: $X_{\eta_+}=X_{\vartheta_+}$ и
$X_{\eta_-}=X_{\vartheta_-}$.

Подробное изложение геометрии уравнений Монжа-Ампера (и не только
двумерных!) можно найти в \cite{KLR2007}.

\renewcommand{\theteorema}{\arabic{teorema}}
\renewcommand{\theequation}{\arabic{equation}}
\setcounter{teorema}{0}
\setcounter{equation}{0}
\nachalo{Е.Н.~Михалкин}{Об амебе дискриминанта алгебраического уравнения}
%markboth{Е.Н.~Михалкин}{Об амебе дискриминанта 
%лгебраического уравнения}
\label{mik-rus-abs}
Рассмотрим общее алгебраическое уравнение $n$-ой степени
 \begin{equation} \label{14.7}
  z^n+x_{n-1}z^{n-1}+\ldots+x_1z-1=0.
 \end{equation}
 Нас будет интересовать дискриминантное множество
 $\nabla=\{\Delta=0\}$ уравнения (\ref{14.7}) (здесь $\Delta$ -- дискриминант этого
 уравнения). Запишем $\nabla$, используя параметризацию
 Пассаре-Циха~\cite{PassTsikh}:
 \begin{equation} \label{14.6}
  x_k(s)=\frac{ns_k}{\langle \alpha,s \rangle}
  \left(-\frac{\langle \alpha,s \rangle}{\langle \beta,s \rangle} \right)^{\frac{k}{n}},\;\:k=1,\ldots,
  n-1;\;\; s\in \mathbb{CP}_{n-2},
  \end{equation}
 где
 \begin{equation} \label{6}
  \alpha=(n\!-\!1,\ldots,2,1), \;\;
 \beta=(1,2,\ldots,n\!-\!1)
 \end{equation}
  -- целочисленные  векторы,
   $\langle \,, \rangle$ -- знак скалярного произведения.

В случае, когда уравнение (\ref{14.7}) содержит лишь один
параметр $x_i$ (в этом случае рассматриваемое уравнение  называется
триномиальным), дискриминантное множество $\nabla$ представляет
собой некоторое подмножество точек комплексной плоскости $\mathbb{C}$ (см.
\cite{PassTsikh}). Но если уравнение (\ref{14.7}) содержит два
параметра  и более, то дискриминантное множество удобно
исследовать, используя амебу дискриминанта (определение амебы
было дано Гельфандом-Капрановым-Зелевинским
\cite{GKZ}). Отметим, что когда интересующее нас уравнение
содержит лишь два параметра $x_i$, то амеба дискриминанта $\mathcal{A}_\nabla$ лежит в $\mathbb{R}^2$.

 В статье \cite{SMZ} из интегральной формулы Меллина \cite{ProdMel} было получено интегральное представление для
 одного из
  решений уравнения (\ref{14.7}) с интегрированием по компакту. Была найдена и область сходимости полученного интеграла. А именно,  доказана
 следующая

 \begin{teorema} \label{23} Ветвь алгебраической функции $z_0(x)$ решения
уравнения (1) с условием $z(0)=1,$ допускает представление в виде
интеграла
% \begin{equation*}
%z_0(x)=1+\frac{1}{2\pi
%in}\int\limits_0^1{t^{\frac{1-n}{n}}{(1-t)}^ 
%{-\frac{1+n}{n}}\:
%\Big[e^{\frac{\pi i}{n}}\:  \mathrm{ln}\: 
%\Big(1-\sum_{k=1}^{n-1}{
%{e^{\frac{k}{n}\pi i}
 %x_kt^{\frac{k}{n}}{(1-t)}^{\frac{n-k}{n}}  }}\Big)-}$$
%$$
%-e^{-\frac{\pi i}{n}}\:  \mathrm{ln}\: \Big(1-\sum_{k=1}^{n-1}
%{{e^{-\frac{k}{n}\pi i}
 %x_kt^{\frac{k}{n}}{(1-t)}^{\frac{n-k}{n}} } \Big) \Big] } \, dt,
 %\end{equation*}
\begin{multline*}
z_0(x)=\\
1+\frac{1}{2\pi
in}\int\limits_0^1{t^{\frac{1-n}{n}}{(1-t)}^ 
{-\frac{1+n}{n}}\:
\Big[e^{\frac{\pi i}{n}}\:  \mathrm{ln}\: 
\Big(1-\sum_{k=1}^{n-1}{
{e^{\frac{k}{n}\pi i}
 x_kt^{\frac{k}{n}}{(1-t)}^{\frac{n-k}{n}}  }}\Big)-}\\
-e^{-\frac{\pi i}{n}}\:  \mathrm{ln}\: \Big(1-\sum_{k=1}^{n-1}
{{e^{-\frac{k}{n}\pi i}
 x_kt^{\frac{k}{n}}{(1-t)}^{\frac{n-k}{n}} } \Big) \Big] } \, dt,
 \end{multline*}
\markboth{Evgeny N.~Mikhalkin}{On amoeba of discriminant (in Russian)}
где
ветви  логарифма определены в  области пространства
$\mathbb{C}^{n-1}$ переменного\\ $x=(x_1, \dots ,x_{n-1}),$
полученной удалением из $\mathbb{C}^{n-1}$ двух семейств
комплексных гиперплоскостей
$$
\begin{array}{l}
\Sigma_-=\bigcup\limits_{t\in [0;1]}  \Big\{
\sum\limits_{k=1}^{n-1}{x_kt^{\frac{k}{n}}{(1-t)}^
{\frac{{n-k}}{n}}e^{-\frac{k}{n}\pi i}}=1 \Big\}, \\
\Sigma_+=\bigcup\limits_{t\in [0;1]}  \Big\{
\sum\limits_{k=1}^{n-1}{x_kt^{\frac{k}{n}}{(1-t)}^
{\frac{{n-k}}{n}}e^{\frac{k}{n}\pi i}}=1 \Big\},\\
\end{array}
$$
и выбираются условием  $\ln1=0.$
 \end{teorema}

 % Следуя  \cite{ProdMel}, ветвь с условием $z_0(0)=1$ была названа главным решением уравнения (\ref{14.7}).

 Поставим задачу исследовать взаимное расположение
 дискриминантного множества $\nabla$ уравнения (\ref{14.7}) и
 семейства
 гиперплоскостей $\Sigma_+$ ($\Sigma_-$). Решение  задачи будет
 более наглядным, если перейти к логарифмической шкале
  $$Log:  (x_1,x_2,\ldots,x_{n-1})\longrightarrow (log|x_1|,log|x_2|,\ldots,log|x_{n-1}|).$$

 Обозначим через
 $$F_{\pm}(x;t)=\sum_{k=0}^{n-1}x_kt^{\frac{k}{n}}{(1-t)}^
{\frac{{n-k}}{n}}e^{\pm\frac{k}{n}\pi i}-1 $$
 -- пару функций, линейных относительно $x$.

 Использовав параметризацию  (\ref{14.6}) дискриминантного множества
 $\nabla$ уравнения (\ref{14.7}),
  а также
   параметризацию нулевого множества функций $F_{\pm}(x;t)$
   $$x_l(\tau)=\frac{\tau_l}{a_{n-1}^{\pm}\tau_{n-1}+\ldots+a_{1}^{\pm}\tau_{1}},\; \tau_l\in \mathbb{C},$$
   при
 $$a_l^{\pm}=t^{\frac{l}{n}} (1-t)^{\frac{n-l}{n}} e^{\pm\frac{l}{n}\pi i},\;\;\; l=1,\ldots, n-1,$$
 можно показать справедливость следующего утверждения.
\markboth{A.V.~Churkin, S.N.~Sergeev}{On amoeba of discriminant (in Russian)}
 \begin{teorema} \label{8}
 Контур амебы дискриминанта уравнения (\ref{14.7}) при $s\in
 \mathbb{R}^+_{n-1}$ является огибающей для семейства амеб
 гиперплоскостей $\Sigma_{\pm}$ при $\arg \tau_l=\mp\frac{l}{n}\pi$.
 Более того, в случае $n=3$ для указанного семейства гиперплоскостей
 является огибающей контур амебы дискриминанта уравнения
 (\ref{14.7}) и при $\frac{\langle \alpha,s \rangle}{\langle \beta,s \rangle}>0.$
 Значения ${\langle \alpha,s \rangle}, \: {\langle \beta,s
 \rangle}$ находятся из равенств~(\ref{6}).
  \end{teorema}

  В дополнение к вышеизложенному, в докладе будет приведена
  геометрическая иллюстрация Теоремы \ref{8} для
  дискриминанта
  $$\Delta(x)=27+4x_1^3-4x_2^3+18x_1x_2-x_1^2x_2^2$$
  кубического уравнения
  \begin{equation} \label{7}
  z^3+x_2z^2+x_1z-1=0.
  \end{equation}
 В дополнение к этому, для уравнения (\ref{7}), в логарифмической шкале
 $$Log:  (x_1,x_2)\longrightarrow (log|x_1|,log|x_2|)$$
 будет найдено уравнение кривой, которая соответствует пересечению дискриминантного
 множества с комплексными прямыми~$\Sigma_{\pm}$. Отметим  некоторые ее свойства:
 это петля, проходящая вокруг
 каспидальной точки, симметричная относительно прямой
 $$u=v, \;\;\mbox{где} \;\; u=log|x_1|,\:v=log|x_2|.$$

\renewcommand{\theteorema}{\thesection.\arabic{teorema}}
\renewcommand{\theequation}{\thesection.\arabic{equation}}
\setcounter{section}{0}
\setcounter{footnote}{0}
\nachaloe{С.Н.~Сергеев и А.В.~Чуркин}{Программа для
демонстрации универсальных алгоритмов решения дискретного 
уравнения Беллмана в различных полукольцах}{Работа выполнена при поддержке грантов
РФФИ  05-01-00824 и 05-01-02807-НЦНИЛ\_а.}
\markboth{A.V.~Churkin, S.N.~Sergeev}{Universal algorithms (in Russian)}
%\markboth{С.Н.~Сергеев и А.В.~Чуркин}{Программа для
%демонстрации универсальных алгоритмов} 
\label{chu-rus-abs}

%Рассматриваются особенности работы программы для демонстрации ряда
%универсальных алгоритмов решения уравнения Беллмана в различных
%полукольцах.

%Обсуждается возможность использования
%универсальных объектно-ориентированных алгоритмов для
%применения в различных задачах математической физики. Работа
%выполнена при финансовой поддержке гранта РФФИ 05-01-02807-НЦНИЛ.

{\it Назначение программы.} Программа предназначена для
демонстрации некоторых универсальных алгоритмов обращения матрицы
и решения уравнения Беллмана в различных полукольцах.
В зависимости от выбора полукольца, программа может либо найти обратную матрицу и решить уравнение $Ax=B$,
где $A$ и $B$ - пользовательская матрица и вектор-столбец соответственно, либо найти матрицу $A^*$ и решить уравнения Беллмана 
$x=A\bigotimes x \bigoplus B$
Перед запуском программы пользователь выбирает одно из полуколец,
требуемую задачу и алгоритм расчета. Затем исходные данные заносятся в матрицу (для наглядности максимальный размер ограничен
величиной 10х10).  
%Непосредственно процедура расчета происходит после нажатия клавиши "calculate".
Результат расчета выводится 
либо в виде матрицы, либо в виде вектора столбца, в зависимости от поставленной задачи.
В процессе разработки программы использовался объектно -
ориентированный подход, позволяющий в полной мере использовать
универсальность предложенных алгоритмов, подключая в качестве
объектов полукольца с определенной арифметикой, актуальной для
решения конкретной задачи. 

{\it Примеры полуколец.} Использование идемпотентных операций $\bigoplus$ и $\bigotimes$ позволяет записать
ряд важных алгоритмов обращения матрицы в универсальном виде.
Выбор пользователем требуемого полукольца определяет тип данных, с которыми будет работать универсальный вычислительный алгоритм.
В программе реализована возможность выбора из следующих полуколец:

1) $\bigoplus="+"$ и $\bigotimes="\times"$ - обычная арифметика.

2) $\bigoplus="max"$ и $\bigotimes="+"$ - арифметика max-plus, в которой операция взятия максимума используется вместо сложения, а сложение - вместо умножения. Такая арифметика часто используется в задачах максимизации, системах автоматического управления и др.

3) $\bigoplus="min"$ и $\bigotimes="+"$ -  арифметика min-plus, в которой операция взятия минимума используется вместо сложения, а сложение - вместо умножения. Используется в задачах нахождения кратчайшего пути, задачах оптимизации.

4) $\bigoplus="max"$ и $\bigotimes="\times"$ - арифметика, в которой  вместо сложения используется операция взятия максимума.

5) $\bigoplus="max"$ и $\bigotimes="min"$ - арифметика max-min, в которой  вместо сложения используется операция взятия максимума,
вместо умножения -  операция взятия минимума. Используется в задачах многокритериальной оптимизации.

6) $\bigoplus="or"$ и $\bigotimes="and"$ - логическая арифметика над булевскими  переменными.

Каждому полукольцу в программе соответствует пользовательский тип данных, определяемый отдельным классом.
Кроме описаний правил сложения и умножения, внутри класса определяются также вид нуля, единицы и правило взятия операции "*".
Такой подход дает возможность дополнять программу новыми типами полуколец, не меняя структуру основной программы и не внося никаких изменений в
ту ее часть, которая занимается вычислительными алгоритмами.

{\it Универсальные алгоритмы.} В последнее время большое количество  работ 
(например, \cite{Carree}-\cite{Sergeeev}) посвящено разработке универсальных версий
алгоритмов линейной алгебры и численного анализа. 
Рассматриваемая программа использует универсальные версии ряда
классических алгоритмов обращения матриц и решения систем линейных уравнений. 
На выбор пользователю предлагаются ряд алгоритмов, как точных,
так и с использованием метода последовательных приближений:

1)Метод исключения по схеме Гаусса;

2)Метод окаймления;

3)Итерационный метод Якоби;

4)Итерационный метод  Гаусса-Зейделя;

5)Алгоритмы для матриц специального вида: симметричных, треугольных, теплицевых и др., в том числе
предложенные в \cite{Glitvinov2} и \cite{Sergeeev}.

В случае выбора обычной арифметики эти алгоритмы ведут себя классическим образом, однако для случая идемпотентного полукольца они позволяют найти матрицу $A^*$ или решить соответствующее уравнение Беллмана. Для полуколец max-plus или min-plus уравнение Беллмана представляет собой  основное функциональное уравнение динамического программирования и выражает принцип оптимальности Беллмана: управление на каждом шаге должно быть оптимальным с точки зрения процесса в целом.

{\it Возможности визуализации.}
Для наглядного представления информации в программе заложена возможность визуализации исходной матрицы в виде графа с соответствующими весами. На отдельной вкладке диалогового окна программы отображается введенная пользователем информация об исходной матрице и результат соответствующего расчета.  В таком режиме работы программы, например, задача о решении уравнения Беллмана в полукольце min-plus будет отображаться как кратчайший путь между узлами заданного пользователем графа.

{\it Использование различных арифметик для контроля точности.} Применение при разработке программы объектно ориентированного подхода позволяет
не только независимо менять полукольца и алгоритмы вычисления, но и управлять базовым типом числовых данных для контроля за точностью вычислений.
В следующей версии программы предполагается реализовать механизм выбора одной из числовых арифметик. Среди них арифметика целых чисел, арифметика чисел с плавающей точкой, дробно-рациональная арифметика с использованием цепных дробей, в том числе с контролируемой точностью. Это позволит сравнить ошибку округления, накопленную в ходе применения того или иного вычислительного алгоритма с ошибкой самого метода (для итерационных алгоритмов), что позволит судить об их итоговой эффективности.

\newpage
\setcounter{section}{0}
\setcounter{footnote}{0}
\nachalo{Роман Ульверт}{Кривые в $\mathbb{C}^2$, амебы которых
определяют фундаментальную группу дополнения}
%\markboth{Роман Ульверт}{Кривые в $\mathbb{C}^2$}
\markboth{Roman Ulvert}{Curves in $\mathbb{C}^2$ (in Russian)}
\label{ulv-rus-abs}

В классической работе ван Кампена \cite{vK}, продолжающей
исследования многих математиков, начиная от Зарисского, был
предъявлен метод вычисления фундаментальной группы дополнения к
плоской комплексной кривой. Результат ван Кампена впоследствии был
переформулирован Б. Мойшезоном и М. Тайхер с использованием
понятия брэйд-монодромии, рассмотренном в \cite{Moishe}.

Дадим краткое описание гомоморфизма 
брэйд-монодромии. Пусть
алгебраическая кривая $C$ 
задана множеством нулей полинома $f(x,y)
\in \mathbb{C}[x,y]$. Будем смотреть на 
$f$ как на полином
Вейерштрасса:
$$f=\alpha_0(y)x^d+\alpha_1(y)x^{d-1}+\dots+
\alpha_d(y).$$ 
Обозначим
$B:=\mathbb{C}\setminus \Delta$, 
где $\Delta = \{y_1,\dots,y_s\}$
--- дискриминант полинома $f$. 
Ограничение проекции $(x,y)\mapsto
y$ на множество 
$E:=\mathbb{C}\times B \setminus C$ 
определяет
локально тривиальное расслоение 
$p \colon E \to B$. Вычисление
фундаментальной группы 
$\pi_1(\mathbb{C}^2\setminus C)$ сводится к
вычислению группы $\pi_1(E)$. 
Выберем $(x_0,y_0) \in E$ и
обозначим через $F$ слой над точкой $y_0$. 
База и слой расслоения
$p$ представляют собой дополнения к конечным наборам точек в
$\mathbb{C}$, поэтому группы $\pi_1(B,y_0)$ и $\pi_1(F,x_0)$
являются свободными группами и точная последовательность
расслоения
$$
1 \longrightarrow \pi_1(F,x_0) \longrightarrow \pi_1(E,(x_0,y_0))
\longrightarrow \pi_1(B,y_0) \longrightarrow 1.
$$
расщепляется, так что группа $\pi_1(E,(x_0,y_0))$ есть полупрямое
произведение групп $\pi_1(B,y_0)$ и $\pi_1(F,x_0)$. Действие
группы $\pi_1(B,y_0)$ на $\pi_1(F,x_0)$ может быть описано в
терминах групп кос. Для этого слой $F$ отождествим с диском $D^2$,
из которого выброшено множество $K=\{x_1,\dots,x_d\}$ различных
точек. Каждый элемент группы $\pi_1(B,y_0)$ определяет биекцию
множества $K$, а следовательно и элемент группы $\Br_d = \Br(D,K)$
кос из $d$ нитей. Таким образом определен гомоморфизм
$\pi_1(B,y_0) \to \Br_d$, носящий название гомоморфизма
брэйд-монодромии. Этот гомоморфизм позволяет выписать соотношения
между образующими группы $\pi_1(E,(x_0,y_0))$.

Переход к рассмотрению амебы $\mathcal{A}_C$ кривой $C$, то есть
образа кривой под действием отображения логарифмической проекции
$(x,y)\mapsto (\ln |x|,\ln |y|)$, ставит вопрос о том, каким
образом знание амебы может помочь в вычислении фундаментальной
группы $\pi_1(\mathbb{C}^2\setminus C)$. Обозначим через $E_{\nu}$
связную компоненту $\mathbb{R}^2 \setminus \mathcal{A}_C$ порядка
$\nu=(\nu_x,\nu_y)$. Определение и свойства порядка $\nu$ можно
найти в статье Форсберга-Пассаре-Циха \cite{FPTs}. Для нас
существенно, что целые числа $\nu_x$ и $\nu_y$ выражают
коэффициенты зацепления петель
$\Gamma_x^{\nu}=\{x=e^{u+it},y=e^v\}$ и
$\Gamma_y^{\nu}=\{x=e^u,y=e^{v+it}\}$ с кривой $C$, причем $\nu_x$
и $\nu_y$ не зависят от выбора точки $(u,v)\in E_{\nu}$. Пусть
$(x_0,y_0)=(e^u,e^v)$, где $(u,v)\in E_{\nu}$. Тогда классы петель
$[\Gamma_x^{\nu}],\; [\Gamma_y^{\nu}] \in
\pi_1(\mathbb{C}^2\setminus C,(x_0,y_0))$ коммутируют. Возникает
вопрос, можно ли описать фундаментальную группу дополнения к
плоской кривой $C$, используя в качестве образующих только классы
петель, получающихся из петель $\Gamma_x^{\nu},\; \Gamma_y^{\nu}$,
где $\nu$ пробегает порядки всех связных компонент дополнения к
амебе кривой?

Чтобы ответить на этот вопрос, мы должны вернуться к описанию
группы $\pi_1(\mathbb{C}^2\setminus C)$ с использованием
брэйд-монодромии. Так как старшая степень, с которой переменная
$x$ входит в полином $f$, равна $d$, то найдется хотя бы одна
связная компонента $E_{\widetilde{\nu}}$ из $\mathbb{R}^2
\setminus \mathcal{A}_C$ порядка ${\widetilde{\nu}} = (d,\nu_y)$.
Выберем $(u,v)\in E_{\widetilde{\nu}}$ так, чтобы $v\neq \ln|y_i|$
для всех $y_i$ из дискриминанта $\Delta$ полинома $f$. Тогда для
всех $t\in [0,2\pi]$ диск $\{\ln|x|\leqslant u,\,y=e^{v+it}\}$
содержит ровно $d$ корней полинома $f(x,e^{v+it})$. Следовательно
определена коса из $d$ нитей, которую можно сопоставить петле
$\Gamma_y^{\widetilde{\nu}}=\{x=e^u,y=e^{v+it}\}$. Отсюда видно,
что петли $\Gamma_x^{\nu},\; \Gamma_y^{\nu}$ естественным образом
связаны с брэйд-монодромией кривой, то есть могут претендовать не
только на роль образующих группы $\pi_1(\mathbb{C}^2\setminus C)$,
но и давать соотношения между образующими.

В качестве простого примера кривой, фундаментальная группа
дополнения которой определяется амебой, рассмотрим дискриминант
$\Delta[2,3]$ кубического уравнения $z^3 + z^2 + xz + y = 0$. Он
задается полиномом $27y^2 + 4x^3 + 4y - 18xy - x^2$ и представляет
собой каспидальную кривую с невырожденной амебой. Фундаментальная
группа $\pi_1(\mathbb{C}^2\setminus \Delta[2,3])$ описывается с
помощью петель $\Gamma_x^{(3,0)}$ и $\Gamma_y^{(0,2)}$.
Брэйд-монодромия дает изоморфизм $\pi_1(\mathbb{C}^2\setminus
\Delta[2,3])\cong \Br_3$.

\newpage
\section*{LIST OF PARTICIPANTS AND AUTHORS}
\vskip0.3cm
\label{list-of-pars}
\thispagestyle{empty}
\avtor{AKIAN Marianne}{\INRIA}{France}{\\Marianne.Akian@inria.fr}
\avtor{BAKLOUTI Ali}{Dep\-art\-ment of Math\-ema\-tics,   Fac\-ulty of
 Sci\-en\-ces at Sfax, Rou\-te de Sou\-kra, 3038, Sfax}{Tunisia}{Ali.Baklouti@fss.rnu.tn}
\avtor{BELAVKIN Viacheslav P.}{School of Math\-ema\-ti\-cal Sci\-en\-ces,
Not\-ting\-ham Uni\-ver\-sity, NG7 2RD}{United Kingdom}{Viacheslav\_Belavkin@nottingham.ac.uk} 
\thispagestyle{empty}
\rusavtor{BENIAMINOV Evgeny M.}%
{Бениаминов Евгений Михайлович}{Russian State Uni\-ver\-sity
for the Hum\-ani\-ties, Fac\-ulty of Math\-ema\-tics,
\mbox{ul. Chayanova 15}, Moscow}%
{Russia}{ebeniamin@yandex.ru}
\rusavtor{CHOURKIN Andrey V.}{Чуркин Андрей Валерьевич}%
{\fizfak}{Russia}%
{churandr@mail.ru}
\thispagestyle{empty}
\rusavtor{DANILOV Vladimir I.}{Данилов Владимир Иванович}%
{CEMI RAS, Nahi\-mov\-ski prosp. 47,
117418, Moscow}{Russia}{\\vdanilov43@mail.ru}
\avtor{FARHI Nadir}{\INRIA}{France}{\\Nadir.Farhi@inria.fr}
\thispagestyle{empty}
\avtor{FAYE Farba}{Dep\-art\-ment of Math\-ema\-tics \& Com\-pu\-ter Sci\-ence, UCAD, Dakar}
{Senegal}{\\lamffaye@yahoo.fr}
\rusavtor{FOCK Vladimir V.}{Фок Владимир Владимирович}%
{ITEP, B.~Che\-re\-mush\-kin\-skaya~25, Mos\-cow}{Russia}{fock@itep.ru}
\avtor{GAUBERT St\'ephane}{\INRIA}{France}{\\Stephane.Gaubert@inria.fr}
\rusavtor{GEL'FAND Alexander M.}{Гельфанд Александр Маркович}{\delfin}
{Russia}{\\gelfand\_a@oaoesp.ru}
\avtor{GONCHAROV Alexander B.}%
{Brown Uni\-ver\-sity, Math. Dep\-art\-ment, \mbox{151 Thayer st.}, Pro\-vid\-ence, RI}%
{USA}{Alexander\_Goncharov@brown.edu}%
\thispagestyle{empty}
\rusavtor{GRBI\'C Tatjana}{Грбич Татjана}{\novisadg }{Serbia}%
{tatjana@uns.ns.ac.yu}
\rusavtor{GULINSKY Oleg V.}{Гулинский Олег Викторович}{IPPI RAS,
B.~Ka\-ret\-ny~1, Mos\-cow}{Russia}{bedelbaeva\_aigul@mail.ru}
\avtor{GUREVICH Dmitry}{LAMAV, Uni\-ver\-site de 
Va\-len\-cien\-nes,\\ Le Mont Houy 59313 Va\-len\-cien\-nes}{France}%
{d.gurevich@free.fr}
\thispagestyle{empty}
\rusavtor{GUTERMAN Alexander E.}{Гутерман Александр Эмильевич}%
{\mehmat}{Russia}{guterman@list.ru}
\avtor{IGONIN Sergey}{Utrecht University, NL-3508 TA Utrecht}{\\The Netherlands}{igonin@mccme.ru}
\avtor{ITENBERG Ilia}{\IRMA}{France}{itenberg@math.u-strasbg.fr}
\rusavtor{KARZANOV Alexander V.}{Карзанов Александр Викторович}%
{ISA RAS, prosp. 60-le\-tiya Ok\-tya\-brya 6, 117312, Mos\-cow}{Russia}%
{\\sasha@cs.isa.ru} 
\avtor{KHARLAMOV Vyacheslav}{\IRMA}{France}{kharlam@math.u-strasbg.fr}
\rusavtor{KIRSHTEYN Boris Kh.}{Кирштейн Борис Хаймович}{\delfin}
{Russia}{\\bkirsh@aha.ru}
\thispagestyle{empty}
\rusavtor{KOLOKOLTSOV Vassili N.}{Колокольцов Василий Никитич}%
{Moscow Institute of Eco\-no\-mics, Pe\-cher\-ska\-ya~ul.~6/1, 129344, \mbox{Moscow},
and the Uni\-ver\-sity of War\-wick, Dept. of Sta\-tis\-tics, Cov\-en\-try CV4 7AL}
{Russia and United Kingdom}{vkolok@fsmail.ru}
\thispagestyle{empty}
\rusavtor{KOSHEVOY Gleb A.}{Кошевой Глеб Алексеевич}%
{\CEMI}{Russia}{\\gleb\_koshevoy@mail.ru}
\rusavtor{KOVAL Gennady V.}{Коваль Геннадий Васильевич}%
{\fizfak}{Russia}
{kovgen@mail.ru}
\rusavtor{KUSHNER Alexey G.}{Кушнер Алексей Гурьевич}%
{Astrakhan State University, ul. Tatischeva 20a, 414056, Astrakhan}%
{Russia}{kushnera@mail.ru}
\rusavtor{KUTATELADZE Simon S.}{Кутателадзе Семен Самсонович}%
{S.L.~Sobolev Math\-ema\-ti\-cal Ins\-ti\-tute, prosp. Kop\-tyu\-ga 4, 630090, 
\mbox{Novosibirsk}}{Russia}{sskut@math.nsc.ru}
\thispagestyle{empty}
\rusavtor{LITVINOV Grigory L.}{Литвинов Григорий Лазаревич}
{\nezavisimy}{\\Russia}{glitvinov@gmail.com}
\rusavtor{MASLOV Victor P.}{Маслов Виктор Павлович}%
{\fizfak}{Russia}{v.p.maslov@mail.ru}

\avtor{MCCAFFREY David}
{7th floor, St. James buil\-ding, 79 Ox\-ford Street,
Man\-che\-ster M1 6SS}{United Kingdom}{david@mccaffrey275.fsnet.co.uk}

\avtor{MCENEANEY William M.}{Dept. of Mech\-ani\-cal and 
Aer\-osp\-ace Engi\-nee\-ring, MC 0411,
Uni\-ver\-sity of Cali\-for\-nia, San Diego
9500 Gil\-man Drive,
La Jolla, CA 92093-0411}{USA}{wmceneaney@ucsd.edu}

\rusavtor{MIKHALKIN Evgeny N.}{Михалкин 
Евгений Николаевич}{\krasped}%
{Russia}{mikhalkin@bk.ru}

\rusavtor{MOLCHANOV Vladimir F.}{Молчанов Владимир
Федорович}{\tambov}{Russia}{molchano@molchano.tstu.ru}

\avtor{NITICA Viorel}{\romanian\ and Dept. of Mathematics, West Chester University, PA}{Romania and USA}
{vnitica@wcupa.edu}

\avtor{PAP Endre}{\novisadps}{Serbia}{pap@im.ns.ac.yu}
\thispagestyle{empty}
\avtor{PERSSON Ulf}{Chalmers University of Tech\-no\-logy, 412 96 G\"{o}\-te\-borg}{Sweden}{\\ulfp@math.chalmers.se}
\avtor{QUADRAT Jean-Pierre}{\INRIA}{France}{\\Jean-Pierre.Quadrat@inria.fr}
\rusavtor{RASHKOVSKII Alexander Yu.}
{Рашковский Александр Юрьевич}{University of Sta\-van\-ger,
4036 Sta\-van\-ger}{Norway}
{\\alexander.rashkovskii@uis.no}

\rusavtor{SERGEEV Serge\u{\i} N.}{Сергеев Сергей Николаевич}
{\fizfak}{Russia}
{sergiej@gmail.com}

\rusavtor{SHPIZ Grigory B.}{Шпиз Григорий Борисович}{\nezavisimy}{\\Russia}{shpiz@theory.sinp.msu.su}

\avtor{SHUSTIN Eugenii}{School of Math\-ema\-ti\-cal Sci\-en\-ces, Tel Aviv Uni\-ver\-sity, Ra\-mat Aviv, Tel Aviv,\\ 69978}%
{Israel}{shustin@post.tau.ac.il}
\thispagestyle{empty}

\rusavtor{SHVEDOV Oleg Yu.}{Шведов Олег Юрьевич}
{\fizfak}{Russia}
{olegshv@mail.ru}
\thispagestyle{empty}

\avtor{SINGER Ivan}{\romanian}{Romania}{ivan.singer@imar.ro}

\rusavtor{SOBOLEVSKI\u{I} Andre\u{\i}}{Соболевский Андрей Николаевич}{\fizfak}{Russia}{ansobol@gmail.com}

\rusavtor{STOYANOVSKY Alexander V.}{Стояновский Александр Васильевич}{\nezavisimy}{Russia}{stoyan@mccme.ru}
\thispagestyle{empty}

\rusavtor{\v{S}TRBOJA Mirjana}{Штрбоjа Мирjана}{\novisadps}{Serbia}%
{mirjanas@im.ns.ac.yu}

\avtor{THIAM Mamadou}{\sengalese}{Senegal}{\\mathiam@netcourier.com}
\thispagestyle{empty}
\avtor{TRUFFET Laurent}{Institut de Recherche en Com\-mu\-ni\-ca\-tions et Cy\-ber\-ne\-tique de Nantes,
IRCCyN UMR-CNRS 6597, 1, rue de la Noe BP 92101 44321 Nantes Cedex}{France}
{Laurent.Truffet@emn.fr}

\rusavtor{TSYKINA Svetlana V.}{Цыкина Светлана Викторовна}
{\tambov}{Russia}%
{tsykinasv@yandex.ru}

\rusavtor{ULVERT Roman V.}{Ульверт Роман Викторович}{\krasfed}{Russia}{ulvertrom@yandex.ru}
\thispagestyle{empty}

\avtor{WAGNEUR Edouard}{\'Ecole Polytechnique de Mon\-tr\'e\-al, C.P.~6079 Succ. Cen\-tre-ville,
Mon\-tr\'eal Qu\'e\-bec H3C 3A7, and
GERAD, Montreal}{Canada}{Edouard.Wagneur@gerad.ca}
\thispagestyle{empty}

\avtor{WALSH Cormac}{\INRIA}{France}
{\\cormac.walsh@inria.fr}
\thispagestyle{empty}

\rusavtor{ZNAMENSKAYA Oksana V.}%
{Знаменская Оксана Витальевна}%
{\krasfed}%
{Russia}{znamensk@lan.krasu.ru}

\thispagestyle{empty}
\end{document}